\documentclass[10pt,twoside]{article}
\usepackage[T1]{fontenc}
\usepackage[utf8]{inputenc}
\usepackage[british]{babel}
\usepackage[tbtags]{amsmath}
\usepackage{amsthm,amssymb}
\usepackage[helvratio=0.88,tighter]{newtxtext}
\usepackage[varvw]{newtxmath}
\usepackage{mathtools,thmtools}
\usepackage{tikz-cd,graphicx,adjustbox}
\usepackage{imakeidx,comment}
\usepackage[shortlabels]{enumitem}
\usepackage[bottom]{footmisc}
\usepackage{titlesec}
\usepackage[hyphens]{url}
\usepackage[unicode,psdextra,hidelinks]{hyperref}
\usepackage{upref}
\usepackage[nameinlink]{cleveref}

\makeatletter
\renewcommand{\normalsize}{%
  \@setfontsize\normalsize{10}{13}%
  \abovedisplayskip 10pt plus 2pt minus 5pt
  \belowdisplayskip\abovedisplayskip
  \abovedisplayshortskip 0pt plus 3pt
  \belowdisplayshortskip 6pt plus 3pt minus 3pt
  \let\@listi\@listI
}
\renewcommand{\small}{%
  \@setfontsize\small{9}{11.5}%
  \abovedisplayskip 8.5pt plus 3pt minus 4pt
  \belowdisplayskip\abovedisplayskip
  \abovedisplayshortskip 0pt plus 2pt
  \belowdisplayshortskip 4pt plus 2pt minus 2pt
}
\renewcommand{\footnotesize}{%
  \@setfontsize\footnotesize{8.4}{10}%
  \abovedisplayskip 6pt plus 2pt minus 4pt
  \belowdisplayskip\abovedisplayskip
  \abovedisplayshortskip 0pt plus 1pt
  \belowdisplayshortskip 3pt plus 1pt minus 2pt
}
\makeatother
\normalsize
\usepackage[
  paperwidth=170mm,paperheight=240mm,
  textwidth=357pt,textheight=542.333pt,
  inner=20mm,top=13mm,includehead,
  headheight=12pt,headsep=16.45pt,footskip=25pt,
  marginparwidth=16mm
]{geometry}
\titleformat{\section}{\normalfont\large\bfseries\boldmath}{\thesection.}{0.5em}{}
\titleformat{\subsection}{\normalfont\bfseries\boldmath}{\thesubsection.}{0.5em}{}
\titlespacing*{\section}{0pt}{26pt}{13pt}
\titlespacing*{\subsection}{0pt}{13pt}{6.5pt}

\newtheoremstyle{plain}{6.5pt}{6.5pt}{\itshape}{}{\bfseries}{.}{0.5em}{}
\newtheoremstyle{definition}{6.5pt}{6.5pt}{\normalfont}{}{\bfseries}{.}{0.5em}{}

\setlist{itemsep=0.25\baselineskip,leftmargin=\parindent,parsep=0pt,
  partopsep=0.5\baselineskip,topsep=0.25\baselineskip,listparindent=\parindent}
\setlist[enumerate]{left=0pt,align=left}
\setlist[enumerate,1]{align=left,labelindent=\parindent,leftmargin=*}
\setlist[itemize]{align=left,labelindent=0pt,leftmargin=\parindent,labelsep=*}
\setlist[description]{font=\mdseries}

\let\unnumberedSubsection\subsection
\RenewDocumentCommand{\subsection}{s o m}{%
  \IfBooleanTF{#1}{%
    \unnumberedSubsection*{#3}%
    \addcontentsline{toc}{subsection}{#3}%
  }{%
    \IfNoValueTF{#2}{\unnumberedSubsection{#3}}{\unnumberedSubsection[#2]{#3}}%
  }%
}

\hypersetup{
  pdftitle={Perfect-Prismatic F-Crystals and p-adic Shtukas in Families},
  pdfauthor={Anton G\"uthge},
  bookmarksnumbered=true
}

\numberwithin{equation}{section}

\tikzset{cong/.style={draw=none,edge node={node [sloped, allow upside down, auto=false]{$\cong$}}},
	iso/.style={draw=none,every to/.append style={edge node={node [sloped, allow upside down, auto=false]{$\cong$}}}}}

\newcommand\isom{
	\xrightarrow{\smash{\raisebox{-0.65ex}{\ensuremath{\scriptstyle\sim}}}}
}

\theoremstyle{plain}
\newtheorem{mythm}{Theorem}[section]

\newtheorem{myprop}[mythm]{Proposition}
\newtheorem{mylemma}[mythm]{Lemma}
\newtheorem{mycor}[mythm]{Corollary}

\theoremstyle{definition}
\newtheorem{mydef}[mythm]{Definition}

\newtheorem{myrem}[mythm]{Remark}
\newtheorem{mynot}[mythm]{Notation}

\newtheorem{mydefprop}[mythm]{Definition + Proposition}

\newtheorem{myconst}[mythm]{Construction}

\crefname{mythm}{Theorem}{Theorems}
\crefname{myprop}{Proposition}{Propositions}
\crefname{mylemma}{Lemma}{Lemmas}

\renewcommand{\theHmythm}{\theHsection.\arabic{mythm}}

\ExplSyntaxOn
\clist_map_inline:nn
  {myprop,mylemma,mycor,myconj,mydef,mytempdef,myrem,mynot,myremdef,mypropdef,mydefprop,mywarn,myconv,myconst,myex}
  {\cs_set:cpn {theH#1} {\theHmythm}}
\ExplSyntaxOff

\DeclareMathOperator{\id}{id}
\DeclareMathOperator{\im}{Im}
\DeclareMathOperator{\Spec}{Spec}

\DeclareMathOperator{\Spa}{Spa}
\DeclareMathOperator{\Spf}{Spf}
\DeclareMathOperator{\Spd}{Spd}
\DeclareMathOperator{\GL}{GL}

\DeclareMathOperator{\Hom}{Hom}

\DeclareMathOperator{\Aut}{Aut}
\DeclareMathOperator{\End}{End}

\DeclareMathOperator{\Perf}{Perf}

\DeclareMathOperator{\perfd}{perfd}

\DeclareMathOperator*{\colim}{colim}
\DeclareMathOperator{\Rep}{Rep}
\DeclareMathOperator{\Frac}{Frac}

\DeclareMathOperator{\Gal}{Gal}

\DeclareMathOperator{\ana}{ana}
\DeclareMathOperator{\disc}{disc}

\newcommand{\ZZ}{\mathbb{Z}}
\newcommand{\NN}{\mathbb{N}}
\newcommand{\RR}{\mathbb{R}}

\newcommand{\QQ}{\mathbb{Q}}
\newcommand{\FF}{\mathbb{F}}

\newcommand{\MCC}{\mathcal{C}}

\newcommand{\MCE}{\mathcal{E}}
\newcommand{\MCF}{\mathcal{F}}
\newcommand{\MCG}{\mathcal{G}}

\newcommand{\MCO}{\mathcal{O}}

\newcommand{\MCS}{\mathcal{S}}

\newcommand{\MCX}{\mathcal{X}}
\newcommand{\MCY}{\mathcal{Y}}

\newcommand{\mfm}{\mathfrak{m}}

\newcommand{\Cech}{\v{C}ech }

\let \phi \varphi
\let \diamond \Diamond

\usepackage{relsize}
\usepackage[bbgreekl]{mathbbol}
\usepackage{amsfonts}
\DeclareSymbolFontAlphabet{\mathbb}{AMSb} 
\DeclareSymbolFontAlphabet{\mathbbl}{bbold}

\DeclareMathOperator{\perf}{perf}
\DeclareMathOperator{\arc}{arc}

\DeclareMathOperator{\eq}{eq}
\DeclareMathOperator{\fin}{fin}
\DeclareMathOperator{\pol}{pol}

\newcommand{\Prism}{{\mathlarger{\mathbbl{\Delta}}}}

\DeclareMathOperator{\Inf}{inf}

\DeclareMathOperator{\Sh}{Sh}

\DeclareMathOperator{\supp}{supp}
\DeclareMathOperator{\Flat}{flat}

\newcommand{\Ainf}{A_{\Inf}}

\newcommand{\pRoots}{{1/p^\infty}}

\newcommand{\ZZip}{\ZZ[p^{-1}]}
\newcommand{\ZZgez}{\ZZ_{\geq 0}}
\newcommand{\ZZgz}{\ZZ_{> 0}}
\renewcommand{\Langle}{\Big\langle}
\renewcommand{\Rangle}{\Big\rangle}

\newcommand{\Ycc}{\MCY_{[0,\infty]}}
\newcommand{\Yco}{\MCY_{[0,\infty)}}
\newcommand{\Yoc}{\MCY_{(0,\infty]}}
\newcommand{\Yoo}{\MCY_{(0,\infty)}}
\newcommand{\OpP}{\MCO^{\perf}_{\Prism}}
\newcommand{\rrrarrows}{\;\substack{\rightarrow\\[-1em] \rightarrow \\[-1em] \rightarrow}\;}

\ExplSyntaxOn
\clist_map_inline:nn
  {myprop,mylemma,mycor,myconj,mydef,mytempdef,myrem,mynot,myremdef,mypropdef,mydefprop,mywarn,myconv,myconst,myex}
  {\AddToHook{env/#1/begin}{\crefalias{mythm}{#1}}}
\ExplSyntaxOff

\title{Perfect-Prismatic F-Crystals and p-adic Shtukas in Families}
\author{Anton G\"uthge\thanks{Department of Mathematics, Technische Universit\"at Darmstadt, Schlossgartenstra\ss e 7, 64289 Darmstadt, Germany. Email: \href{mailto:anton.guethge.math@gmail.com}{\nolinkurl{anton.guethge.math@gmail.com}}.}}
\date{}

\begin{document}
\maketitle

\begin{abstract}
We construct an equivalence between the two categories in the title, thus establishing a link between Frobenius-linear objects of formal (schematic) and analytic (adic) nature. We will do this for arbitrary \(I\)-complete rings with \(I\) finitely generated and containing \(p\), for objects equipped with an action by an arbitrary affine flat group scheme of finite type over \(\ZZ_p\) and without making use of the Frobenius-linear structure.
\end{abstract}

\begingroup
\par
\endgroup

\tableofcontents

\section{Introduction}
A slightly special case of our main theorem is the following:
\begin{mythm}\emph{(\Cref{BKFEquivalence})}
    Let \(R^+\) be an integral perfectoid ring. Then pullback defines an exact tensor equivalence
    \[\emph{BKF-modules on \(R^+\)}\isom \emph{families of shtukas on \(\Spd(R^+)\),}\]
    where by \emph{families of shtukas} we mean the following: For every map from an affinoid perfectoid \(\Spa(B,B^+)\to \Spd(R^+)\), we specify a shtuka on \(\Spd(B,B^+)\) such that the untilt \((B^\sharp, B^{\sharp+})\) defining the leg of the shtuka receives a map \((R^+,R^+)\to (B^\sharp, B^{\sharp+})\).
\end{mythm}
Thus, the theorem provides a comparison between Frobenius-linear data of schematic and adic nature. The fact that the first category naturally lives in mixed characteristic (i.e., \(\Spec(R^+)\)) and the second one in characteristic \(p\) (i.e., \(\Spd(R^+)\)) requires us to fix the leg in the second one. This change in characteristic is not the main concern, however -- the main point is the transition from schematic to adic.\par
The most general versions of our main result are stated in \Cref{equivalenceCorIntroduction} and \Cref{FEquivalenceThmIntroduction}. To formulate them, we will need to introduce some new terminology, but first, we want to give some perspective on the result.

\subsection*{v-locally trivial implies Zariski-locally trivial: The general perspective}

For a fixed scheme \(F\), an \(F\)-bundle on a scheme \(B\) is a scheme map \(E\to B\) that locally on \(B\) looks like the projection \(B\times F\to B\). A common question is then whether the category of \(F\)-bundles stays the same when changing what ``locally'' means in this context, i.e., when passing to a coarser or finer topology. For example, vector bundles of rank \(n\) are \(\mathbb{A}^n\)-bundles (with linear transition maps), and it does not matter if we ask for them to be trivializable in the Zariski or fpqc topology: The two categories are equivalent. Before going into the precise definition of our objects of study, let us examine our main theorem from this point of view.\par
Our base scheme \(B\) will (for now) be \(\Spec(R)\) for a perfect ring \(R\) of characteristic \(p\), and we want to consider objects that locally look like \(W(R)^{n}\), where \(W(R)\) is the ring of \(p\)-typical Witt vectors over \(R\). The category on the side of the equivalence using the coarser topology will consist of objects that are trivializable Zariski-locally on \(\Spec(R)\). We give two equivalent definitions. The first is most in line with the philosophy put forth in this introduction, while the second is the one that is actually used.

\begin{myrem}
\((W(R),(p))\) being a henselian pair implies that the following pieces of data are equivalent:
    \begin{enumerate}
        \item A sheaf of \(W(\MCO_{\Spec(R)})\)-modules on \(\Spec(R)\) that Zariski-locally\footnote{By \cite{Takaya} (which will also be used later in the text), one could replace ``Zariski-local'' by ``arc-local'' here. This is another instance where the notion of ``Witt vector bundles'' stays the same when passing to a different topology, but it differs from our result in that it stays in the world of schemes while our ``finer topology'' lives in the world of adic spaces.} looks like a finite direct sum of copies of \(W(\MCO_{\Spec(R)})\). 
        \item A vector bundle on \(\Spec(W(R))\).
    \end{enumerate}
\end{myrem}
    
Adic spaces can be thought of as a common generalization of schemes, formal schemes and rigid analytic spaces, where the \(\Spec\) and \(\Spf\) functors are replaced by a functor called \(\Spa\), taking as input, e.g., complete topological rings carrying the \(p\)-adic topology (in general, the input will be certain pairs \((R,R^+)\) of topological rings). For precise definitions, see \cite[Sections 1-3]{SW20}. Some useful facts about adic spaces are:
\begin{itemize}
    \item They have a nice subcategory consisting of so-called \emph{affinoid perfectoids} (\cite[Definition 7.1.2 and Remark 7.1.3]{SW20}), which carry a (very fine) topology called the v-topology. One can then define a functor from adic spaces to \emph{v-sheaves on affinoid perfectoids in characteristic \(p\)} (\cite[Section 18.1]{SW20}), or just \emph{v-sheaves} for short, written as \(X\mapsto X^\diamond\).
    \item The category of adic spaces\footnote{We always assume adic spaces to be \(p\)-complete in the sense that for all Huber pairs \(p\) is contained in some ideal of definition.} has a final object \(\Spa(\ZZ_p)\), and the functor from adic spaces to v-sheaves over \(\Spa(\ZZ_p)^\diamond\) is fully faithful for certain nice adic spaces (\cite[Proposition 10.2.3, Proposition 18.3.1 and Proposition 18.4.1]{SW20}).
    \item For \(R\) a perfect ring of characteristic \(p\), interpreted as a topological ring with the discrete topology, the definition of \(\Spd\) and the universal property of Witt vectors imply that
    \[\Spd(W(R))\cong \Spd(R)\times \Spd(\ZZ_p).\]
    This is the way in which the Witt vectors can be understood as a literal product of v-sheaves referenced above.
\end{itemize}

The finer topology on the other side of the equivalence in our main theorem is the following: We view \(R\) as a v-sheaf \(\Spa(R)^\diamond\), and look at the site of v-sheaves over \(\Spa(R)^\diamond\) with the topology generated by surjective maps. Now for a v-sheaf \(\MCX\to \Spa(R)^\diamond\), the v-sheaf \(\MCX\times\Spa(\ZZ_p)^\diamond\) (which we view as ``\(W(\MCX)\)" through the last bullet point above) is not always representable by an adic space, even if \(\MCX\) is. However, if \(\MCX\) is a representable v-sheaf, i.e. \(\MCX = X^\diamond\) for \(X\) affinoid perfectoid, \(X^\diamond\times\Spa(\ZZ_p)^\diamond\) is represented by an adic space usually called \(\Yco(X)\) (\cite[Proposition 11.2.1]{SW20}). As in any category of sheaves on a site (satisfying appropriate finiteness conditions), representable sheaves form a basis of the topology. We can now define the category on the other side of the equivalence:
\begin{mydef}
    An integral \(\MCY\)-bundle on a quasi-compact v-sheaf \(\MCX\) is a cover \(X\to \MCX\) by an affinoid perfectoid \(X\) together with a vector bundle \(\MCE\) on \(\Yco(X)\) and a (hyper-)descent datum of \(\MCE\) over \(\MCX\), up to the usual notion of equivalence.
\end{mydef}

\begin{myrem}
    Note that \(\Yco(X)\) is not quasi-compact and that it is not at all clear whether every vector bundle on \(\Yco(X)\) is already trivializable on \(X\). Even for \(X=\Spa(K,K^+)\) a perfectoid field, I am uncertain if all vector bundles on \(\Yco(\Spa(K,K^+))\) are trivial, though our main result implies this a fortiori for those vector bundles on \(\Yco(K,K^+)\) which arise as part of the datum of an integral \(\MCY\)-bundle on \(\Spa(R^+)\) for a \(p\)-complete ring \(R^+\).
\end{myrem}

Now a special version of our main theorem reads as follows: 

\begin{mythm}
    For \(R\) perfect of characteristic \(p\), the natural functor from vector bundles on \(\Spec(W(R))\) to integral \(\MCY\)-bundles on \(\Spa(R)^\diamond\) is an equivalence of categories.
\end{mythm}

\begin{myrem}
    Another perspective on this suggested by the referee is that of seeing the main result as a version of GAGA, linking the (\(p\)-adic) analytic and algebraic worlds. This viewpoint would be in line with \cite[8.7.7]{KL15}, which shows equivalence of vector bundles on the schematic and adic Fargues--Fontaine curve, or \cite[2.7.7]{KL15}, which shows that for a sheafy analytic Huber pair \((A,A^+)\) the vector bundles on \(\Spec(A)\) and \(\Spa(A,A^+)\) are equivalent.
\end{myrem}

\subsection*{Shtukas and Perfect-Prismatic F-Crystals}

Let's leave this standpoint of v-locally trivial implies Zariski-locally trivial and instead motivate the result starting from the categories in the title of this paper. From this point on, I assume familiarity with the language of adic spaces and v-sheaves as in \cite{SW20} as well as integral perfectoids and \(\Ainf\) as defined in \cite[Section 3.1-3.2]{BMS19}.\par
\(p\)-adic shtukas were first defined in \cite[11.4.1]{SW20} over perfectoid spaces:
\begin{mydef}
Let \((B,B^+)\) be a perfectoid Huber pair in characteristic \(p\). A \(\GL_n\)\emph{-shtuka} over \((B,B^+)\) with leg in some untilt \((B^\sharp, B^{+\sharp})\) consists of:
\begin{itemize}
	\item A vector bundle \(\MCE\) of rank \(n\) on \(\Yco(B,B^+)=\Spa(B,B^+)\dot\times \Spa(\ZZ_p)\).
	\item An isomorphism 
    \[\Phi\colon\phi^*\MCE|_{\Yco(B,B^+)-V(\xi)}\longrightarrow \MCE|_{\Yco(B,B^+)-V(\xi)},\]
    where \(\phi\colon \Yco\to\Yco\) is the Frobenius isomorphism induced from the Witt vector Frobenius \(W(B^+)\to W(B^+)\) and \(\xi\in \Ainf(B^{+\sharp})=W(B^+)\) defines the divisor cutting out \(\Spa(B^\sharp,B^{+\sharp})\) from \(\Yco(B,B^+)\). Furthermore, we ask \(\Phi\) to be meromorphic around \(\xi\) (see loc. cit. for details).
\end{itemize}
\end{mydef}
The notion of shtukas as defined above has a globalization from perfectoid spaces to arbitrary (not necessarily analytic) adic spaces \(X\), arriving at so-called \emph{shtukas in families on \(X\)}, either by v-descent from a cover \(\Spa(B,B^+)\to X\) by a perfectoid space (if \(X\) is quasi-compact) or, equivalently, as compatible families of shtukas on \((B,B^+)\) for every map \(\Spa(B,B^+)\to X\). This is done, e.g., in \cite{PR}, and it is the incarnation of shtukas that we will consider in the present article.

On the other side of the equivalence, perfect-prismatic\footnote{A hyphen is used to make clear that these are to be understood as \(F\)-crystals on the perfect-prismatic site -- not \(F\)-crystals on the prismatic side which happen to be perfect in some sense.} \(F\)-crystals are a variant of the prismatic \(F\)-crystals from, e.g., \cite[Definition 4.1]{BS21b}: Instead of the (absolute) prismatic site over some \(I\)-complete ring \(A^+\), \(p\in I\), one uses the (absolute) perfect-prismatic site over \(A^+\), which admits an equivalent definition as the site of all integral perfectoid \(A^+\)-algebras. If \(A^+\) is integral perfectoid itself, both the prismatic and perfect-prismatic sites over \(A^+\) have \((\Ainf(A^+),(d))\) (where \(d\in \Ainf(A^+)\) is a distinguished element, i.e., \(\Ainf(A^+)/(d)=A^+\)) as a final object. In this case, both the category of perfect-prismatic \(F\)-crystals and the category of usual prismatic \(F\)-crystals are equivalent to the category of Breuil-Kisin-Fargues modules, an integral mixed-characteristic version of \(F\)-isocrystals:
\begin{mydef}
    Let \(A^+\) be an integral perfectoid ring. A \emph{BKF-module} on \(A^+\) is given by:
    \begin{itemize}
        \item A finite projective \(\Ainf(A^+)=W(A^{+\flat})\)-module \(\MCE\).
        \item An isomorphism 
        \[\Phi\colon\phi^*\MCE[1/d]\longrightarrow\MCE[1/d],\]
        where \(\phi\colon \Ainf(A^+)\to \Ainf(A^+)\) is the Witt vector Frobenius and \(d\in \Ainf(A^+)\) is a distinguished element such that \(\Ainf(A^+)/d=A^+\).
    \end{itemize}
\end{mydef}
So in particular, prismatic \(F\)-crystals and perfect-prismatic \(F\)-crystals agree on integral perfectoids. Where they differ is in the fact that the category of perfect-prismatic \(F\)-crystals on any \(I\)-complete ring \(A^+\) is already completely determined by what happens on integral perfectoid \(A^+\)-algebras. Thus, all the prismatic terminology introduced in \cite{BS21} is not strictly necessary to define perfect-prismatic \(F\)-crystals: they can be completely described using the slightly more classical notions of integral perfectoid rings and BKF-modules.

The strong similarity between the two definitions above should be apparent even to someone reading them for the first time. And indeed, a number of equivalences between various Frobenius-linear objects of adic and schematic nature have been proved over the years, for example:
\begin{enumerate}
   \item For algebraically closed fields \(k\), a bijection on isomorphism classes between isocrystals and vector bundles on the Fargues-Fontaine curve in \cite[Théorème 8.2.10 (4)]{FF18} (combine loc. cit. with the Dieudonné-Manin classification of isocrystals), the latter of which has been shown to have a description by adic spaces later.
   \item Building on that, an equivalence between shtukas on \((C,C^\circ)\) and BKF-modules on \(C^{\circ}\) in \cite[Theorem 14.1.1, 1. and 5.]{SW20} (attributed to Fargues) for \(C\) an algebraically closed non-Archimedean field.
\end{enumerate}
On the other hand, shtukas and BKF-modules are objects living over fairly different kinds of spaces, as perfectoid spaces are \emph{analytic} and integral perfectoids are \emph{formal} in nature. This is also visible in the two equivalences above as these pass from objects on a (formal) scheme to objects on some version of analytic locus of the associated adic space. This passing to a smaller space gives rise to numerous nuisances:
\begin{itemize}
\item Restrictions on the input space, e.g., that \(C\) be an algebraically closed field or that \(C^+=C^\circ\),
\item Making it necessary to use the Frobenius structure in an essential way, i.e., one does not get an equivalence of the categories of vector bundles without such Frobenius structure,
\item The equivalences are usually not exact. As an example of the technical problems caused by this, it was proved in \cite{Ked19} that vector bundles on \(\Spec(\Ainf)\) and its analytic locus \(\Spec(\Ainf)-V(\varpi, p)\) (see below for notation) are equivalent. If this equivalence were exact, it would immediately extend to an equivalence of \(G\)-bundles for any affine flat group scheme \(G\) of finite type via the Tannakian formalism. The fact that it is not is the main obstacle in \cite{Ans18}, where the equivalence was proved for parahoric group schemes through careful technical considerations.
\end{itemize}

A way to get around this problem is to carry a formal scheme to the analytic world not by passing to its analytic locus, but by considering its associated v-sheaf on the category of perfectoid spaces in characteristic \(p\). 
\begin{enumerate}[resume]
\item In \cite{Ans22}, this technique is used to establish an equivalence between isocrystals on some algebraically closed field \(k\) in characteristic \(p\) and vector bundles on a ``family of Fargues-Fontaine curves" over \(\Spa(k)\).
\item In \cite{PR}, Pappas-Rapoport generalize this idea both to arbitrary schemes in characteristic \(p\) \cite[Theorem 2.3.8]{PR} and perfectoid fields in mixed characteristic \cite[Proposition 2.3.11]{PR}. Their approach using Sen theory has strongly inspired ours.
\item In \cite[Section 3 and Theorem 7.14]{GI23}, Gleason, Ivanov and Zillinger obtain a classicality result for rational bundles (i.e., excluding \(p=0\)) and, for shtukas with a leg in \(p=0\), also for integral ones.\footnote{This differs from our ``rational" version of the equivalence in \Cref{analyticLocuss} in that they invert \(p\) in a meromorphic way, instead of just removing its vanishing locus from \(\Spa(\Ainf)\) which results in \(p\) being inverted in an analytic way.}
\end{enumerate}
We want to generalize (4) while eliminating the three nuisances mentioned above: Our equivalence will work for arbitrary integral perfectoids and even (trivially by descent) for arbitrary \(I\)-complete rings with \(I\) finitely generated and \(p\in I\). It produces an exact tensor equivalence, so, in particular, we get the relevant statement for \(\mathcal{G}\)-bundles for any affine flat group scheme of finite type for free and without making use of the Frobenius structure at all:
\begin{mythm}\emph{(\Cref{GEquivalenceTheorem})}\label{equivalenceCorIntroduction}
Let \(R^+\) be a complete topological ring carrying the \(I\)-adic topology for some finitely generated ideal \(I\subset R^+\) containing \(p\), and let \(\MCG/\ZZ_p\) be any affine flat group scheme of finite type. There is a canonical equivalence of groupoids
\[\text{\emph{\(\MCG\)-\(\OpP\)-bundles on }}R^+
\longrightarrow\text{\emph{integral \(\MCG\)-\(\MCY\)-bundles on }}\Spd(R^+).\]
\end{mythm}
The \(\OpP\)-bundles can be seen as ``perfect-prismatic \(F\)-crystals without Frobenius structure", while integral \(\MCY\)-bundles are ``shtukas in families without Frobenius structure". Adding these Frobenius structures, we arrive at:
\begin{mythm}\emph{(\Cref{newFGEquivalenceThm})}\label{FEquivalenceThmIntroduction}
The functor defined above induces an equivalence of groupoids
\[\emph{perfect-prismatic \(\MCG\)-crystals on }R^+
\longrightarrow\emph{ \(\MCG\)-shtukas with fixed leg on }\Spd(R^+).\]
\end{mythm}
See \Cref{DefShtukas} for a precise definition of shtukas with fixed leg. Roughly, it means that the leg of the shtuka over \(\Spd(R^+)\) should lie in \(\Spa(R^+)\subset \Spa(\Ainf(R^+))\).

\subsection*{Proof Strategy}
For the rest of this introduction we will only talk about the linear objects with their respective Frobenius structure (so BKF-modules and families of shtukas instead of \(\OpP\)-bundles and integral \(\MCY\)-bundles) as their terminology should be more familiar to most readers, but all statements hold in the generality of \Cref{equivalenceCorIntroduction}.\par
While both \cite{Ans22} and \cite[Theorem 2.3.8]{PR} only cover the (discrete) characteristic \(p\) case, \cite[Proposition 2.3.11]{PR} deals with a ring of integral elements in a Huber pair, i.e., where ``most of the adic space is analytic". Both approaches will be relevant to us as our main result generalizes both special cases. The key to bringing them together is the following:

\begin{myprop}\emph{(\Cref{lemmaCoveringg})}
Let \(S^+\) be any \(I\)-complete ring with \(I\) finitely generated and containing \(p\). Then there exists a map \(S^+\to A^+\) which is simultaneously an \(I\)-complete arc cover and a v-cover of the associated v-sheaves with \(A^+=A^+_{\disc}\times A^+_{\ana}\) where
\begin{itemize}
    \item \(A^+_{\disc}\) is a perfect ring of characteristic \(p\).
    \item \(A^+_{\ana}\) is a ring of integral elements in an affinoid perfectoid -- in fact, of a very simple kind, a so-called \emph{product of points}.
\end{itemize} 

\end{myprop}

As perfect-prismatic \(F\)-crystals descend along \(I\)-complete arc covers (\Cref{arcDescentForFCrys}) and shtukas descend along v-covers of v-sheaves (\Cref{vDescentForShtukas}), this allows us to reduce our statement to rings of one of those basic forms.\par
For any \(A^+\), the easy direction is to produce from a perfect-prismatic \(F\)-crystal (or, equivalently, a BKF-module if \(A^+\) is perfectoid) a family of shtukas, indexed by maps to perfectoid Huber pairs \((A^+,A^+)\to (B,B^+)\). Conceptually\footnote{This does not strictly work as it is not clear whether \(\Spa(\Ainf(A^+))\) exists as an adic space in general, see \Cref{liuAinf}}, this is done simply by pullback along
\[\Yco(B,B^+)\longrightarrow \Spa(\Ainf(B^+))\longrightarrow \text{``\(\Spa(\Ainf(A^+))\)"},\]
see \Cref{equivalenceConst} for a thorough definition. The challenge then is to use the data of all these shtukas as \((B,B^+)\) varies to reconstruct an object living on \(A^+\).\par
The case \(A^+=A^+_{\disc}\) builds crucially on \cite[Corollary 5.3]{Kim}. In this case the one shtuka out of the family appearing most prominently in the proof is the one living on the perfectoid Huber pair \((B,B^+)=(A^+((t^{1/p^\infty})),A^+[[t^{1/p^\infty}]])\), which covers \(\Spd(A^+,A^+)\). Loc. cit. gives us descent for vector bundles along \(\Spa(B,B^+)\to \Spa(A^+,A^+)\), a result which we extend in \Cref{KimExtTheorem} to descent of certain Witt vector bundles by an induction and limit argument. Now the families of shtukas sit somewhere between the BKF-modules (in this case, as \(p=0\) in \(A^+\), these are just classical Frobenius crystals) on \(A^+\) and those Witt vector bundles, the map to which is shown via explicit calculations to be fully faithful, thus showing the result.\par
On the other hand, for \(A^+=A^+_{\ana}\), the most prominent shtuka out of the family is the one on \((A,A^+)\), essentially a vector bundle on the space \(\Yco(A,A^+)\), which, heuristically, feels like it should be dense in \(\Spa(\Ainf(A^+))\). Now all we need to do to construct our inverse functor is to extend \(\MCE|_{\Yco(A,A^+)}\) back to \(\Spa(\Ainf(A^+))\). We proceed in two steps: First, we extend our shtuka to \(\Ycc(A,A^+)=\Spa(\Ainf(A^+)) - V(p,[\varpi])\). Usually this is done using the Frobenius structure and some extra assumptions on \((A,A^+)\), but instead we will get this by making use of all the other evaluations of the family of shtukas (i.e., on \(\Spa(B,B^+)\) for any map to a perfectoid Huber pair \((A^+,A^+)\to (B,B^+)\)) and methods of Sen, strongly inspired by \cite[Theorem 2.3.8]{PR}. In the second step, using that \(A^+\) is a product of points, we get the extension from \(\Ycc\) to all of \(\Spa(\Ainf(A^+))\) (or, equivalently, to \(\Spec(\Ainf(A^+))\)) by a result of Kedlaya, generalized by Gleason.\par
In this sense, even though for both \(A^+_{\disc}\) and \(A^+_{\ana}\) we need the full descent datum of shtukas, one could argue that for \(A^+_{\disc}\) the proof uses more of an approach \emph{from above} (descending the vector bundle from a cover) while the proof for \(A^+_{\ana}\) uses an approach \emph{from below} (extending the vector bundle from a dense subspace).

\subsection*{Structure of the paper}

The first two sections after the introduction serve to establish the basics for both the schematic and adic sides. All proofs are straightforward, but we introduce three new definitions: \(\MCO_{\Prism}^{\perf}\)-bundles in \Cref{equivalentDefsOfPrismaticBundles} and \(\MCY\)-bundles\footnote{Related categories of vector bundles are studied in \cite[Section 3]{GI23}.} in \Cref{defYBundle}, which encode the underlying vector bundle data (i.e., without the Frobenius structure) of perfect-prismatic \(F\)-crystals (the third new definition) and shtukas in families, respectively. Notably, the \(\MCY\)-bundles sit somewhere between the well-established notions of vector bundles and v-bundles on an adic space, but the exact relationship between those three notions remains subtle. See \Cref{baseChangeDiagrams} for an overview of the situation.\par
All the interesting proofs are in Section 4. In the first subsection we extend \Cref{Kim}, a result from \cite{Kim}, to the Witt vector bundles from Subsection 3.4. The second subsection then formulates the main theorem and proves it for the \(A^+_{\disc}\) case described above (\Cref{mainThmDiscrete}). Turning to the case of \(A^+_{\ana}\), Subsections 4.3--4.5 tackle the extension from \(\Yco\) to \(\Ycc\) and the Sen theory mentioned above: Subsection 4.3 contains \Cref{proposalGeometricSemilinear} and \Cref{proposalLocalTriviality}, describing v-bundles via continuous semilinear actions and giving an effectivity criterion in terms of Galois cocycles. This result will be combined with methods of Sen in Subsection 4.4 to prove the general descent result \Cref{SenDescent} giving a criterion for when effectivity of a certain pro-étale descent datum of vector bundles can be tested both generically and pointwise. Both of these results should be interesting in their own right. \Cref{SenDescent} is then applied in Subsection 4.5 to show a classicality result on \(\Yoc(A,A^+)\) which can be seen as a rational version of our main result, at least for those integral perfectoids arising as rings of integral elements of an affinoid perfectoid (and thus in particular for \(A^+_{\ana}\)). Subsection 4.6 then deals with what happens in \(p=[\varpi]=0\), or in other words, the extension from \(\Ycc(A,A^+)\) to \(\Spa(\Ainf)\), finishing the proof for \(A^+_{\ana}\) (\Cref{mainThmAnalytic}). The cases of \(A^+_{\disc}\) and \(A^+_{\ana}\) are combined in the last subsection to prove our main result \Cref{equivalenceThm} for general \(I\)-complete rings.\par
In the last section, we add the group structure and Frobenius structure to both categories and show that they are respected by the equivalence. For the group structure, we get this for free by exactness, while the proof for the Frobenius structure is some straightforward linear algebra.

\subsection*{Notation and conventions}

We use the notion of integral perfectoid rings from \cite{BMS19}. For an integral perfectoid \(R^+\), an element \(\varpi\in R^+\) will always be an element such that \(R^+\) is \(\varpi\)-adically complete and \(\varpi^p|p\). On the other hand, for a Tate Huber pair \((A,A^+)\), an element \(\varpi\in A\) will refer to a topologically nilpotent unit. In the latter case, we refer to such an element as a \emph{pseudouniformizer}. Note that for integral perfectoids, \(\varpi\) might not be a non-zero divisor, e.g., in a perfect ring of characteristic \(p\) it can always be chosen as \(0\), so the term pseudouniformizer does not seem appropriate here.\par
We furthermore make use of the whole formalism of adic and perfectoid spaces from \cite{SW20}, always assumed to be \(p\)-complete in the sense that \(p\) is contained in an ideal of definition. For \((B,B^+)\) a perfectoid Huber pair in characteristic \(p\), we will sometimes use \(\Spa(B,B^+)\) and \(\Spd(B,B^+)\) interchangeably. We call a \(v\)-sheaf \emph{representable} if it is of this form. By definition, \((X\dot\times Y)^{\diamond}=X^\diamond \times Y^\diamond\). For a perfectoid Huber pair \((A,A^+)\) with pseudouniformizer \(\varpi\), we will make use of the following adic spaces:
\begin{align*}
\Spa(A^+)\dot\times \Spa(\ZZ_p)&=\Spa(\Ainf(A^+))\\
\Spa(A^+)\dot\times \Spa(\QQ_p)&=\Yoc(A,A^+)=\Spa(\Ainf(A^+))-V(p)\\
\Spa(A,A^+)\dot\times \Spa(\ZZ_p)&=\Yco(A,A^+)=\Spa(\Ainf(A^+))-V([\varpi])\\
\Spa(A,A^+)\dot\times \Spa(\QQ_p)&=\Yoo(A,A^+)=\Spa(\Ainf(A^+))-V(p[\varpi])\\
&\quad\;\;\Ycc(A,A^+)=\Spa(\Ainf(A^+))-V(p,[\varpi]).
\end{align*}
Except for the first one, they are all non-affinoid with an open cover by sousperfectoid analytic affinoids and are thus sheafy (by \cite[Proposition 6.3.4]{SW20}) with a good theory of vector bundles (by \cite[Theorem 5.2.8]{SW20}). The first space \(\Spa(\Ainf(A^+))\) is by definition affinoid, but non-analytic, and it is not clear whether \(\Ainf(A^+)\) is sheafy in general. We will make use of this space as an adic space only in the setting of \Cref{liuAinf}. For discrete perfect rings, we instead use finite projective Witt-vector modules and the associated v-sheaves.\par
Each of these spaces carries a Frobenius endomorphism, which we will always call \(\phi\), while \(\phi\)-linear maps will be named \(\Phi\). There is the common convention to write \(\Spa(A)\) for \(\Spa(A,A^+)\) if \(A^+=A^\circ\) and the topology is understood. For clarity, we will make use of this only if \(A=A^+=A^\circ\) or \(A=\QQ_p\).\par
An element \(d\in \Ainf(S^+)\) for some integral perfectoid \(S^+\) will refer to a distinguished element, i.e., \(\Ainf(S^+)/(d)=S^+\). The adic version of this is 
    \[\xi\in \Gamma(\Yco(B,B^+),\MCO_{\Yco(B,B^+)}),\]
    so that \(V(\xi)=\Spa(B,B^+)\).

\newpage

\section{The schematic side}
\subsection {\texorpdfstring{The perfect-prismatic site and the \(I\)-complete arc topology}
  {The perfect-prismatic site and the I-complete arc topology}}
The following variant of the arc topology was already considered in \cite{BS21} (for \(I=(p)\)), in \cite{Ito} (for \(I=(\varpi)\) and integral perfectoids) and in \cite[Section 4.3]{Takaya} (for varying \(I\)).

\begin{mydef}\label{defArcSite}(The \(I\)-complete arc site)
Let \(R^+\) be an \(I\)-complete ring for some finitely generated ideal \(I\subset R^+\) containing \(p\) and let \(\MCC\) be a class of \(I\)-complete \(R^+\)-algebras, stable under \(I\)-complete tensor products\footnote{In our case, \(\MCC\) will always be either the category of all \(I\)-complete rings or of \(I\)-complete integral perfectoids}. Then the \emph{\(I\)-complete arc topology} on \(\MCC\) has covers given by maps \(A^+\to B^+\) which satisfy the arc lifting property with respect to \(I\)-complete rank-one valuation rings \(V\), i.e., for every map \(A^+\to V\) we must be able to find a commutative diagram

\begin{center}
\begin{tikzcd}
B^+ \arrow[r]           & W                 \\
A^+ \arrow[u] \arrow[r] & V \arrow[u, hook]
\end{tikzcd}
\end{center}

with \(V\to W\) a faithfully flat map (equivalently, an injective local homomorphism) of \(I\)-complete rank-one valuation rings.
\end{mydef}

We use the language of prisms from \cite{BS21}. Recall that a prism is called \emph{perfect} if its Frobenius endomorphism is bijective. Perfect prisms are equivalent to integral perfectoid rings by \((A,J)\mapsto A/J\) (\cite[Theorem 3.10]{BS21}).

\begin{mydef}(The perfect-prismatic site \((R^+)_{\Prism}^{\perf}\))\label{defPrismaticSite}
Let \(R^+\) be a complete topological ring carrying the \(I\)-adic topology for some finitely generated ideal \(I\subset R^+\) containing \(p\). In this setting, we want to define the \emph{perfect-prismatic site} \((R^+)_{\Prism}^{\perf}\) as follows:
\begin{enumerate}
    \item (Underlying category) The underlying category of the perfect-prismatic site over \(R^+\) is the opposite of the category of \(I\)-complete integral perfectoid \(R^+\)-algebras. Equivalently, it is the restriction of the prismatic site over \(R^+\) to the sub-site containing only \((\overline I,d)\)-complete perfect prisms \((\Ainf(A^+),(d))\). Here, \(d\in \Ainf(A^+)\) is a distinguished element (so that \(\Ainf(A^+)/(d)=A^+\)) and \(\overline I\subset \Ainf(A^+)\) is the preimage of \(I\subset A^+\) under the projection.
    \item (Base changes) For a diagram \((C,JC)\leftarrow (A,J)\rightarrow (B,JB)\), it is shown in \cite[Corollary 3.12]{BS21} that the derived \((p,J)\)-completion of \(C\otimes_A^{L} B\) again defines a prism and that this is the pushout in the category of prisms. It is clear that it is perfect if the three prisms in the diagram are, so it serves as the pushout in the category of perfect prisms as well. Finally, by quotienting out \(J\), we see that for perfectoid rings \(R''\leftarrow R \rightarrow R'\), the derived \(p\)-completion of \(R''\otimes_R^{L} R'\) is again perfectoid and thus provides an equivalent definition of base changes in the category of perfect prisms.
    \item (Topology) We can equip this site with two different topologies:
\begin{itemize}
    \item The \emph{\((\overline I,d)\)-complete flat topology} (or just \emph{flat topology} for short) as in \cite[Definition 4.1]{BS21}, obtained by restriction of the topology on the usual prismatic site; here, covers are given by maps \((A,J)\to (B,J')\) such that \(A\to B\) is \((p,\overline I)\)-completely faithfully flat. 
    \item We can also consider the \emph{\(I\)-complete arc topology} on integral perfectoids, i.e., covers are given by maps \((A,J)\to (B,J')\) such that \(A/J\to B/J'\) is an \(I\)-complete arc cover. 
\end{itemize}
The \(I\)-complete arc topology is finer than the flat topology (\Cref{arcFinerThanFlat}). The two descriptions of the base change above show that the covers in both topologies are stable under base change. Strictly speaking, we will use \((R^+)_{\Prism}^{\perf}\) to refer only to the underlying category of the site and will always explicitly mention which topology we put on it whenever it is relevant. 
    \item (Structure sheaf) The usual prismatic site carries a structure sheaf \(\MCO_\Prism\). Restriction to the perfect-prismatic site yields a structure sheaf
\[\MCO_{\Prism}^{\perf}\colon(A,J)\mapsto A.\]
\Cref{takaya} implies that this is even a sheaf for the \(I\)-complete arc topology.
\end{enumerate}
\end{mydef}

We made a claim in the definition which we should prove.

\begin{mylemma}\label{arcFinerThanFlat}
In the context of \Cref{defPrismaticSite}, the \(I\)-complete arc topology is finer than the flat topology.
\end{mylemma}

\begin{proof}
Let \((A,(d))\to (B,(d'))\) be a \((\overline I,d)\)-completely faithfully flat map of perfect prisms. Then \(A/(d, (AI)^n)\to B/(d', (BI)^n)\) is a faithfully flat map of rings for any \(n\in \ZZgez\). By the considerations in \cite[Example (2) at the beginning of Section 2.2.1]{CS21}, this implies that \(A/(d)\to B/(d')\) is an \(I\)-complete arc cover.
\end{proof}
\subsection {\texorpdfstring{Locally free sheaves on the perfect-prismatic site: \(\MCO^{\perf}_{\mathbbl{\Delta}}\)-bundles}
  {Locally free sheaves on the perfect-prismatic site: OΔperf-bundles}}
\begin{mydef}
Let \(R^+\) be an \(I\)-complete ring, \(p\in I\). We say that a presheaf \(\MCE\) of \(\MCO_{\Prism}^{\perf}\)-modules on \((R^+)_{\Prism}^{\perf}\) satisfies the \emph{crystal condition} if for every map of perfect prisms \((A,J)\to (B,J')\) over \(R^+\) the natural map
\[\MCE(A,J)\otimes_A B\longrightarrow\MCE(B,J')\]
is an isomorphism.
\end{mydef}

The following descent result was first proved in \cite[Theorem 1.2]{Ito} (for \(I=(\varpi)\)) and then later expanded to general \(I\) in \cite{Takaya}. In both cases, the result was proved for perfect complexes. We only need the special case of finite projective modules.

\begin{mythm}\label{takaya}\emph{\cite[Proposition 4.17]{Takaya}} Let \(R^+\) be a ring complete with the \(I\)-adic topology for some ideal \(p\in I\subset R^+\). Sending an integral perfectoid \(I\)-complete \(R^+\)-algebra \(S^+\) to the category of finite projective \(\Ainf(S^+)\)-modules satisfies \(I\)-complete arc descent.
\end{mythm}

\begin{mydefprop}\label{equivalentDefsOfPrismaticBundles}
Let \(R^+\) again be any \(I\)-complete ring. An \emph{\(\MCO^{\perf}_{\Prism}\)-bundle} over \(R^+\) is a presheaf \(\MCE\) of \(\OpP\)-modules on \((R^+)_{\Prism}^{\perf}\) subject to one of the following equivalent conditions:

\begin{enumerate}
    \item \(\MCE\) is a sheaf for the topology \(\tau\) and finite locally free in the topology \(\tau'\), for any choice \((\tau,\tau')\in \{(\Flat,\Flat),(\arc,\arc),(\arc,\Flat)\}\).
    \item \(\MCE\) sends any prism \((A,J)\) to a finite projective \(A\)-module and satisfies the crystal condition.
\end{enumerate}
We equip the category of \(\MCO^{\perf}_{\Prism}\)-bundles with the exact structure coming from the flat topology, so that a complex 
\[0\longrightarrow\MCE_1\longrightarrow\MCE_2\longrightarrow\MCE_3\longrightarrow0\] is exact if it is exact as a sequence of sheaves on the perfect-prismatic site using the flat topology -- In \Cref{exactnessArcLocallyGeneral} we will see that, once again, this is equivalent to defining it through the arc topology.
\end{mydefprop}

\begin{proof}
    First, every rule as in (2) is an arc sheaf by \Cref{takaya}, and flat-locally free as it is Zariski-locally free on \(A\), and a completed Zariski-localization of a perfect prism is again a perfect prism (see \cite[Remark 2.16]{BS21} for localization of delta rings) and the completed localization map will certainly be a \((p,I)\)-complete flat cover.\par
    Now every arc sheaf is a flat sheaf and every sheaf that is flat-locally free is also arc-locally free, so the variant \((\tau,\tau')=(\arc,\Flat)\) from the previous paragraph implies the other two.\par
    Finally, we have to show that a \(\tau\)-sheaf \(\MCE\) that is \(\tau=\tau'\)-locally free only admits values in finite projective modules and satisfies the crystal condition. The former is a direct consequence of \Cref{takaya} (or its variant for the flat topology, following from \Cref{arcFinerThanFlat}), the latter is a standard argument: Let \((A,I)\to (B,J)\) be a map in \((R^+)_{\Prism}^{\perf}\) and let \((A,I)\to( A', I')\) be a \(\tau\)-cover trivializing \(\MCE\). Write
\begin{align*}
(A'',I'')&=( A', I')\hat\otimes_{(A,I)}( A', I')\\
( B', J')&=( A', I')\hat\otimes_{(A,I)}(B,J)\\
( B'',J'')&=(A'',I'')\hat\otimes_{(A,I)}(B,J)
\end{align*}

for the base changes as defined in \Cref{defPrismaticSite} (2). Then \(\MCE\) is also free on all those. We now have a commutative diagram

\begin{center}
\begin{tikzcd}
0 \arrow[r] & {\MCE(A,I)\otimes_A B} \arrow[d] \arrow[r] &  \MCE( A',I')\otimes_{A'} B' \arrow[d] \arrow[r] & { \MCE(A'',I'')\otimes_{A''} B''} \arrow[d] \\
0 \arrow[r] & {\MCE(B,J)} \arrow[r]                      & { \MCE( B', J')} \arrow[r]                               & \MCE(B'',J'').                        
\end{tikzcd}
\end{center}

The rows are exact by the sheaf condition (here it is relevant that \(\tau=\tau'\)), and the two vertical arrows on the right are isomorphisms because we already know \(\MCE\) to be free on those prisms and because tensor products commute with finite products. The result follows by the five lemma. 
\end{proof}

\begin{mycor}\label{CrysForPerfdRings}
For \(R^+\) an \(I\)-complete integral perfectoid ring, there are mutually inverse equivalences between \(\OpP\)-bundles on \(R^+\) and finite projective \(\Ainf(R^+)\)-modules given by the obvious constructions, namely: evaluating the \(\OpP\)-bundle at \((\Ainf(R^+),(d))\) resp. sending an \(\Ainf(R^+)\)-module \(M\) to \((A,J)\mapsto M\otimes_{\Ainf(R^+)} A\).
\end{mycor}

\begin{proof}
Recall that by the equivalent description as perfectoid \(R^+\)-algebras in \Cref{defPrismaticSite}, the perfect-prismatic site over a perfectoid \(R^+\) has \(\Ainf(R^+)\) as final object. The equivalence of categories then follows immediately from characterization (2) in \Cref{equivalentDefsOfPrismaticBundles}.
\end{proof}

For the final descent result of this subsection, we need a slight refinement of a result by Scholze-\v{C}esnavicius, a more precise version of which we will also prove later in \Cref{lemmaCoveringg}.

\begin{mylemma}\label{LemmaSC}\emph{\cite[Lemma 2.2.3]{CS21}} 
    Let \(R^+\) be an \(I\)-complete ring for some finitely generated ideal \(I\subset R^+\) containing \(p\). Then \(R^+\) has an \(I\)-complete arc cover \(R^+\to A^+\) where \(A^+=\prod_i V_i\) is a product of \(I\)-complete valuation rings with algebraically closed field of fractions. In particular, \(A^+\) is integral perfectoid.
\end{mylemma}

\begin{myprop}\label{arcDescentForCrysNew}
    Let \(R^+\) be an \(I\)-complete ring. Sending an \(I\)-complete \(R^+\)-algebra \(S^+\) to the category of \(\OpP\)-bundles on \(S^+\) satisfies \(I\)-complete arc descent.
\end{myprop}

\begin{proof}
    This follows from \Cref{LemmaSC}, \Cref{CrysForPerfdRings} and \Cref{takaya}: Let 
    \[R^+\longrightarrow S^+\rightrightarrows S'^+   \rrrarrows  S''^+\]
    be a truncated \(I\)-complete arc hypercover with a descent datum of \(\OpP\)-bundles. Let \(R^+\to R^{+}_{\perfd}\) be a map to an \(I\)-complete integral perfectoid. By successively taking base changes and \(I\)-complete arc covers as in \Cref{LemmaSC}, we get a map of truncated hypercovers:
    \begin{center}
\begin{tikzcd}
R^+_{\perfd} \arrow[r]   & S^+_{\perfd} \arrow[r, shift left] \arrow[r, shift right]  & S'^+_{\perfd} \arrow[r, shift left=2] \arrow[r] \arrow[r, shift right=2]  & S''^+_{\perfd}   \\
R^+ \arrow[u] \arrow[r] & S^+ \arrow[u] \arrow[r, shift left] \arrow[r, shift right] & S'^+ \arrow[u] \arrow[r, shift left=2] \arrow[r] \arrow[r, shift right=2] & S''^+ \arrow[u]
\end{tikzcd}
    \end{center}
where the top row consists of \(I\)-complete integral perfectoids. Pulling back the descent datum of \(\OpP\)-bundles, we get (by \Cref{CrysForPerfdRings}) a descent datum of finite projective modules, which (by \Cref{takaya}) descends to a single finite projective module on \(R^+_{\perfd}\). It is now easy to check that this construction satisfies the crystal condition as \(R^+_{\perfd}\) varies, showing the result.
\end{proof}
\subsection{\texorpdfstring{Exact structure on \(\MCO^{\perf}_{\mathbbl{\Delta}}\)-bundles}
  {Exact structure on OΔperf-bundles}}
\begin{mylemma}\label{exactnessArcLocally}
    Let \(R^+\to S^+\) be an \(I\)-complete arc cover between \(I\)-complete integral perfectoid rings for some \(I\subset R^+\) containing \(p\). Then a complex
    \[0\longrightarrow M_1\longrightarrow M_2\longrightarrow M_3\longrightarrow0\]
    of finite projective \(\Ainf(R^+)\)-modules is exact if and only if its base change to \(\Ainf(S^+)\) is.
\end{mylemma}

\begin{proof}
A short exact sequence of finite projective modules splits, so it remains short exact after base change. For the converse, we follow the argumentation in \cite[Lemma 2.3.6]{PR}.
For any ring \(R\), a complex 
\[0\longrightarrow N_1\longrightarrow N_2\longrightarrow N_3\longrightarrow0\]
of finite projective \(R\)-modules is short exact if for every maximal ideal \(\mfm\subset R\) the base change to the residue field is short exact. Indeed, let \(x\in \Spec(R)\) be a closed point. Then Nakayama's lemma implies surjectivity of the second map in an open neighbourhood \(U\ni x\). Shrinking \(U\) if necessary, we obtain a splitting of this map, so we can assume that the kernel is finite projective. Now the first map factors through this finite projective kernel and is, by assumption, an isomorphism on residue fields. Another application of Nakayama's lemma and splitting shows that it is an isomorphism on a neighbourhood of \(x\) and, by varying \(x\), on \(\Spec(R)\).\par
Using this and faithful flatness of field extensions, it therefore suffices to show that every closed point of \(\Spec(\Ainf(R^+))\) lifts to a point of \(\Spec(\Ainf(S^+))\). Choose a distinguished element \(d\in\Ainf(R^+)\) with \(\Ainf(R^+)/(d)=R^+\). Since \(\Ainf(R^+)\) is \((p,d)\)-adically complete, \(d\) lies in its Jacobson radical, and thus all its closed points lie in \(\Spec(R^+)\). Furthermore, by \(I\)-completeness, every closed point of \(\Spec(R^+)\) lies in \(\Spec(R^+/I)\). The map \(R^+/I\to S^+/IS^+\) is an arc cover in the usual sense, so the associated map of schemes is surjective. Hence a closed point of \(\Spec(R^+/I)\) lifts to \(\Spec(S^+)\) and thus to \(\Spec(\Ainf(S^+))\).
\end{proof}

Using this, the following facts are easily verified. The proofs work in the same way as those in Section 3.3, namely \Cref{exactnessVBundleYBundleEquivalence} and \Cref{exactnessIsVLocal}, though simpler, as by definition each object in \((R^+)_{\Prism}^{\perf}\) is qcqs in the sheaf-theoretic sense.

\begin{mylemma}\label{exactnessArcLocallyGeneral}
    \begin{enumerate}
        \item Let \(R^+\) be an integral perfectoid ring. Then a complex \(0\to M_1\to M_2\to M_3\to0\) of finite projective \(\Ainf(R^+)\)-modules is short exact if and only if the associated complex of \(\OpP\)-modules on \(R^+\) is short exact as a sequence of sheaves on \((R^+)_{\Prism}^{\perf}\), equipped with the \(I\)-complete arc topology. The same equivalence holds for the flat topology.
        \item Let \(R^+\to A^+\) be an \(I\)-complete arc cover of \(I\)-complete rings. Then a complex \(0\to\MCE_1\to\MCE_2\to\MCE_3\to0\) of \(\OpP\)-bundles on \(R^+\) is short exact if and only if its base change to \(A^+\) is.
    \end{enumerate}
    \end{mylemma}

\section{The adic side}
\subsection{The v-site}
We now start dipping into the adic world, namely v-sheaves on the category \(\Perf\) of perfectoid Huber pairs in characteristic \(p\) and certain locally free sheaves on them.\par

\begin{mynot}\label{notationVSheaves}(Letters for v-sheaves)
Just for this section, we use the following notational convention: We will denote a v-sheaf that comes with a structure map to \(\Spd(\ZZ_p)\) by an ordinary letter with a diamond subscript, such as \(X_{\diamond}\). It may or may not (but in our applications always will) be given as \(X^\diamond\) for some adic space \(X\). In particular, whenever we have a map from a representable \(\Spa(B,B^+)\to X_{\diamond}\), we get a canonical untilt \((B^{\sharp},B^{+\sharp})\) supplied by the composition \(\Spa(B,B^+)\to X_{\diamond}\to \Spd(\ZZ_p)\). In contrast, v-sheaves without a structure map to \(\Spd(\ZZ_p)\) will be denoted by cursive letters like \(\MCX\).
\end{mynot}

\begin{mydef}
    A map \(\MCX'\to \MCX\) of v-sheaves is called a \emph{v-cover} or simply \emph{cover} if it is surjective in the sheaf-theoretic sense. Explicitly, (using that affinoid perfectoids are qcqs in the sheaf-theoretic sense) for any map \(\Spa(B,B^+)\to \MCX\) from a representable one has to be able to find a commutative diagram
    \begin{center}
        \begin{tikzcd}
        {\Spa(B',B'^+)} \arrow[d] \arrow[r] & \MCX' \arrow[d] \\
        {\Spa(B,B^+)} \arrow[r]             & \MCX           
        \end{tikzcd}
    \end{center}
    with \(\Spa(B',B'^+)\) also representable and the left map being a v-cover of affinoid perfectoids.
\end{mydef}

Recall that a v-sheaf is called \emph{small} if it admits a cover by a perfectoid space, which is really just a set-theoretic bound and holds for most v-sheaves arising in practice. For v-sheaves arising from \(I\)-complete rings one can even cover them by a single representable, which will be useful later:

\begin{mylemma}\label{coverByRepresentable}
Let \(A^+\) be a complete topological ring carrying the \(I\)-adic topology for some finitely generated ideal \(I\subset A^+\) containing \(p\). Then \(\Spd(A^+)\) is quasi-compact (in the sheaf-theoretic sense). Equivalently, there is a v-cover by a representable
\[\Spa(B,B^+)\longrightarrow \Spd(A^+).\]
\end{mylemma}

\begin{proof}
    It is enough to find a cover of \(\Spd(A^+)\) by a quasi-compact v-sheaf. We know from \cite[Chapter 15]{etCohDiamonds} that the v-sheaves associated to affinoid analytic adic spaces are spatial diamonds, so in particular quasi-compact in the sheaf-theoretic sense, so it suffices to find a cover by a finite number of these.\par
    By properties (1), (2) and (4) of \Cref{lemmaCoveringg}, we can assume that \(A^+\) is either discrete or that the \(I\)-adic topology coincides with the \(\varpi\)-adic one for some non-zero divisor \(\varpi\in A^+\) and \(A=A^+[1/\varpi]\). In the first case,
    \[(A^+,A^+)\longrightarrow(A^+((t)),A^+[[t]])\] is clearly a v-cover by an analytic adic space (using that every affinoid perfectoid has a topologically nilpotent unit to which we can send \(t\), lifting a map from \((A^+,A^+)\)). In the second case, give \(A^+[[t]]\) the \((\varpi,t)\)-adic topology and consider
    \[U\coloneqq\Spa(A^+[[t]])\setminus V(\varpi,t).\]
    This has a cover by the two analytic affinoid charts \(U_{t\leq\varpi\neq0}\) and \(U_{\varpi\leq t\neq0}\). Now every map \((A^+,A^+)\to(B,B^+)\) to a perfectoid Huber pair lifts to \(U\) by sending \(t\) to a pseudouniformizer of \(B\). Thus the union of the two charts gives the required quasi-compact v-cover.
\end{proof}

The following are sheaves by \cite[Theorem 17.1.3]{SW20}:

\begin{mydef}\label{defVSheaf}
    Let \(\MCX\) be a v-sheaf. Then on the v-site \(\Perf_{/\MCX}\) of characteristic \(p\) affinoid perfectoids over \(\MCX\) one defines a \emph{tilted structure sheaf}
    \[\MCO^{\flat\Perf}_{/\MCX}\colon (\Spa(B,B^+)\to \MCX)\mapsto B.\]
    If \(\MCX\) comes equipped with a map to \(\Spd(\ZZ_p)\), so we would call it \(X_\diamond\) in line with \Cref{notationVSheaves}, we furthermore have an \emph{untilted structure sheaf}:
        \[\MCO^{\sharp\Perf}_{/X_\diamond}\colon (\Spa(B,B^+)\to X_{\diamond})\mapsto B^\sharp.\]
\end{mydef}

Let's finish this subsection with a little aside on \(\Spa(\Ainf(A^+))\). In general, it is not clear whether \(\Spa(\Ainf(A^+))\) exists as an adic space (i.e., whether \((\Ainf(A^+), \Ainf(A^+))\) is sheafy) or if it has a good theory of vector bundles. For the analytic case in the proof of our main theorem, we use the following result:

\begin{mylemma}\label{liuAinf}
    Let \(A^+\) be an integral perfectoid topological ring carrying the \(\varpi\)-adic topology for some non-zero divisor \(\varpi\in A^+\) such that \(\varpi^p|p\).
    Then \((\Ainf(A^+), \Ainf(A^+))\) is sheafy and vector bundles on \(\Spa(\Ainf(A^+))\) are in bijection with finite projective \(\Ainf(A^+)\)-modules.
\end{mylemma}

\begin{proof}
    Sheafiness is \cite[Theorem 1.1]{Liu24}, and the statement about vector bundles is \cite[Theorem 1.4]{Liu24}.
\end{proof}

\begin{myrem}\label{discreteWittSheafiness}
    For a discrete perfect ring \(R\) of characteristic \(p\), we work with finite projective \(W(R)\)-modules and the associated v-sheaves instead of vector bundles on \(\Spa(W(R))\), allowing us to avoid the sheafiness assumption. The same might be possible to do for \(A^+\) as above, though being able to speak about vector bundles on \(\Spa(\Ainf(A^+))\) simplifies some proofs.
\end{myrem}
\subsection{\texorpdfstring{Locally free sheaves on the v-site: v-bundles and \(\MCY\)-bundles}
  {Locally free sheaves on the v-site: v-bundles and Y-bundles}}
The following notion seems to be well-known, e.g., a special case is studied in \cite{Heu21}:

\begin{mydef}\label{defVBundle}
Let \(X_{\diamond}\) be a small v-sheaf over \(\Spd(\ZZ_p)\). Then a \emph{v-bundle} \(\MCE\) on \(X_\diamond\) is one of the following equivalent pieces of data:
\begin{enumerate}
    \item A finite locally free \(\MCO^{\sharp\Perf}_{/X_\diamond}\)-module,
    \item A rule assigning to any \((B,B^+)\in \Perf\) and \((B,B^+)\)-valued point \(\alpha \colon\Spa(B,B^+)\to X_{\diamond}\) with associated untilt \((B^{\sharp},B^{+\sharp})\) a vector bundle on \(\Spa(B^{\sharp},B^{+{\sharp}})\) which we will call \(\alpha^*\MCE\), or, if \(\alpha\) is clear from the context, \(\MCE|_{\Spa(B^{\sharp},B^{+{\sharp}})}\). This rule has to be compatible with pullbacks, more precisely: For any map \(f\colon\Spa(A,A^+)\to \Spa(B,B^+)\) in \(\Perf_{/X_\diamond}\), which also gives us a map \(f^{\sharp}\colon \Spa(A^{\sharp},A^{+\sharp})\to \Spa(B^{\sharp},B^{+\sharp})\), we want \(f^{\sharp*}(\alpha^*\MCE)=(\alpha\circ f)^* \MCE\). Here \(f^{\sharp*}\) just denotes the usual pullback of vector bundles along maps of adic spaces.
    \item A 3-truncated v-hypercover
    \[X'''^\diamond\;\substack{\rightarrow\\[-1em] \rightarrow \\[-1em] \rightarrow}\; X''^\diamond\rightrightarrows X'^\diamond\longrightarrow X_{\diamond}\]
    by v-sheaves associated to perfectoid spaces not necessarily of characteristic \(p\) (always possible by smallness of \(X_{\diamond}\)), with the maps assumed to be compatible with the structure maps to \(\Spd(\ZZ_p)\),
     together with a descent datum of vector bundles on 
    \[X''' \rrrarrows X''\rightrightarrows X'\]
    modulo the obvious equivalence relation.
\end{enumerate}
\end{mydef}

The equivalence of these definitions follows from v-descent of vector bundles on perfectoid spaces (\cite[Lemma 17.1.8]{SW20}) through a simple proof very similar to \Cref{equivalentDefsOfPrismaticBundles} which will not be repeated here. The second definition makes it clear that we can pull back v-bundles along maps of v-sheaves over \(\Spd(\ZZ_p)\).\par

Even more relevant to us will be a slightly different notion:

\begin{mydef}\label{defYBundle}
Let \(\MCX\) be any small v-sheaf and set \(\MCS=\ZZ_p\) or \(\MCS=\QQ_p\). Then an \emph{integral} (if \(\MCS=\ZZ_p\)) resp. a \emph{rational} (if \(\MCS=\QQ_p\)) \(\MCY\)-bundle on \(\MCX\) is given by one of the following equivalent pieces of data:
\begin{enumerate}
    \item A rule associating to any map \(\Perf\ni\Spa(B,B^+)\to \MCX\) a vector bundle on \[\Spa(B,B^+)\dot\times\Spa(\MCS)=
\begin{cases}
\Yco(B,B^+)\quad &\text{if }\MCS=\ZZ_p\\
\Yoo(B,B^+)\quad &\text{if }\MCS=\QQ_p,
\end{cases}\]
compatible with pullbacks along maps \(\Spa(A,A^+)\to \Spa(B,B^+)\) in \(\Perf_{/\MCX}\), cf. \Cref{defVBundle}.
\item A 3-truncated v-hypercover 
    \[X''' \rrrarrows X''\rightrightarrows X'\longrightarrow \MCX\]
    by perfectoid spaces in characteristic \(p\) with a descent datum of vector bundles on 
    \[X'''\dot\times\Spa(\MCS) \rrrarrows X''\dot\times\Spa(\MCS)\rightrightarrows  X'\dot\times\Spa(\MCS).\]
\end{enumerate}
\end{mydef}

The equivalence is again proved in the same way as \Cref{equivalentDefsOfPrismaticBundles}, now using v-descent of vector bundles on open subspaces of \(\Yco\) (\cite[Proposition 19.5.3]{SW20}). Again, the first definition gives us a pullback functor from rational/integral \(\MCY\)-bundles along maps of v-sheaves \(\MCX\to \MCX'\).

\begin{myrem}
    In contrast to \(\OpP\)-bundles (\Cref{equivalentDefsOfPrismaticBundles}) and v-bundles (\Cref{defVBundle}), integral/rational \(\MCY\)-bundles have no definition as locally free sheaves on some site. This is because \(\Yco(R,R^+)\) (resp. \(\Yoo(R,R^+)\)) is not affinoid, so vector bundles on it do not directly correspond to modules over some ring. Also, it is not even clear to me if for an algebraically closed field \(C\) every vector bundle on \(\Yoo(C,C^+)\) is trivial, so there seems to be little hope of trivializing these bundles locally in the v-topology on \(\MCX\).
\end{myrem}

\begin{myrem}\label{ComparisonGI}
    In \cite[Section 3]{GI23}, Gleason, Ivanov and Zillinger study related categories of vector bundles on relative adic spaces and over Witt vector rings. Their notation differs from ours.
\end{myrem}

\begin{myprop}\label{vDescentForVBundles}
The functors
\begin{enumerate}
    \item associating to a small v-sheaf \(X_{\diamond}\) over \(\Spd(\ZZ_p)\) the category of v-bundles on \(X_\diamond\),
    \item associating to a small v-sheaf \(\MCX\) the category of integral/rational \(\MCY\)-bundles on \(\MCX\)
\end{enumerate}
 satisfy descent along covers of small v-sheaves (assumed to be compatible with the maps to \(\Spd(\ZZ_p)\) in the first case).
\end{myprop}

\begin{proof}
By the already cited descent results \cite[17.1.8]{SW20} and \cite[19.5.3]{SW20} both are sheaves on affinoid perfectoids. The proof now proceeds like the one of \Cref{arcDescentForCrysNew} with some added bookkeeping because of the non quasi-separatedness of general small v-sheaves.
\end{proof}

In cases where this makes sense, there are base change functors from vector bundles to v-bundles, from \(\MCY\)-bundles to v-bundles and from vector bundles to \(\MCY\)-bundles; the constructions work in somewhat different generalities, so let's name all three. A diagram of those functors will follow at the end of the subsection.

\begin{myconst}\label{baseChangFunctors}
 \begin{enumerate}
    \item (vector bundles to v-bundles) Let \(X\) be any adic space over \(\ZZ_p\), and consider its associated v-sheaf \(X^\diamond\). We get a functor from vector bundles on \(X\) to v-bundles on \(X^\diamond\) as follows: Let \(\MCE\) be a vector bundle on \(X\). The \(\Spa(B,B^+)\)-valued points of \(X^\diamond\) are by definition given by an untilt \((B^{\sharp},B^{+\sharp})\) together with a map \(\Spa(B^{\sharp},B^{+\sharp})\to X\). By definition of the map \(X^\diamond\to\Spd(\ZZ_p)\), this untilt agrees with the one given by the composition
    \[\Spa(B,B^+)\longrightarrow X^\diamond\longrightarrow \Spd(\ZZ_p).\]
    Pullback of \(\MCE\) along \(\Spa(B^\sharp,B^{+\sharp})\to X\) defines a rule for a v-bundle, and compatibility with pullbacks is obvious by the definition.
\item (\(\MCY\)-bundles to v-bundles) Now consider any small v-sheaf \(\MCX\). In this context, there is a functor from integral (resp. rational) \(\MCY\)-bundles on \(\MCX\) to v-bundles on \(\MCX\times\Spd(\MCS)\) for \(\MCS=\ZZ_p\) (resp. \(\MCS=\QQ_p\)), where \(\MCX\times\Spd(\MCS)\) is a v-sheaf over \(\Spd(\ZZ_p)\) through projection to the second component. Let \(\MCE\) be such a \(\MCY\)-bundle. We need to assign to every map from a representable \(\Perf\ni T\xrightarrow{(g,f)} \MCX\times\Spd(\MCS)\) a vector bundle on \(T^\sharp\). But \((g,f)\) factors as
\[T \xrightarrow{(\id,f)} T\times \Spd(\MCS) \xrightarrow {(g,\id)} \MCX\times\Spd(\MCS)\]
in the category of v-sheaves over \(\Spd(\ZZ_p)\), and we can untilt the left part of the diagram to get a map
\[T^\sharp \longrightarrow T^\sharp\dot\times\Spa(\MCS).\]
The datum of an integral/rational \(\MCY\)-bundle on \(T\) now gives us a vector bundle on \(T^\sharp\dot\times\Spa(\MCS)\), which we can pull back to \(T^\sharp\). Compatibility with maps \(T\to T'\) is straightforward to check.
    \item (vector bundles to \(\MCY\)-bundles) Now keep the assumption from (2) and assume furthermore that \(Y_{\diamond}=\MCX\times \Spd(\MCS)\) is given as \(Y^\diamond\) for a suitable adic space \(Y\) -- this is the case notably if \(\MCX\) is given by
    \begin{itemize}
        \item a perfectoid Huber pair \(\Spa(B,B^+)\) in characteristic \(p\)
        \item or \(\Spd(A^+)\) for an integral perfectoid \(A^+\) equipped with the \(\varpi\)-adic topology for a suitable pseudouniformizer and \(\MCS=\QQ_p\).\footnote{Note the following subtlety: If \(A^+\) is perfect of characteristic \(p\) with the discrete topology, then \(Y=\Spa(W(A^+)[1/p],W(A^+))\); meanwhile, if \(A^+\) carries the \(\varpi\)-adic topology for some pseudouniformizer \(\varpi\) which is a non-zero divisor, so that \((A^+[1/\varpi],A^+)\) is a perfectoid Huber pair, then \(Y=\Spa(A^+)\dot\times\Spa(\QQ_p)=\Yoc(A^+[1/\varpi],A^+)\); i.e., the inversion of \(p\) can happen in either a rational or an analytic way, and which one occurs might depend completely on the choice of \(\varpi\) and, with it, the choice of topology on \(A^+\).}
    \end{itemize}
    In this context we can associate to any vector bundle \(\MCE\) on \(Y\) an integral/rational \(\MCY\)-bundle on \(\MCX\), very similarly to (1): Let \(\Spa(B,B^+)\to \MCX\) be any map from a representable. Then we can simply associate to this map the pullback of \(\MCE\) along
    \[\Spa(B,B^+)\dot\times\Spa(\MCS)\longrightarrow Y\]
    obtained by untilting the map
    \[\Spa(B,B^+)\times\Spd(\MCS)\longrightarrow \MCX\times\Spd(\MCS).\]
        \item (integral to rational \(\MCY\)-bundles) For any small v-stack \(\MCX\), there is also an obvious functor from integral to rational \(\MCY\)-bundles.
\end{enumerate}
In cases where it makes sense, these constructions commute.
\end{myconst}

\begin{myrem}\label{baseChangeDiagrams}
Let's try to summarize the situation in a single commutative diagram. To make sure all the categories make sense, we only consider the following cases. For the rest of the remark, \((B,B^+)\) will always be a perfectoid Huber pair and \(A^+\) an integral perfectoid.\footnote{See \Cref{liuAinf} for why we cannot allow \(X=\Spa(A^+)\) and \(\MCS=\ZZ_p\).}
\begin{itemize}
    \item \(X=\Spa(B,B^+)\) and \(\MCS=\QQ_p\).
    \item \(X=\Spa(B,B^+)\) and \(\MCS=\ZZ_p\).
    \item \(X=\Spa(A^+)\) and \(\MCS=\QQ_p\).
    \item \(X=\Spa(B^+)\) and \(\MCS=\ZZ_p\).
\end{itemize}

We have a commutative diagram:
    
\begin{center}
    \begin{tikzcd}
\text{vector bundles on } X\dot\times\Spa(\MCS) 
\arrow[dd, "\MCF_{\text{vec} \to \text{v}}"', bend right=49, shift right=22] 
\arrow[d, "\MCF_{\text{vec} \to \MCY}"] \\
\text{ integral/rational \(\MCY\)-bundles on } X 
\arrow[d, "\MCF_{\MCY \to \text{v}}"] \\
\text{v-bundles on } (X\dot\times\Spa\MCS)^\diamond                                                                 
\end{tikzcd}
\end{center}

Those maps have the following properties:
\begin{itemize}
    \item \(\MCF_{\MCY \to \text{v}}\) is fully faithful for both \(\MCS=\ZZ_p\) and \(\MCS=\QQ_p\), by \Cref{baseChangeFunctorsProperties} (2). For \(X=\Spa(B,B^+)\), the same holds for \(\MCF_{\text{vec} \to \text{v}}\), by (1). For \(X=\Spa(B^+)\), full faithfulness also follows from the equivalences below. These functors need not be essentially surjective.
   \item \(\MCF^{}_{\text{vec} \to \MCY}\) is an equivalence of categories in the following cases: for \(X=\Spa(B,B^+)\) this is trivial; for \(X=\Spa(B^+)\) and \(\MCS=\QQ_p\) it is \Cref{analyticLocuss}; and for \(X=\Spa(B^+)\) and \(\MCS=\ZZ_p\) it is a special case of our main result.
\end{itemize}
\end{myrem}

\begin{mylemma}\label{lemmaVComplete}
For any affinoid perfectoid \((R,R^+)\) of characteristic \(p\), the affinoid analytic adic spaces
\[X_n\coloneqq\MCY_{[0,p^n]}(R,R^+)\]
form an increasing open cover of \(\Yco(R,R^+)\) and are v-complete, i.e., the natural map
\[H^0(X_n,\MCO_{X_n})\longrightarrow H^0((X_n^\diamond)_v,\MCO_v)\]
is an isomorphism for every \(n\).
\end{mylemma}

This can be seen as a slight extension of a result from \cite[Proposition 9.9]{KimIgusa}, proving stable v-completeness (i.e., open subsets are also v-complete) of \(\Yoo(R,R^+)=\Yco(R,R^+)\setminus V(p)\), the locus away from \(p=0\).

\begin{proof}
Choose a pseudouniformizer \(\varpi\in R^+\) and define rings
    \[
        A^+=W(R^+),\qquad K_\infty=\widehat{\QQ_p(p^{1/p^\infty})},\qquad
        B^+=A^+\otimes_{\ZZ_p}\MCO_{K_\infty},
    \]
    and, with \(t_n=[\varpi^{1/p^n}]\), for \(n\geq0\) and \(m\geq1\)
    \[
        A_{n,0}=\bigl(A^+[p/t_n]\bigr)^{\wedge(t_n)},\qquad
		A_n=A_{n,0}[1/t_n],\qquad
		A_{n,m}=A_n\langle t_n^m/p\rangle,
	\]
	\[
        B_{n,0}=\bigl(B^+[p/t_n]\bigr)^{\wedge(t_n)},\qquad
		B_n=B_{n,0}[1/t_n],\qquad
		B_{n,m}=B_n\langle t_n^m/p\rangle,
    \]
    where \((\cdot)^{\wedge(t_n)}\) denotes \(t_n\)-adic completion.
	Setting \(A_n^+\) as the integral closure of \(A_{n,0}\) inside \(A_n\) (and similarly for the other rings), we get increasing affinoid-analytic maps 
	\[\Spa(B_n,B_n^+)\longrightarrow X_n=\Spa(A_n,A_n^+)=\MCY_{[0,p^n]}(R,R^+)\subset \Yco(R,R^+) \]
	\[\Spa(B_{n,m},B_{n,m}^+)\longrightarrow \Spa(A_{n,m},A_{n,m}^+)=\MCY_{[p^n/m,p^n]}(R,R^+)\subset X_n(R,R^+) \]
		where the \(X_n\) cover \(\Yco(R,R^+)\) as \(n\) varies, and the annular charts cover \(X_n\setminus V(p)\) as \(m\) varies for fixed \(n\),
with the \(B\)-versions furthermore affinoid-perfectoid.\par
Now base change of the \(\ZZ_p\)-linear retraction \(\lambda:\MCO_{K_\infty}\to\ZZ_p\) (obtained by restriction of a retraction of \(\QQ_p\to K_\infty\) sending the basis to elements of \(\ZZ_p\) of norm \(\leq 1\)) to \(A_n\) and \(A_{n,m}\) yields retractions
    \[
        \sigma_n:B_n\longrightarrow A_n,\qquad
        \sigma_{n,m}:B_{n,m}\longrightarrow A_{n,m}
    \]
	compatible with the restriction maps for different \(n,m\) (cf. \cite[proof of Proposition II.1.1]{FS21}). \par
    For these perfectoid covers and any fixed \(m\geq1\), restriction gives an injection
    \[B_n\longrightarrow B_{n,m}\]
    by the perfectoid-disc description in \cite[Proposition II.1.1 and its proof]{FS21} and \cite[Theorem 5.2.1]{SW20}. Now let \(s\in H^0((X_n^\diamond)_v,\MCO_v)\), pull it back to \(b_n\in B_n\), and put \(a_n=\sigma_n(b_n)\in A_n\). By the stable v-completeness of \(\Yoo(R,R^+)\) (\cite[Proposition 9.9]{KimIgusa}), the restriction of \(s\) to \(\Spa(A_{n,m},A_{n,m}^+)^\diamond\) comes from a unique \(a_{n,m}\in A_{n,m}\). Compatibility with pullback implies that the image of \(a_{n,m}\) in \(B_{n,m}\) is the restriction of \(b_n\). The retractions commute with restriction, as in the diagram
\begin{center}
\begin{tikzcd}
A_n \arrow[d, hook, shift left] \arrow[r] & {A_{n,m}} \arrow[d, hook, shift left] \\
B_n \arrow[r, hook] \arrow[u, "\sigma_n", shift left] & {B_{n,m}} \arrow[u, "{\sigma_{n,m}}", shift left]
\end{tikzcd}
\end{center}
so \(a_n|_{A_{n,m}}=\sigma_{n,m}(b_n|_{B_{n,m}})=a_{n,m}\). Thus \(b_n-a_n\) restricts to zero in \(B_{n,m}\), and injectivity of the bottom map gives \(b_n=a_n\). Since \(\Spd(B_n,B_n^+)\to X_n^\diamond\) is a v-cover, \(a_n\) induces \(s\). Injectivity of the scalar comparison follows from the split injection \(A_n\to B_n\), proving the assertion for every \(n\).
\end{proof}

\begin{myprop}\label{baseChangeFunctorsProperties}
\begin{enumerate}
    \item The functor from vector bundles on \(X\) to v-bundles on \(X^\diamond\) is an equivalence of categories for \(X\) perfectoid. It is fully faithful for \(X=\Yco(S)\) or \(X=\Yoo(S)\), where \(S\) is a perfectoid space of characteristic \(p\).
    \item For any small v-sheaf \(\MCX\), the functor from integral (resp. rational) \(\MCY\)-bundles on \(\MCX\) to v-bundles on \(\MCX\times\Spd(\ZZ_p)\) (resp. \(\MCX\times\Spd(\QQ_p)\)) is fully faithful.
    \item If \(X\) is a perfectoid adic space, then the functor from vector bundles on \(\Yco(X)\) (resp. \(\Yoo(X)\)) to integral (resp. rational) \(\MCY\)-bundles on \(X\) is an equivalence of categories.
\end{enumerate}
\end{myprop}

\begin{proof}
\begin{enumerate}
	\item[(1)] The perfectoid case follows from v-descent of vector bundles, \cite[Lemma 17.1.8]{SW20}. For the remaining statements, analytic gluing allows us to assume that \(S=\Spa(R,R^+)\) is affinoid perfectoid. For \(X=\Yco(S)\), use the increasing affinoid cover \((X_n)_n\) from \Cref{lemmaVComplete}. Given vector bundles \(\MCE,\MCF\) on \(X\), their Hom-bundle \(\MCE^\vee\otimes\MCF\) restricts on each \(X_n\) to the bundle associated with a finite projective \(H^0(X_n,\MCO_{X_n})\)-module. Passing to direct summands of finite free modules therefore extends the scalar comparison in \Cref{lemmaVComplete} to this Hom-bundle, proving full faithfulness on each \(X_n\). Uniqueness makes the lifted morphisms compatible on the increasing cover, so they glue to a unique morphism on \(X\). For \(X=\Yoo(S)\), the same argument on an affinoid open cover uses stable v-completeness from \cite[Proposition 9.9]{KimIgusa}.

\item[(2)] For any affinoid perfectoid \(S\to\MCX\), the restriction
    of a morphism of the associated v-bundles to
    \(S\times\Spd(\ZZ_p)\) (resp. \(S\times\Spd(\QQ_p)\)) comes from
    a unique vector-bundle morphism on \(\Yco(S)\) (resp. \(\Yoo(S)\)),
    by (1) and the identification of diamonds in
    \cite[Proposition II.1.2]{FS21}. For a map \(T\to S\) over \(\MCX\),
    uniqueness in (1), applied to \(T\), makes these morphisms compatible
    with pullback. By \Cref{defYBundle}, they therefore define a unique
    morphism of integral (resp. rational) \(\MCY\)-bundles on \(\MCX\).

    \item[(3)] is clear for \(X=\Spa(B,B^+)\) affinoid perfectoid as the v-site over \(\Spa(B,B^+)\) has \(\Spa(B,B^+)\) itself as a final object and our \(\MCY\)-bundles are assumed to be compatible with pullback. The statement for general perfectoid spaces \(X\) follows from analytic gluing.
\end{enumerate}
\end{proof}
\subsection{\texorpdfstring{Exact structure on \(\MCY\)-bundles}{Exact structure on Y-bundles}}
We equip vector bundles on an adic space \(X\) with the exact structure given by short exact sequences of the underlying \(\MCO_X\)-modules. Thus a complex
\[0\longrightarrow\MCE_1\xrightarrow{\alpha_1}\MCE_2\xrightarrow{\alpha_2}\MCE_3\longrightarrow0,\qquad \alpha_2\alpha_1=0,\]
is short exact if \(\alpha_1\) is injective, \(\alpha_2\) is surjective as a morphism of sheaves, and \(\im(\alpha_1)=\ker(\alpha_2)\). Such sequences are locally split because \(\MCE_3\) is finite locally free, and hence remain short exact under arbitrary pullback.

\begin{mylemma}\emph{\cite[Proposition 7.51]{Wed19}}\label{adicPointsFromClassicalPoints}
    Let \((B,B^+)\) be a complete Huber pair and let \(\mfm\subset B\) be a maximal ideal. Then \(\mfm\) is closed and there exists \(v\in \Spa(B,B^+)\) with \(\supp v = \mfm\).
\end{mylemma}

\begin{mylemma}\label{exactnessAdicLocal}
Let \(X\) be a (not necessarily analytic) adic space and
\[0\longrightarrow\MCE_1\longrightarrow\MCE_2\longrightarrow\MCE_3\longrightarrow0\]
a complex of finite locally free \(\MCO_X\)-modules. Then it is short exact if and only if for every \(x\in X\) its base change to the completed residue field \(\hat\kappa(x)\) is exact.
\end{mylemma}

\begin{proof}
    A short exact sequence of vector bundles is locally split, so its base change to every completed residue field is exact. For the converse, fix \(x\in X\), put \(A=\MCO_{X,x}\), and let \(k_x=\MCO_{X,x}/\mfm_x\) be the (uncompleted) residue field. Now the natural map \(A\to\hat\kappa(x)\) induces a field extension \(k_x\to\hat\kappa(x)\). By assumption, the base change of the complex of finite projective \(A\)-modules
\[0\longrightarrow P_1\longrightarrow P_2\longrightarrow P_3\longrightarrow0\tag{*}\label{exactnessStalkSequence}\]
defined by \(P_j\coloneqq\MCE_{j,x}\) is exact after tensoring with \(\hat\kappa(x)\), and, as field extensions are faith\-fully-flat, it is exact already after tensoring with \(k_x\). From here, the proof proceeds like the one of \Cref{exactnessArcLocally}, applying Nakayama and splitting twice to obtain that the sequence of stalks \eqref{exactnessStalkSequence} is already exact. Since \(x\) was arbitrary, this shows that the sequence of sheaves is exact.
\end{proof}

\begin{mycor}\label{exactnessAdicLocalCorollary}
    Let \(X'\to X\) be a map of adic spaces such that \(|X'|\to |X|\) is surjective. Then a complex of finite locally free \(\MCO_X\)-modules
    \[0\longrightarrow\MCE_1\longrightarrow\MCE_2\longrightarrow\MCE_3\longrightarrow0\]
    is short exact if and only if its pullback to \(X'\) is.
\end{mycor}

\begin{mylemma}\label{YcoSurjectiveness}
    Let \(\Spa(A,A^+)\to \Spa(B,B^+)\) be a v-cover of characteristic \(p\) affinoid perfectoids\footnote{By definition, this is equivalent to being surjective on topological spaces, see \cite[Definition 17.1.1]{SW20}}. Then for \(\MCS\in \{\ZZ_p,\QQ_p\}\), the map
    \[|\Spa(A,A^+)\dot\times \Spa(\MCS)|\longrightarrow |\Spa(B,B^+)\dot\times \Spa(\MCS)|\]
    is surjective.
\end{mylemma}

\begin{proof}
We have a diagram:
\begin{center}
    \begin{tikzcd}
{|\Spa(A,A^+)\dot\times \Spa(\MCS)|} \arrow[r] \arrow[d, "\thicksim"] & {|\Spa(B,B^+)\dot\times \Spa(\MCS)|} \arrow[d, "\thicksim"] \\
{|\Spd(A,A^+)\times \Spd(\MCS)|} \arrow[r]                                       & {|\Spd(B,B^+)\times \Spd(\MCS)|.}                            
\end{tikzcd}
\end{center}
The vertical maps are isomorphisms by \cite[Lemma 15.6]{etCohDiamonds} as the spaces on the top row are analytic. The top map is the one we want to show surjectivity of. Thus we have to show that the bottom map is surjective.\par
Some \(x\in |\Spd(B,B^+)\times \Spd(\MCS)|\) is represented by a map \((B,B^+)\to (C,C^+)\) to a characteristic \(p\) perfectoid field together with an untilt \((C^\sharp,C^{\sharp +})\) (assumed to have characteristic \(0\) if \(\MCS=\QQ_p\)). By surjectivity of
    \[|\Spa(A,A^+)|\longrightarrow |\Spa(B,B^+)|,\]
    we can lift this to a map \((A,A^+)\to (\tilde C,\tilde C^+)\) to some perfectoid field \((\tilde C, \tilde C^+)\). Now the tilting equivalence for perfectoid \(C\)-algebras (\cite[Theorem 6.2.7]{SW20}) produces an untilt \((\tilde C^\sharp, \tilde C^{+\sharp})\) representing a preimage of \(x\) in \(|\Spa(A,A^+)\times \Spd(\MCS)|\).
\end{proof}

\begin{mycor}\label{exactnessVLocalBetweenRepresentables}
        Let \(\Spd(A,A^+)\to \Spd(B,B^+)\) be a v-cover of representables. Then a complex \(0\to\MCE_1\to\MCE_2\to\MCE_3\to0\) of vector bundles on \(\Yco(B,B^+)\) (resp. \(\Yoo(B,B^+)\)) is short exact if and only if its pullback to \(\Yco(A,A^+)\) (resp. \(\Yoo(A,A^+)\)) is.
\end{mycor}

\begin{proof}
Combine \Cref{exactnessAdicLocalCorollary} and \Cref{YcoSurjectiveness}.
\end{proof}

\begin{mydef}\label{defExactnessOfVBundles}
A complex \(0\to\MCE_1\to\MCE_2\to\MCE_3\to0\) of integral (resp. rational) \(\MCY\)-bundles on a small v-sheaf \(\MCX\) is called short exact if there is a v-cover by representables
\[(\Spa(B_i,B^+_i))_{i\in I}\longrightarrow \MCX\]
such that for any \(i\in I\) the evaluation 
    \[0\longrightarrow\MCE_1(\Spa(B_i,B^+_i))\longrightarrow\MCE_2(\Spa(B_i,B^+_i))\longrightarrow\MCE_3(\Spa(B_i,B^+_i))\longrightarrow0\]
    is short exact as a sequence of vector bundles on \(\Yco(B_i,B_i^+)\) (resp. \(\Yoo(B_i,B_i^+)\)).

\end{mydef}

These short exact sequences define the exact structure on integral/rational \(\MCY\)-bundles. They are stable under pullbacks along maps of small v-sheaves, since the evaluated short exact sequences of vector bundles are locally split. We will use this freely in the proofs below.

\begin{mylemma}\label{exactnessVBundleYBundleEquivalence}
     For a representable v-sheaf \(\Spd(B,B^+)\), a complex \(0\to\MCE_1\to\MCE_2\to\MCE_3\to0\) of vector bundles on \(\Yco(B,B^+)\) (resp. \(\Yoo(B,B^+)\)) is short exact if and only if it is short exact as a sequence of integral (resp. rational) \(\MCY\)-bundles.
\end{mylemma}

\begin{proof}
The ``only if'' direction is clear. Assume \(\MCE_1\to\MCE_2\to\MCE_3\) is an exact sequence of integral/rational \(\MCY\)-bundles on \(\Spd(B,B^+)\) and let  
 \[\Spd(B_i,B_i^+)_{i\in I}\longrightarrow \Spd(B,B^+)\] be a cover by representables such that 
 \[\MCE_1(\Spd(B_i,B_i^+))\longrightarrow \MCE_2(\Spd(B_i,B_i^+))\longrightarrow \MCE_3(\Spd(B_i,B_i^+))\]
 is an exact sequence of vector bundles on \(\Yco(B_i,B_i^+)\) (resp. \(\Yoo(B_i,B_i^+)\)) for each \(i\in I\). Then, by quasi-compactness of \(\Spd(B,B^+)\) as a v-sheaf (see \cite[Chapter 15]{etCohDiamonds}), there is a finite \(I_{\fin}\subseteq I\) such that
 \[\Spd(B_i,B_i^+)_{i\in I_{\fin}}\longrightarrow \Spd(B,B^+)\]
 is still a cover. Now 
 \[\Spd(\tilde B,\tilde B^+)\coloneqq\Spd\left(\bigoplus_{i\in I_{\fin}} B_i, \bigoplus_{i\in I_{\fin}} B_i^+\right)=\bigsqcup_{i\in I_{\fin}} \Spd(B_i,B_i^+) \longrightarrow \Spd(B,B^+) \]
 is a cover of representables, and by assumption our exact sequence pulls back to an exact sequence of vector bundles on \(\Yco(\tilde B,\tilde B^+)\) (resp. \(\Yoo(\tilde B,\tilde B^+)\)) and thus we can conclude by \Cref{exactnessVLocalBetweenRepresentables}.
\end{proof}

\begin{mycor}\label{exactnessYBundleEquivalentDefinition}
    A complex \(0\to\MCE_1\to\MCE_2\to\MCE_3\to0\) of integral (resp. rational) \(\MCY\)-bundles on a small v-sheaf \(\MCX\) is short exact if and only if for any map \(\Spd(B,B^+)\to \MCX\) from a representable the evaluation
    \[0\longrightarrow\MCE_1|_{\Spa(B,B^+)}\longrightarrow\MCE_2|_{\Spa(B,B^+)}\longrightarrow\MCE_3|_{\Spa(B,B^+)}\longrightarrow0\]
    is short exact as a sequence of vector bundles on \(\Yco(B,B^+)\) (resp. \(\Yoo(B,B^+)\)).
\end{mycor}

\begin{proof}
    The if direction is clear. For the ``only if'' direction, observe that the sequence of integral/rational \(\MCY\)-bundles remains exact as a sequence of integral/rational \(\MCY\)-bundles when pulled back to \(\Spd(B,B^+)\), and thus by \Cref{exactnessVBundleYBundleEquivalence} also as a sequence of vector bundles on \(\Yco(B,B^+)\) (resp. \(\Yoo(B,B^+)\)).
\end{proof}

\begin{myprop}\label{exactnessIsVLocal}
    Let \(\MCX'\to \MCX\) be a v-cover of small v-sheaves. Then a complex
    \[0\longrightarrow\MCE_1\longrightarrow\MCE_2\longrightarrow\MCE_3\longrightarrow0\]
    of integral/rational \(\MCY\)-bundles on \(\MCX\) is short exact if and only if its pullback to \(\MCX'\) is.
\end{myprop}

\begin{proof}
 Assume \(\MCE_1\to\MCE_2\to\MCE_3\) is an exact sequence of integral/rational \(\MCY\)-bundles on \(\MCX\) and let  
 \[\Spd(B_i,B_i^+)_{i\in I}\longrightarrow \MCX\] be a cover by representables such that 
 \[\MCE_1(\Spd(B_i,B_i^+))\longrightarrow \MCE_2(\Spd(B_i,B_i^+))\longrightarrow \MCE_3(\Spd(B_i,B_i^+))\]
 is an exact sequence of vector bundles on \(\Yco(B_i,B_i^+)\) (resp. \(\Yoo(B_i,B_i^+)\)) for each \(i\in I\). Choose a cover
 \[\bigsqcup_j\Spd(B'_j,B'^+_j)\longrightarrow\left(\bigsqcup_i \Spd(B_i,B^+_i)\right)\times_{\MCX}^{} \MCX'\]
 by representables. We get a commutative diagram:
    \begin{center}
        \begin{tikzcd}
{\bigsqcup_j\Spd(B'_j, B'^+_j)} \arrow[d, two heads] \arrow[r] & {\bigsqcup_i\Spd(B_i,B^+_i)} \arrow[d, two heads] \\
{\MCX'} \arrow[r]                                                                                          & \MCX                                             
\end{tikzcd}
    \end{center} 
    Now the pullback of our exact sequence on \(\MCX\) to the top-right and hence also the top-left v-sheaf is exact (when viewed as vector bundles on \(\Yco(\cdot)\) resp. \(\Yoo(\cdot)\)) by assumption, and the left vertical map is a v-cover, thus by definition 
    \[\MCE_1|_{\MCX'}\longrightarrow \MCE_2|_{\MCX'}\longrightarrow \MCE_3|_{\MCX'}\]
    is exact.\par
    Conversely, suppose that the pullback of the sequence to \(\MCX'\) is exact. By \Cref{defExactnessOfVBundles}, there is a v-cover by representables
    \[\bigsqcup_i\Spd(B_i,B_i^+)\longrightarrow\MCX'\]
    on which the evaluated sequences of vector bundles on
    \[\Yco(B_i,B_i^+)\quad\text{(resp. }\Yoo(B_i,B_i^+)\text{)}\]
    are exact. Since composites of v-covers are v-covers, the composite
    \[\bigsqcup_i\Spd(B_i,B_i^+)\longrightarrow\MCX'\longrightarrow\MCX\]
    is again a v-cover of \(\MCX\) by representables. The evaluations of the original sequence along these composite maps agree with those of its pullback to \(\MCX'\), so they are exact. Hence the original sequence is exact by \Cref{defExactnessOfVBundles}.
    \end{proof}
\subsection{\texorpdfstring{\(W(\MCO^\flat)\)-bundles on the v-site}
  {W(O♭)-bundles on the v-site}}
In this section, we will define \(W(\MCO^\flat)\)-bundles, a third kind of locally free sheaf on the adic site, and cite the descent result \Cref{Kim}. This result will be extended in section 4.1 and used in section 4.2.

\begin{mydef}\label{WBundleDef}
    \begin{enumerate}
        \item For any \(r\in \NN\cup \{\infty\}\), define a sheaf of rings
        \[W_r(\MCO^\flat)\colon (B,B^+)\longmapsto W_r(B)\]
        on \(\Perf^{\emph v}\), where \(W_\infty(B)\coloneqq W(B)\).
        \item  A \emph{\(W_r(\MCO^\flat)\)-bundle} on a small v-sheaf \(\MCX\) is one of the following equivalent pieces of data:
        \begin{itemize}
            \item A locally finite free sheaf of \(W_r(\MCO^\flat)\)-modules on \(\Perf^\emph{v}_{/\MCX}\).
            \item A 3-truncated v-hypercover
            \[X''' \rrrarrows X''\rightrightarrows X'\longrightarrow \MCX\]
            by perfectoid spaces together with a descent datum of finite projective modules on 
            \[W_r(\Gamma(X''',\MCO_{X'''}))\;\substack{\leftarrow\\[-1em] \leftarrow \\[-1em] \leftarrow}\; 
            W_r(\Gamma(X'',\MCO_{X''}))\leftleftarrows W_r(\Gamma(X',\MCO_{X'})).\]
        \end{itemize}

    \end{enumerate}
\end{mydef}

To verify that the \(W_r(\MCO^\flat)\) are actually sheaves for the v-topology, note that 
\[\MCO^\flat\colon(B,B^+)\mapsto B=B^\flat\]
indeed defines a sheaf (see \Cref{defVSheaf}). Hence, by induction on \(r\) and using the short exact sequence
\[0\longrightarrow B=W_1(B)\xrightarrow{\cdot\, p^r} W_{r+1}(B)\longrightarrow W_{r}(B)\longrightarrow 0,\]
\(W_r(\MCO^\flat)\) for \(r<\infty\) is a sheaf; finally, because limits of sheaves are sheaves, \(W(\MCO^\flat)\) is a sheaf. This also implies that the two definitions of \(W_r(\MCO^\flat)\)-bundles given above are equivalent.

If we equip \(\MCX\) with a map to \(\Spd(\ZZ_p)\) via \(\MCX \to \Spd(\FF_p)\to \Spd(\ZZ_p)\), the \(W_1(\MCO^\flat)\)-bundles defined here are just the v-bundles from \Cref{defVBundle}. The only case of \(W_r(\MCO^\flat)\)-bundles for \(r>1\) that will be relevant for us is the following:

\begin{myrem}\label{WBundleSpecialCase}
    Let \(R^+\) be a perfect ring in characteristic \(p\) equipped with the discrete topology. Then a \(W_r(\MCO^\flat)\)-bundle on \(\Spa(R^+,R^+)\) is given by a descent datum of finite projective modules on
    \[ W_r(R^+((t^{1/p^\infty}))) \rightrightarrows W_r(R^+((t^{1/p^\infty},s^{1/p^\infty})))  \rrrarrows W_r(R^+((t^{1/p^\infty},s^{1/p^\infty},u^{1/p^\infty}))).\]
    Indeed, \(\Spa(R^+((t^{\pRoots})),R^+[[t^{\pRoots}]])\to \Spa(R^+)\) is a v-cover whose \v{C}ech nerve has global sections
    \[R^+((t^{1/p^\infty})) \rightrightarrows R^+((t^{1/p^\infty},s^{1/p^\infty})) \rrrarrows R^+((t^{1/p^\infty},s^{1/p^\infty},u^{1/p^\infty})),\]
    see \cite[Lemma 5.2]{Kim} and the proof of \cite[Corollary 5.3]{Kim}. The case for general \(r\) now follows by definition.
\end{myrem}

For a precise definition of \(R^+((t_0^{1/p^\infty},\dots, t_n^{1/p^\infty}))\), see \cite[Definition 2.4]{Kim} or \Cref{KimExtLaurentDef}.

\section{Equivalence of the two categories}
\subsection{Descent of Witt vector bundles}

\begin{mythm}\label{Kim} (\cite[Theorem 4.9]{Kim})
    Every descent datum of finite projective modules on
    \[R^+((t^{1/p^\infty})) \rightrightarrows R^+((t^{1/p^\infty},s^{1/p^\infty}))  \rrrarrows R^+((t^{1/p^\infty},s^{1/p^\infty},u^{1/p^\infty})).\]
    descends to a finite projective \(R^+\)-module
\end{mythm}
In particular, this implies that v-bundles and vector bundles on \(\Spd(R^+)\) are equivalent (\cite[Corollary 5.3]{Kim}).
In this subsection, we want to extend this result to more general \(W_r(\MCO^\flat)\)-bundles, including \(r=\infty\). We will need three lemmas containing quite a bit of calculations, but first, let's properly define the rings of interest:

\begin{mydef}\cite[Definition 2.4]{Kim}\label{KimExtLaurentDef}
    Define:
    \[R^+((T_1^{\pRoots},\dots, T_n^{\pRoots}))\coloneqq
    \left\{
    \sum_{\mathbf{d}=(d_1,\dots,d_n)\in I}a_{\mathbf{d}}^{}T_1^{d_1}\cdots T_n^{d_n}
    \;\middle|\;
    \substack{I\subset \ZZ[p^{-1}]^n, \; a_{\mathbf{d}}\in R^+, \\
    \#\{(d_1,\dots,d_n)\in I \;\mid\; \sum_{j=1}^n s_jd_j<C\}<\infty \\
    \text{  for all \(C\in \RR\) and \(s_1,\dots,s_n\in \RR_{>0}\)}
    }
    \right\}\]
\end{mydef}
Informally, the condition can be phrased as: Some \(a=\sum_{\mathbf{d}\in \ZZip^n}a_{\mathbf{d}}^{}T_1^{d_1}\cdots T_n^{d_n}\) is in \(R^+((T_1^{\pRoots},\dots, T_n^{\pRoots}))\) if and only if for any (oriented)
hyperplane \(H\subset \ZZip^n\) defined by \(s_1,\dots,s_n\in \RR_{>0},C\in \RR\) there are only finitely many non-zero summands \(a_{\mathbf{d}}^{}T_1^{d_1}\cdots T_n^{d_n}\) for which \(\mathbf d\in \ZZip^n\) lies on the negative side of \(H\).\par

Let us rephrase the definition of \(R^+((T_1^{\pRoots},\dots,T_n^{\pRoots}))\) in a more convoluted way which will be useful for the computation of \Cech cohomology in \Cref{KimExtCechCalculation}. The guiding idea is to find an analogue of the decomposition (as \(R^+\)-modules) 
\[R^+[T_1,\dots,T_n]=\bigoplus_{d\in \ZZgez} R^+[T_1,\dots,T_n]_d,\]
where \(R^+[T_1,\dots,T_n]_d\) are the homogeneous polynomials of degree \(d\).

\begin{mylemma}\label{KimExtLaurentLemma}
    Let \(P\) be the set of partial functions \(\ZZip\to \ZZip\), i.e., pairs \((D,m)\) where \(D\subset \ZZip\) and \(m\colon D\to \ZZip, \; d\mapsto m_d\). We will generally say \(D=D_m\) and regard it as part of the datum of \(m\). Let \(S\) be the subset
    \[\left\{m\in P
    \;\middle|\;\substack{m_d\leq \min(0, \frac dn ) \text{ for each \(d\in D_m\)},\\
    \#\left\{d\in D_m \;\mid\;  \sum_{i=1}^n s_i\left(m_d + \delta_{ij}(d-nm_d)\right)<C \text{ for some } j\in \{1,\dots,n\}\right\} \\
    \text{ is finite for all \(C\in \RR\), \(s_1,\dots,s_n\in \RR_{>0}.\)}
    }
    \right\}.\]
    We define a partial order on \(P\) (and hence on \(S\)) by 
    \[m\leq m' \quad \iff \quad D_m\subseteq D_{m'} \emph{ and } \forall d\in D_m\colon m_d\geq m_d'.\]
    Note the change of directions in the inequalities. Then there is an isomorphism of \(R^+\)-modules:
    \[R^+((T_1^{\pRoots},\dots, T_n^{\pRoots}))=\colim_{m\in S} \;\;\prod_{d\in D_m}
    \;\bigoplus_{\substack{d_1,\dots,d_n\in \ZZ[p^{-1}]\\
    \forall i : d_i\geq m_d\\
    d_1+\cdots+d_n=d}}
    R^+ T_1^{d_1}\cdots T_n^{d_n}\]
\end{mylemma}

\begin{proof}
    First, let's remark on the definition of \(S\): For any pair \((d,m_d)\) the set of tuples \((d_1,\dots,d_n)\) over which we take the direct sum forms a scaled and shifted \(n-1\)-standard-simplex lying in the (\(\sum d_i = d\))-plane. Indeed, at least if \(m_d\leq \frac dn\) (which we enforce in the definition of S), it is the convex hull of the scaled unit vectors 
    \[(d-nm_d)e_1,\dots, (d-nm_d)e_n\in \ZZip^n,\]
    shifted by \(m_d\sum_{i=1}^n e_i\). So through this lens, some \(m\in P\) can be seen as a set of disjoint shifted standard simplices in \(\ZZip^n\), indexed by \(d\in D_m\). The partial order is defined so that \(m\leq m'\) if and only if the set of shifted simplices defined by \(m\) is contained in the set of shifted simplices defined by \(m'\) (as subsets of \(\ZZip^n\)).\par
    Now if any point \((d_1,\dots,d_n)\in \ZZip^n\) of a shifted simplex defined by \((d, m_d)\) lies on the negative side of a hyperplane \(H\) specified by \(s_1,\dots,s_n,C\), then by convexity the same must be already true for one of its \(n\) corner points 
    \[\sum_i(m_d + \delta_{ij}(d-nm_d))e_i.\] 
    Thus for an \(m\in P\) to be in \(S\), for any hyperplane \(H\) defined by \(s_1,\dots,s_n,C\), only finitely many simplices may have a corner point (or, equivalently, any point) lying on the negative side of \(H\). The condition \(m_d\leq \frac 1n d\) guarantees that the simplices are non-empty, while \(m_d\leq 0\) enforces that the (shifted) simplex is large enough not to be fully contained in the \(d_i>0\)-quadrant, a technical assumption that will be useful in \Cref{KimExtCechCalculation} - this result would also be valid without it.\par
    Now the two descriptions of \(R^+((T^{\pRoots}_1,\dots,T_n^{\pRoots}))\) that we want to compare in this lemma clearly lie in the set of \emph{all} formal sums
    \[\sum_{\mathbf{d}=(d_1,\dots,d_n)\in \ZZip^n}a_{\mathbf{d}}^{}T_1^{d_1}\cdots T_n^{d_n},\]
    so we can show double inclusion. We will use the shorthand \(|\mathbf d|\coloneqq \sum_{i=1}^n d_i\).
    \begin{itemize}
        \item["\(\supseteq\)":]  Every element \(a\) in the group on the right side of the isomorphism lies in
        \[\prod_{d\in D_m}
        \;\bigoplus_{\substack{d_1,\dots,d_n\in \ZZ[p^{-1}]\\
        d_1,\dots,d_n\geq m_d\\
        |\mathbf d|=d}}
        R^+ T_1^{d_1}\cdots T_n^{d_n}\]
        for some \(m\in S\). Let \(H\) be an oriented hyperplane defined by \(s_1,\dots,s_n,C\). Then by the discussion about \(S\) above only finitely many \(d\in D_m\) define shifted simplices having at least one point lying on the negative side of \(H\), and any such \((d,m_d)\) can contribute only finitely many summands, so in total, \(a\) has only finitely many non-zero summands \(a_{\mathbf{d}}^{}T_1^{d_1}\cdots T_n^{d_n}\) such that \(\mathbf{d}\) lies on the negative side of \(H\). Thus \(a\in R^+((T_1^{\pRoots},\dots,T_n^{\pRoots}))\).
        \item["\(\subseteq\)":] Assume 
        \[a=\sum_{\mathbf{d}\in \ZZip^n} a_{\mathbf{d}}^{}T_1^{d_1}\cdots T_n^{d_n}\in R^+((T_1^{\pRoots},\dots,T_n^{\pRoots})).\]
        We define 
        \[D_m\coloneqq\{|\mathbf d|\;:\; a_{\mathbf d}\neq 0\}\]
        and for each \(d\in D_m\)
        \[m_d\coloneqq \min\Big(0,\;
        \min_{\substack{i\in\{1,\dots,n\}\\\mathbf d\in \ZZip^n\\|\mathbf d|=d,\;a_{\mathbf d}\neq 0}} d_i\Big).\]
        This automatically enforces \(m_d\leq 0\) and \(m_d\leq \frac dn\). The inner minimum is actually taken over a finite set: Indeed, by the definition of \(R^+((T_1^{\pRoots},\dots, T_n^{\pRoots}))\), there are only finitely many non-zero summands \(a_{\mathbf d}T_1^{d_1}\cdots T_n^{d_n}\) for which \(|\mathbf d|<d+1\) (by choosing \(s_1=\cdots=s_n=1, C= d+1\)). In particular, \(m_d\) is well-defined. This finiteness also shows that
        \[a\in\prod_{d\in D_m}
        \;\;\bigoplus_{\substack{d_1,\dots,d_n\in \ZZ[p^{-1}]\\
        d_i\geq m_d\\
        |\mathbf d|=d}}
        R^+ T_1^{d_1}\cdots T_n^{d_n}.\]
        Now all that is left to show is \(m\in S\)\footnote{The relevant picture here is the following: A priori it could happen that for some \(a\in R^+((T_1^{\pRoots},\dots, T_n^{\pRoots}))\) the subset \(\{\mathbf d \mid a_{\mathbf d}\neq 0\}\subset \ZZip^n\) has only finitely many points on the negative side of every oriented hyperplane defined as above, but that once we take the ``simplicial hull" (as in, the set of simplices defined by \(m\colon D_m\to \RR\)) of \(\{\mathbf d \mid a_{\mathbf d}\neq 0\}\), those simplices have infinitely many corners on the negative side of some hyperplane. We have to show that this can, in fact, not happen}. Assume \(m\notin S\), i.e., that there is some choice of a hyperplane defined by \(C\in \RR,s_1,\dots,s_n\in \RR_{>0}\) and an infinite subset \(D_<\subseteq D_m\) such that for each \(d\in D_<\) there is some \(j\in \{1,\dots,n\}\) for which
        \begin{align*}
            &\sum_{i=1}^n s_i\left(m_d + \delta_{ij}(d-nm_d)\right)\\
        =&m_d\left(\sum_{i=1}^n s_i \right)+ (d-nm_d)s_{j}
        <C\tag{*}\label{KimExtSupportBound}.
        \end{align*}
        First, as there are only finitely many choices for \(j\), there is some fixed \(j_0\) and some smaller but still infinite \(D_{<}'\subseteq D_{<}\) such that \eqref{KimExtSupportBound} holds true for \(j=j_0\) and each \(d\in D_{<}'\). Also, from the second form of the expression \eqref{KimExtSupportBound}, we see that the condition only depends on \(\sum_{i=1}^n s_i\) and \(s_{j_0}\), and that since \(d-nm_d>0\), decreasing \(s_{j_0}\) without altering \(\sum_{i=1}^n s_i\) will only make the expression smaller, so we can change the hyperplane so that \(s_{j_0}<\frac 1n\sum_{i=1}^n s_i\) and \eqref{KimExtSupportBound} still holds true.
        Rearranging the equation a different way
        \begin{align*}&\sum_{i=1}^n s_i\left(m_d + \delta_{ij_0}(d-nm_d)\right)\\=
        &s_{j_0}d+
        \left(-ns_{j_0}+\sum_{i=1}^ns_i\right)m_d<C
        \end{align*}
        we get
        \[m_d < -cd+b\]
        where
        \[c=s_{j_0}\left(-ns_{j_0}+\sum_{i=1}^ns_i\right)^{-1}, \qquad b=C\left(-ns_{j_0}+\sum_{i=1}^ns_i\right)^{-1},\]
        using that, by our assumption on \(s_{j_0}\) we have \(-ns_{j_0}+\sum_{i=1}^ns_i>0\) so we can divide by it without changing the direction of the inequality. The same assumption also shows \(c>0\). \par
        We have already used above that for any bound \(C'\) there are only finitely many \(\mathbf d\in \ZZip^n\) for which \(a_{\mathbf d} \neq 0\) and \(|\mathbf d|<C'\). Now for all \(d\geq C'\coloneqq\frac{b}{c}\), we have \(m_d < -cd+b \leq 0\), so only the (finitely many) \(d<C'\) can have \(m_d=0\). Thus, for each \(d\) in the infinite set \(D''_<\coloneqq D'_<\cap [C',\infty)\) we have 
        \[m_d=\min_{\substack{\mathbf d\in \ZZip^n\\|\mathbf d|=d,\;a_{\mathbf d}\neq 0\\i\in\{1,\dots,n\}}} d_i.\]
        Let \(\mathbf D''_< \subset \{\mathbf{d}\in\ZZip^n\mid a_\mathbf{d} \neq 0\}\) be the set of \(\mathbf d=(d_1,\dots, d_n)\) where that minimum is attained, i.e. 
        \[d=|\mathbf d|\in D''_<, \quad m_d = \min_{i\in \{1,\dots, n\}} d_i.\]
        Because there are only finitely many choices for \(i\) in the minimum above, we can assume that there is some \(i_0\) and some infinite subset \(\mathbf D'''_<\subseteq \mathbf D''_<\) such that for each \(\mathbf d\in \mathbf D'''_<\), the minimum is attained at \(i_0\), i.e. 
        \[d_{i_0}=m_d < -cd+b.\]
        Setting
        \[s'_i=\begin{cases}1 \quad &\text{ if } i\neq i_0\\1+\frac 1c&\text{ if } i=i_0\end{cases}, \qquad C'=\frac bc\]
        we get 
        \begin{align*}
            \sum_{i=1}^n s'_id_i
            &=\left(\sum_{i=1}^n d_i \right)+ \frac{d_{i_0}}c\\
            &< d+\frac{(-cd+b)}c\\
            &=\frac b c \\
            &=C'
        \end{align*}
        for each of the infinitely many \(\mathbf d\in \mathbf{D}'''_<\), each of which has \(a_{\mathbf d}\neq 0\). This shows that \(a\notin R^+((T_1^{\pRoots},\dots, T_n^{\pRoots}))\).
    \end{itemize}
\end{proof}

\begin{mylemma}\label{KimExtCechCalculation}
    The \Cech complex
    \[C_\infty\colon \quad0 \longrightarrow R^+((T_0^{\pRoots}))\longrightarrow  R^+((T_0^{\pRoots}, T_1^{\pRoots}))\longrightarrow\cdots\]
    with differentials defined by
    \begin{align*}\tag{**}\label{KimExtDifferential}
        T_0^{d_0}\cdots T_n^{d_n}\mapsto \sum_{i=0}^{n+1} (-1)^i T_0^{d_0}\cdots T_{i-1}^{d_{i-1}}\cdot T_{i+1}^{d_i}\cdots T_{n+1}^{d_n}
    \end{align*}
    is acyclic in positive degrees.
\end{mylemma}

\begin{proof}
    We will gradually build up the complex \(C_\infty\) from simpler ones. To begin with, the complex associated to 
    \[C_{\pol}\colon \quad 0\longrightarrow R^+ [T_0^{\pm \pRoots}] \rightrightarrows R^+[ T_0^{\pm \pRoots},T_1^{\pm \pRoots}]  \rrrarrows \cdots\]
    is acyclic in positive degrees since this is just the \Cech nerve of \(\Spec(R^+[T^{\pm 1/p^\infty}])\to \Spec(R^+)\) which is an affine faithfully flat cover of an affine scheme. Note that as \(p=0\) in \(R^+\), no \(p\)-adic completion is necessary after adding the \(p\)-th power roots. The complex associated to this is
    \[0\longrightarrow R^+ [T_0^{\pm\pRoots}] \xrightarrow{d_1} R^+[T_0^{\pm\pRoots},T_1^{\pm\pRoots}] \xrightarrow{d_2}\cdots\]
    with differentials defined by the same formula \eqref{KimExtDifferential}. This obviously respects the degree of polynomials, so the complex decomposes as
    \[C_{\pol}=\bigoplus_{d\in \ZZip} C_{d}\]
    where
    \[C_{d}\colon\quad\left(0\longrightarrow R^+ [T_0^{\pm\pRoots}]_d \longrightarrow R^+[T_0^{\pm\pRoots},T_1^{\pm\pRoots}]_d \longrightarrow\cdots\right),\]
    with 
    \[R^+[T_0^{\pm\pRoots},\dots,T_n^{\pm\pRoots}]_d\subset R^+[T_0^{\pm\pRoots},\dots,T_n^{\pm\pRoots}]\]
    the homogeneous (fractional) polynomials of degree \(d\). Similarly, each of these further decomposes as 
    \[C_{d}= C_ {d,\geq m} \oplus C_{d,\ngeq m}\]
    where  \(m\in \ZZip_{\leq 0}\) and
    \begin{align*}
        &C_{d,\geq m}\colon \quad \left(0\longrightarrow R^+[T_0^{\pm\pRoots}]_{d,\geq m}\longrightarrow R^+[T_0^{\pm\pRoots}, T_1^{\pm\pRoots}]_{d,\geq m}\longrightarrow\cdots\right) \\\\
        &C_{d,\ngeq m}\colon \quad\left(0\longrightarrow R^+[T_0^{\pm\pRoots}]_{d,\ngeq m}\longrightarrow R^+[T_0^{\pm\pRoots}, T_1^{\pm\pRoots}]_{d,\ngeq m}\longrightarrow\cdots\right) 
    \end{align*}
    with
    \begin{align*}
        R^+[T_0^{\pm\pRoots},\dots,T_n^{\pm\pRoots}]_{d,\geq m} 
        &= 
        \bigoplus_{\substack{d_0,\dots,d_n\in \ZZ[p^{-1}]\\
        d_i\geq m\\
        |\mathbf d|=d}}
        R^+ T_0^{d_0}\cdots T_n^{d_n}.
        \\
        R^+[T_0^{\pm\pRoots},\dots,T_n^{\pm\pRoots}]_{d,\ngeq m} 
        &= 
        \bigoplus_{\substack{d_0,\dots,d_n\in \ZZ[p^{-1}]\\
        \exists i:d_i<m\\
        |\mathbf d|=d}}
        R^+ T_0^{d_0}\cdots T_n^{d_n}.
    \end{align*}
    Note that this would be wrong for \(m>0\) as the insertion of a \(1\) in the \(i\)-th position in the formula of the differential would not preserve 
    \[R^+[T_0^{\pm\pRoots},\dots,T_n^{\pm\pRoots}]_{d,\geq m} \longrightarrow R^+[T_0^{\pm\pRoots},\dots,T_{n+1}^{\pm\pRoots}]_{d,\geq m}.\]
    All of this is to say that \(C_{d,\geq m}\) is acyclic in positive degrees. Now, by \Cref{KimExtLaurentLemma},
    \[R^+((T_0^{\pRoots},\dots, T_n^{\pRoots}))=\colim_{m\in S} \;\;\prod_{d\in D_m}
    R^+[T_0^{\pm\pRoots},\dots,T_n^{\pm\pRoots}]_{d,\geq m_d},\]
    and as the differentials are still computed by the same formula, one automatically gets
    \[C_{\infty}=\colim_{m\in S}\prod_{d\in D_m} C_{d,\geq m_d}.\]
    Since products and filtered colimits of abelian groups are exact, \(C_{\infty}\) is acyclic in positive degrees, showing the statement.
\end{proof}

\begin{mylemma}\label{KimExtShortLemma}
    Let \(R^+\) be as above. Define
    \[Y\coloneqq\Spa(R^+((T^{1/p^\infty})),R^+[[T^{1/p^\infty}]])\longrightarrow X\coloneqq \Spa(R^+)\]
    and let \(Y^{(n)}\coloneqq Y\times_X\dots\times_X Y\) be the \(n\)-fold product. Then
    \[H^i_{\Perf^v}(Y^{(n)},\MCO_{Y^{(n)}})=0\]
    for all \(i,n>0\).
\end{mylemma}

\begin{proof}
\(Y^{(1)}=Y\) is affinoid perfectoid, so the case \(n=1\) follows from \cite[Theorem 17.1.3]{SW20}. For \(n>1\), the argument is similar to that at the beginning of Section II.2 in \cite{FS21}: We have a cover
\[Y^{(n)}=\cup_{m\in \NN} \Spa(B_m,B^+_m)\]
with
\begin{align*}
B_m^+&\coloneqq R^+[[T_1^{1/p^\infty},\dots,T_n^{1/p^\infty}]]\langle \frac{T_1^m}{\hat T_1\cdots T_n},\cdots, \frac{T_n^m}{T_1\cdots \hat T_n}\rangle,\\
B_m&=B_m^+[1/(T_1\cdots T_n)].
\end{align*}
As noted, the \(\Spa(B_m,B^+_m)\) have vanishing higher cohomology. The maps \(B_{m+1}\to B_{m}\) are not surjective, but they have dense image, so by \cite[(0.13.2.4)]{EGA3} \(R\lim_m B_m\) is concentrated in degree 0. Combining these facts yields the result.
\end{proof}

\begin{mythm}\label{KimExtTheorem}
Let \(R^+\) be a perfect ring of characteristic \(p\) carrying the discrete topology and let \(r\in \ZZ_{>0}\cup \{\infty\}\).
\begin{enumerate}
\item Let \(\MCE\) be a finite projective \(W_r(R^+)\)-module and \(\MCE_v\) the associated \(W_r(\MCO^\flat)\)-bundle on \(\Spd(R^+,R^+)\). Then
\[R\Gamma(\Perf_{/\Spd(R^+)}^{v},\MCE_v)=\MCE.\]
\item The pullback functor
\[\emph{finite projective \(W_r(R^+)\)-modules} \longrightarrow \emph{\(W_r(\MCO^\flat)\)-bundles on \(\Spa(R^+)\)}\]
is an equivalence of categories.
\end{enumerate}
\end{mythm}

\begin{proof}
\begin{enumerate}
    \item The heart of the proof is the case \(r=1, \MCE = R^+\), which can be computed using \Cech cohomology. Let \(X=\Spd(R^+)\) which we identify with the representing adic space \(X=\Spa(R^+)\).
    As \Cech cohomology requires a single object as input and \(\Perf^v_{/X}\) has no final object, we consider the larger site \(\Sh(\Perf_{/X}^v)\) of small v-sheaves over \(X\) with covers given by surjective maps of sheaves. As every small v-sheaf has a cover by a family of affinoid perfectoids, 
    \[H^i(\Perf^v_{/X},\MCO^\flat)=H^i_{\Sh(\Perf_{/X}^v)}(X,\MCO^\flat),\]
    where on the right side \(\MCO^\flat\) is extended to non-representable small v-sheaves in the natural way.\par
    Now let \(Y=\Spa(R^+((T^{1/p^\infty})), R^+[[T^{1/p^\infty}]])\) and consider the \Cech nerve
    \[\cdots Y\times_X Y \rightrightarrows Y \longrightarrow X.\]
    As in \Cref{WBundleSpecialCase},
    \[H^0(Y\times_X Y,\MCO_{Y\times_X Y})=R^+((T_0^{1/p^\infty},T_1^{1/p^\infty})),\]
    and similarly for higher fiber products. By \Cref{KimExtCechCalculation}, all higher \Cech cohomology groups associated to
        \[0\longrightarrow R^+((T_0^{1/p^\infty})) \rightrightarrows R^+((T_0^{1/p^\infty},T_1^{1/p^\infty}))  \rrrarrows \cdots \]
    vanish. Now by \Cref{KimExtShortLemma} the condition in \cite[03F7]{Stacks} is satisfied, which implies that \(H^i(\Perf^v_{/X},\MCO^\flat)=\check H^i_{\Sh(\Perf^v_{/X})}(Y\to X,\MCO^\flat)=0\) for \(i>0\).\par
    The case of \(\MCE=W_r(R^+)\) and general \(r<\infty\) follows by induction from the long exact cohomology sequence associated to
    \[0\longrightarrow R^+=W_1(R^+)\xrightarrow{\cdot p^r} W_{r+1}(R^+) \longrightarrow W_r(R^+)\longrightarrow 0\]
    and the case \(r=\infty\) by taking the limit (which is again exact as the transition maps are surjective).\par
    Finally, we can write any finite projective \(W_r(R^+)\)-module as a direct summand of a free one, showing the result for general \(\MCE\).
    \item Given \Cref{Kim} and the vanishing of cohomology in (1) this is fairly formal; a similar argument appears, e.g., in \cite[Theorem 4.1]{projectivity}. Full faithfulness follows from (1) applied to \(P^\vee\otimes Q\) for finite projective \(W_r(R^+)\)-modules \(P,Q\).\par
{\emergencystretch=1em For essential surjectivity at finite levels we do induction on \(r\). The case \(r=1\) is \Cref{Kim}. Now assume descent holds for some \(r\) and consider a \(W_{r+1}(\MCO^\flat)\)-bundle \(\MCE\). By induction hypothesis, \(\MCE/p^r\) descends to a finite projective \(W_r(R^+)\)-module \(M_r\). Lifting finite projective modules over nilpotent thickenings (\cite[0D4A]{Stacks}) produces a finite projective \(W_{r+1}(R^+)\)-module \(P\) with \(P/p^rP\cong M_r\), which, setting \(\MCF=P_v\) for the associated v-sheaf, translates to \(\overline\alpha:\MCF/p^r\xrightarrow{\sim}\MCE/p^r\). We have a short exact sequence\par}
    \[0\longrightarrow \MCE/p\xrightarrow{\cdot\, p^r} \MCE\longrightarrow \MCE/p^r\longrightarrow 0.\]
    Now applying \(\underline{\Hom}_{W_{r+1}(\MCO^\flat)}(\MCF,-)\) to this sequence preserves exactness, and by the equivalence for \(r=1\), \(\underline{\Hom}_{\MCO^\flat}(\MCF/p,\MCE/p)\) descends to a finite projective \(R^+\)-module. Now by (1), 
\[H^1(\Perf^v_{/X},\underline{\Hom}_{\MCO^\flat}(\MCF/p,\MCE/p))=0,\]
so the map
\[\underline{\Hom}_{W_{r+1}(\MCO^\flat)}(\MCF,\MCE) \longrightarrow \underline{\Hom}_{W_{r+1}(\MCO^\flat)}(\MCF,\MCE/p^r)\]
is also surjective on global sections and we can lift \(\overline\alpha\) to a global map \(\alpha:\MCF\to\MCE\), which is an isomorphism since it is so modulo the nilpotent ideal \((p)\). This shows descent for \(r+1\).\par
 For \(r=\infty\), first descend each finite step \(\MCE/p^r\) individually to some finite projective \(W_r(R^+)\)-module \(P_r\). Full faithfulness makes these compatible for varying \(r\), so \cite[0D4B]{Stacks} produces a finite projective module \(P=\varprojlim_rP_r\) with \(P/p^rP\cong P_r\), and 
\[P_v\cong\varprojlim_r(P_r)_v\cong\MCE\]
by \(p\)-adic completeness.
\end{enumerate}
\end{proof}
\subsection{\texorpdfstring{Formulation of the main theorem and proof of Case I: \(\varpi = p = 0\)}
  {Formulation of the main theorem and proof of Case I: ϖ=p=0}}
Let us formulate our main theorem:

\begin{myconst}\label{equivalenceConst}
Let \(R^+\) be a complete topological ring carrying the \(I\)-adic topology for some finitely generated ideal \(I\subset R^+\) containing \(p\). Then there is a functor 
\[\text{\(\OpP\)-bundles on } R^+
\longrightarrow\text{integral \(\MCY\)-bundles on }\Spd(R^+)\]
defined as follows: For any map from a representable \(\Spa(B,B^+)\to \Spd(R^+)\), we get an untilt \((B^\sharp,B^{+\sharp})\) with a map \(R^+ \to B^{+\sharp}\). As \(B^{+\sharp}\) itself is integral perfectoid, an \(\OpP\)-bundle defines a finite projective \(\Ainf(B^{+\sharp})=W(B^+)\)-module and hence a vector bundle on \(\Yco(B,B^+)\).
Upon varying \((B,B^+)\), this defines a \(\MCY\)-bundle on \(\Spd(R^+)\).
\end{myconst}

\begin{mythm}\label{equivalenceThm}
Let \(R^+\) be a complete topological ring carrying the \(I\)-adic topology for some finitely generated ideal \(I\subset R^+\) containing \(p\). Then the functor
\[\text{\emph{\(\OpP\)-bundles on }} R^+
\longrightarrow\text{\emph{integral \(\MCY\)-bundles on }}\Spd(R^+)\]
constructed above is an exact tensor equivalence. In particular, if \(R^+\) is integral perfectoid, finite projective \(\Ainf(R^+)\)-modules are equivalent to integral \(\MCY\)-bundles on \(\Spd(R^+)\) by \Cref{CrysForPerfdRings}.
\end{mythm}

\begin{myrem}\label{equivalenceTheoremRemark}
For \(R^+\) integral perfectoid, by \Cref{CrysForPerfdRings} \(\OpP\)-bundles on \(R^+\) are just finite projective \(\Ainf(R^+)\)-modules. In particular, the category of \(\OpP\)-bundles does not depend on the topology on \(R^+\) (i.e., a choice of \(I\subset R^+\)). Thus \Cref{equivalenceThm} implies that the category of integral \(\MCY\)-bundles on \(\Spd(R^+)\) is independent of the topology as well. This came as quite a surprise to me, and I do not know of a more direct way to see this. Also, I do not know whether this remains true for non-perfectoid \(R^+\) as in that case, a change of topology might alter the category of \(\OpP\)-bundles on \(R^+\) as well (see \Cref{defPrismaticSite}).
\end{myrem}

The proof of \Cref{equivalenceThm} will be completed in \Cref{secFinishing}. In this section, we want to prove (without the exactness) the case where \(R^+\) is a discrete perfect ring of characteristic \(p\) that is a product of valuation rings. We will do this by comparing integral \(\MCY\)-bundles to the \(W(\MCO^\flat)\)-bundles from Subsection 3.4, for which a similar descent result is already known. First, we will work out a general criterion for the comparison of global sections which we will use. The setup is the following:\par
For any ring \(R^+\), consider a level-wise injective map of 2-truncated cosimplicial rings over \(R^+\)
\begin{center}
\begin{tikzcd}
A^2 \arrow[r, "f_2", hook]                                                                 & B^2                                                                 \\
A^1 \arrow[u, "i^A_2"', shift right] \arrow[u, "i^A_1", shift left] \arrow[r, "f_1", hook] & B^1 \arrow[u, "i_2^B"', shift right] \arrow[u, "i^B_1", shift left]\\
R^+ \arrow[r, "\id"] \arrow[u] & R^+ \arrow[u]
\end{tikzcd}
\end{center}
and a 2-truncated descent datum of finite projective modules on \(A^1\rightrightarrows A^2\), i.e.,
\begin{itemize}
    \item A finite projective \(A^1\)-module \(M_A\).
    \item An isomorphism of \(A^2\)-modules
    \[\psi_{A}\colon M_A\otimes_{A^1, i_1^A} A^2\isom M_A\otimes_{A^1, i_2^A}A^2.\]
\end{itemize}
Let \((M_B, \psi_{B})\) be the base change of \((M_A, \psi_{A})\) to \(B^1\rightrightarrows B^2\). Denote by \(H^0(M_A,\psi_A)\) the elements \(x\in M_A\) which intertwine \(i_1^A, i_2^A\) and \(\psi_A\) in the obvious way, and similarly for \(H^0(M_B,\psi_B)\).

\begin{mylemma}\label{caseILemmaGeneralStatement}
Assume that in the setup above, the following properties hold: 

\begin{enumerate}
    \item The map \(R^+\to \eq(B^1\rightrightarrows B^2)\) is an isomorphism.
    \item \((M_B,\psi_B)\) is the trivial rank \(n\) descent datum, i.e., obtained as base change of \((R^+)^n\).
    \item There exists an integer \(m_0\geq n\) which is at least twice the size of a generating family of \(M_A\) over \(A^1\), such that for every \(C\in M_{n\times m_0}(B^1)\) and \(D\in\GL_{m_0}(A^2)\) satisfying
        \[i_2^B(C)=i^B_1(C)\cdot f_2(D),\]
    \(C\) already lies in the image of \(M_{n\times m_0}(A^1)\) under \(f_1\).
\end{enumerate}
Then the map
\[H^0(M_A,\psi_A)\longrightarrow H^0(M_B,\psi_B)\]
is an isomorphism.
\end{mylemma}

Note that assumption (1) (together with injectivity of \(f_1\) and \(f_2\)) implies injectivity of \(R^+\to A^1\) and \(R^+\to B^1\) as well as \(\eq(A^1\rightrightarrows A^2) =R^+\).

\begin{proof}
    As finite projective modules are flat, base change induces an injection 
    \[M_A\hookrightarrow M_A\otimes_{A^1}B^1=M_B,\] 
    by which we can understand \(M_A\) as an \(A^1\)-sublattice of \(M_B\). This already shows injectivity of
    \[H^0(M_A,\psi_A)\longrightarrow H^0(M_B,\psi_B).\]
    For surjectivity, we will show that every \(\psi_B\)-invariant section of \(M_B\) takes values in \(A^1\) under all elements of \(M_A^\vee\); by biduality, it then belongs to \(M_A\). Consider the dual datum
    \[L_A\coloneqq M_A^\vee,\qquad L_B\coloneqq M_B^\vee=L_A\otimes_{A^1}B^1,
    \qquad \theta_A\coloneqq(\psi_A^{-1})^\vee,\qquad \theta_B\coloneqq(\psi_B^{-1})^\vee.\]
For later reference, note that if \(M_A\) can be generated by \(m\) elements, the same is true for \(L_A\), as can be seen by splitting a surjection \((A^1)^m\to M_A\) and dualizing that splitting. Choose a presentation
    \begin{equation*}\tag{*}\label{caseIPresentation}
        (A^1)^m \overset{\pi}{\twoheadrightarrow} L_A\hookrightarrow L_B.
    \end{equation*}
    Consider the following diagram:
    \begin{center}
\begin{tikzcd}
{(A^1)^m\otimes_{A^1,i_1^A} A^2} \arrow[d, "\thicksim"] \arrow[r, two heads] & {L_A\otimes_{A^1,i_1^A} A^2} \arrow[r, hook] \arrow[d, "\theta_A"] & {L_B\otimes_{B^1,i_1^B} B^2} \arrow[d, "\theta_B"] \\
{(A^1)^m\otimes_{A^1,i_2^A} A^2} \arrow[r, two heads]           & {L_A\otimes_{A^1,i_2^A} A^2} \arrow[r, hook]                     & {L_B\otimes_{B^1,i_2^B} B^2}.
\end{tikzcd}
    \end{center}

The horizontal maps are the base changes of the composition \eqref{caseIPresentation}. To ensure that the left map can be chosen to be invertible (a not-necessarily invertible lift exists because the source \((A^1)^m\otimes_{A^1,i_1^A}A^2\cong(A^2)^m\) is projective as an \(A^2\)-module), we have to enlarge our free presentation of \(L_A\) by replacing \(\pi\) from \eqref{caseIPresentation} by
\[\widetilde\pi\colon (A^1)^{2m}=(A^1)^m\oplus(A^1)^m\twoheadrightarrow L_A,
\qquad (x,y)\mapsto\pi(x).\]
Now denote by \(N_k=L_A\otimes_{A^1,i_k^A}A^2\) the terms in the middle of the diagram and let \(\pi_k\colon\allowbreak (A^2)^m\twoheadrightarrow N_k\) be the base change of \(\pi\). The new presentations over \(A^2\) then become
\[\widetilde\pi_k\colon (A^2)^{2m}\twoheadrightarrow N_k,
\qquad (x,y)\mapsto\pi_k(x).\]
After identifying \(N_1\) with \(N_2\) via \(\theta_A\), we apply Schanuel's lemma (\cite[00O3]{Stacks}) to the two short exact sequences of \(A^2\)-modules
\[
\begin{aligned}
0&\longrightarrow\ker(\pi_1)\longrightarrow (A^2)^m\xrightarrow{\theta_A\circ\pi_1}N_2\longrightarrow 0,\\
0&\longrightarrow\ker(\pi_2)\longrightarrow (A^2)^m\xrightarrow{\pi_2}N_2\longrightarrow 0.
\end{aligned}
\]
Both sequences have the same quotient \(N_2\) and free middle terms \((A^2)^m\). The diagrammatic form of the lemma in loc. cit., interchanging the two free summands in its bottom row, produces an isomorphism
\[F\colon (A^2)^{2m}\isom (A^2)^{2m},
\qquad \widetilde\pi_2\circ F=\theta_A\circ\widetilde\pi_1\]
and thus the required invertible lift. Since \(2m\leq m_0\), we may without loss of generality replace \((A^2)^{2m}\) by \((A^2)^{m_0}\), sending the additional generators to 0 inside \(N_k\) and setting \(F\) as the identity on them. From now on, write \(m=m_0\).

Let \(e_1,\dots,e_n\) be the basis of \(M_B\) induced by the trivialization in assumption (2), and let \(e_1^\vee,\dots,e_n^\vee\) be its dual basis in \(L_B\). Denote by \(\lambda_i\in L_A\), \(1\leq i\leq m\), the images under \(\pi\) of the standard basis vectors of \((A^1)^m\). Define \(C\in M_{n\times m}(B^1)\) by
\[\lambda_i=\sum_{j=1}^n C_{ji}e_j^\vee\quad\text{in }L_B.\]
Thus \(C\) is the matrix of the base change of \eqref{caseIPresentation} to \(B^1\). Now equip:
\begin{itemize}
    \item \((A^1)^m\otimes_{A^1,i_k^A} A^2\) with the standard basis.
    \item \(L_B\otimes_{B^1,i^B_k} B^2\) with the basis given by the images of the \((e_i^\vee)_i\) under \(L_B\to L_B\otimes_{B^1,i_k^B}B^2\),
\end{itemize}
for \(k\in \{1,2\}\) respectively. Then we get the following coordinate descriptions of the maps:
\begin{itemize}
    \item By definition, the horizontal compositions in the diagram above are given as \(x\mapsto (i_k^B(C))x\) for \(k\in \{1,2\}\).
    \item The right vertical map is given as \(x\mapsto \mathbf{1}_nx\) because of assumption (2).
    \item The left vertical map is given as \(x\mapsto D^{-1}x\) for some \(D\in \GL_m(A^2)\).
\end{itemize}

Now commutativity of the diagram implies
\[i_2^B(C)=i^B_1(C)\cdot f_2(D),\]
so by assumption (3), \(C\in M_{n\times m}(A^1)\). By assumptions (1) and (2), any \(x\in H^0(M_B,\psi_B)\) has the form
\[x=\sum_{j=1}^n a_je_j,\qquad a_j\in R^+,\]
hence
\[\lambda_i(x)=\sum_{j=1}^n C_{ji}a_j\in A^1.\]
Since the \(\lambda_i\) generate \(L_A=M_A^\vee\), evaluation at \(x\) defines an element of \(L_A^\vee\cong M_A\) whose image in \(M_B\) is \(x\) and thus \(x\in M_A\). Its compatibility with \(\psi_A\) follows from its compatibility with \(\psi_B\), since the maps
\[M_A\otimes_{A^1,i_k^A}A^2\hookrightarrow M_B\otimes_{B^1,i_k^B}B^2\]
are injective by finite projectivity and injectivity of \(f_2\). This proves surjectivity, and hence bijectivity, of
\[H^0(M_A,\psi_A)\longrightarrow H^0(M_B,\psi_B),\]
concluding the proof.
\end{proof}

Now let's define the rings to which we want to apply this general criterion.
\begin{mydef}\label{caseIDefinition}
Let \(R^+\) be a perfect ring of characteristic \(p\). Define the following topological rings:
\begin{align*}
    R^{i}&\coloneqq R^+((T_1^{\pRoots},\dots,T_i^{\pRoots}))\\
    R^{i+}&\coloneqq R^+[[T_1^{\pRoots},\dots,T_i^{\pRoots}]]\\
\end{align*}
Let 
\begin{align*}
X^i&\coloneqq\Spa(R^1,R^{1+})\times_{\Spa(R^+)}\dots\times_{\Spa(R^+)}\Spa(R^1,R^{1+})\\
Y^i&\coloneqq X^i\dot \times \Spa(\ZZ_p).
\end{align*}
We want to take a closer look at \(Y^1\) and \(Y^2\). Let \(m\) be a positive integer.
\begin{enumerate}
    \item Let
    \[Y^1_{[0,m]}\coloneqq \MCY_{[0,m]}(R^1,R^{1+})\subset Y^1=\Yco(R^1,R^{1+}).\]
    Then 
    \[Y^1_{[0,m]}=\Spa(A^1_{[0,m]}, A^{1+}_{[0,m]})\]
    where
    \begin{align*}
        A^{1+}_{[0,m]}
    &= W(R^+[[T^{\pRoots}]])\Langle\frac {p^m} T\Rangle\\
    A^1_{[0,m]}&=A_{[0,m]}^{1+}[T^{-1}].
    \end{align*}

    \item One computes
        \[X^2=\bigcup_{m'\in \ZZ_{>0}}\Spa(R_{m'}^2, R_{m'}^{2+})\]
        where
        \[R_{m'}^2\coloneqq R^+[[T^{\pRoots},S^{\pRoots}]]\Langle\frac{T^{m'}}S,\frac{S^{m'}}{T}\Rangle [(TS)^{-1}].\]
        so that \(R^2=\bigcap_{m'} R^2_{m'}\). This fact was already used in the definition of \(W(\MCO^\flat)\)-bundles, see the proof of \cite[Lemma 5.2]{Kim}. Now
        \[Y^2=\bigcup_{m'\in \ZZgz} \Yco(R^2_{m'}, R^{2+}_{m'}),\]
        and, writing \(\operatorname{pr}_1,\operatorname{pr}_2\colon Y^2\rightrightarrows Y^1\) for the two projections, we define
    \[
    Y^2_{[0,m]}\coloneqq \operatorname{pr}_1^{-1}(Y^1_{[0,m]})\cap\operatorname{pr}_2^{-1}(Y^1_{[0,m]})\subset Y^2.
    \]
    This guarantees that both projections land in \(Y^1_{[0,m]}\). For varying \(m\), these domains form an increasing open cover of \(Y^2\).
    The ring \(A^2_{[0,m]}\coloneqq\Gamma(Y^2_{[0,m]},\MCO_{Y^2})\) has the description
    \[A^2_{[0,m]}
    =\bigcap_{m'\in \ZZgz} 
    W\Big(R^+[[T^{\pRoots},S^{\pRoots}]]\Big)
    \Langle\frac{T^{m'}}S,\frac{S^{m'}}{T},\frac {p^{m}}{T},\frac {p^{m}}{S}\Rangle
    [(TS)^{-1}]\]
\end{enumerate}
Finally, define
\[B^i\coloneqq W(R^i).\]
\end{mydef}

\begin{mylemma}\label{caseILemmaConvergenceCondition}
    \begin{enumerate} 
    \item We have injections of rings:
    \[ \dots \hookrightarrow A^1_{[0,3]}\hookrightarrow A^1_{[0,2]}\hookrightarrow A^1_{[0,1]}\hookrightarrow B^1\]
    and 
    \[ \dots \hookrightarrow A^2_{[0,3]}\hookrightarrow A^2_{[0,2]}\hookrightarrow A^2_{[0,1]}\hookrightarrow B^2.\]
    \item For some 
    \[a=\sum_{\substack{k\in \ZZip\\ n \in\ZZgez}} a_{n,k}T^kp^n\in B^1\]
    the following are equivalent:
    \begin{enumerate}[label=(\roman*)]
        \item \(a\in A^1_{[0,m]}\).
        \item For each \(M\in \RR\) there are only finitely many non-zero summands \(a_{n,k}T^kp^n\) for which \(mk+n<M\).
        \item The sequence 
        \[n \mapsto m(\min_{\substack{k\in \ZZip\\a_{n,k}\neq 0}} k) +n \]
        converges to \(\infty\).
    \end{enumerate}
        \item For some
        \[a=\sum_{\substack{k,l\in \ZZip\\ n \in\ZZgez}} a_{n,k,l}T^kS^lp^n\in B^2,\]
        consider the following statements:
        \begin{enumerate}[label=(\roman*)]
            \item \(a\in A^2_{[0,m]}\).
            \item For any \(M\in \RR\) there are only finitely many non-zero summands \(a_{n,k,l}T^kS^lp^n\) such that
            \[a_{n,k,l}\neq 0,\quad ml+mk+n <M.\]
            \item The sequence
            \[n \mapsto \min_{\substack{k,l\in \ZZip\\ a_{n,k,l}\neq 0}}(mk+ml+n)\]
            converges to \(\infty\).
        \end{enumerate}
        Then \((i)\implies (ii)\implies (iii)\).
    \end{enumerate}
    \end{mylemma}
    Here minima over empty supports are understood to be \(+\infty\).
    \begin{proof}
        \begin{enumerate}
        \item Is clear from the definitions. For the fact that \(A^2_{[0,m]}\) is contained in \(B^2\), note once again that
        \[R^+((T^{\pRoots}, S^{\pRoots}))=\bigcap_{m'\in \ZZgz}R^+[[T^{\pRoots},S^{\pRoots}]]\Langle \frac{T^{m'}}S, \frac {S^{m'}}{T}\Rangle[(TS)^{-1}].\]
        \item\begin{enumerate}[align=left]
                    \item[(ii) \(\implies\) (i)] Write \(n=mq+r\) with \(0\leq r<m\). Then
                    \[T^kp^n=T^{k+q}p^r\left(\frac{p^m}{T}\right)^q,\qquad k+q=\frac{mk+n-r}{m}.\]
                    Condition (ii) therefore implies that the series converges in the \(T\)-adic topology of the rational localization after inverting \(T\), so it defines an element of \(A^1_{[0,m]}\).
                \item[(i) \(\implies\) (ii)]
                Let \[D=W(R^+)[T^{\pRoots}][p^m/T],\] so that
                \[A^1_{[0,m]}=D^{\wedge(T)}[1/T].\]
                A monomial \(T^kp^n\) is divisible by \(T^r\) inside \(D\) if and only if \(k+\left\lfloor\frac{n}{m}\right\rfloor\geq r\). For any given \(r\), all but finitely many non-zero monomials in the \(T\)-adic completion \(D^{\wedge(T)}\) must satisfy this, and this stays true inside \(A^1_{[0,m]}=D^{\wedge(T)}[1/T]\). Multiplying the inequality by \(m\) and setting \(M=mr\) yields our result.
                    \item[(ii) \(\implies\) (iii)] This is clear.
                \item[(iii) \(\implies\) (ii)] First, by the convergence in (iii), there can only be finitely many \(n\in \NN\) for which there exists some \(k\in \ZZip\) such that \(a_{n,k}\neq 0\) and \(mk+n<M\). Second, by the definition \(B^1\coloneqq W(R^+((T^{\pRoots})))\), for each \(n\) there are only finitely many \(k\) such that \(a_{n,k}\neq 0\) and \(k<\frac{M-n}{m}\). Combining these two facts yields the result.
                \end{enumerate}
        \item\begin{enumerate}[align=left]
            \item[(i) \(\implies\) (ii)]
                This works similarly as before. By the definition of \(A^2_{[0,m]}\) as an intersection, it suffices to show the statement for an element inside a single
                \[E^{\wedge(T)}[1/(TS)], \quad \text{for} \quad E\coloneqq W(R^+)[T^{\pRoots}, S^{\pRoots}][\frac{p^m}{T},\frac{p^m}{S},\frac{T}{S},\frac{S}{T}].\]
                This corresponds to the case \(m'=1\) in the intersection. Now for any given \(n\), all but finitely many non-zero monomials \(T^lS^kp^r\) in the expansion of an element of \(E^{\wedge(T)}\) are divisible by \((ST)^n\) inside \(E\), which is equivalent to
                \[\lfloor l\rfloor+\lfloor k\rfloor+\left\lfloor\frac{r}{m}\right\rfloor\geq 2n,\]
                which stays true again in the localization. Setting \(M=2nm\) yields the result.
         \item[(ii) \(\implies\) (iii)] This is again clear.
            \end{enumerate}
        \end{enumerate}
        \end{proof}
\begin{mylemma}\label{caseILemmaAssumption3}
   The rings
    \begin{center}
    \begin{tikzcd}
    A^2_{[0,m]} \arrow[r, "f_2", hook]                                                                 & B^2                  \\
    A^1_{[0,m]} \arrow[u, "i^A_2"', shift right] \arrow[u, "i^A_1", shift left] \arrow[r, "f_1", hook] & B^1 \arrow[u, "i_2^B"', shift right] \arrow[u, "i^B_1", shift left]\\
    W(R^+) \arrow[r, "\id"] \arrow[u] & W(R^+) \arrow[u]
    \end{tikzcd}
    \end{center}
    from \Cref{caseIDefinition} satisfy condition (3) from \Cref{caseILemmaGeneralStatement}.
\end{mylemma}

\begin{proof}
        Consider an equation
        \[i_2^B(C)=i^B_1(C)f_2(D)\]
        as in \Cref{caseILemmaGeneralStatement} (3). For notational convenience, set \(C_1\coloneqq i_1^B(C), C_2\coloneqq i_2^B(C)\) and denote \(f_2(D)\) simply by \(D\), so the equation becomes \(C_2=C_1D.\) Write \(r\) for the number of columns of \(C\), so that \(D\in\GL_r(A^2_{[0,m]})\).
        Write the entries of \(C_1\), \(C_2\), and \(D\) as
        \[c^{ij}_1=\sum_{k\in \ZZip, n\in \ZZgez} c^{ij}_{n,k} T^k p^n\]
        \[c^{ij}_2=\sum_{l\in \ZZip, n\in \ZZgez} c^{ij}_{n,l} S^l p^n\]
        \[d_{ij}=\sum_{k,l\in \ZZip, n\in \ZZgez} d^{ij}_{n,k,l} T^k S^l p^n.\]
        Note that by the symmetry between \(C_1\) and \(C_2\), we can indeed use the same coefficients for \(c^{ij}_1\) and \(c^{ij}_2\).\par
        For \(N\in \ZZgez\), let
        \[L_N\coloneqq \min_{\substack {\exists i,j,l\\c^{ij}_{N,l}\neq 0}}\;l\]
        Set \(L_N=+\infty\) when the support at degree \(N\) is empty. Then, using \(C_2=C_1D\) and the coefficient-multiplication rule, for any \(N\) with \(L_N<+\infty\) there must be non-zero summands
        \[d^{i_Nj_N}_{n_N,k_N,L_N}T^{k_N}S^{L_N}p^{n_N}, \quad 
        c^{i'_Nj'_N}_{n'_N,k'_N} T^{k'_N} p^{n'_N}\]
        of entries of \(D\) and \(C_1\), respectively, such that
        \[ n'_N+n_N\leq N,\quad k'_N+k_N=0.\]
        The inequality allows for \(p\)-adic carries.
        Furthermore, by definition of \(L_N\) and the symmetries between \(C_1\) and \(C_2\)
        \[k'_N\geq L_{n'_N}\]
        for each \(N\). Thus we get
        \begin{align*}
            mL_N+mk_N+n_N
            &\leq mL_N-mk'_N+(N-n'_N)\\
            &\leq m(L_N-L_{n'_N})+(N-n'_N)\\
            &=mL_N+N - (mL_{n'_N}+n'_N)\\ 
            &=a_N-a_{n'_N}
        \end{align*}
        where \(a_N\coloneqq mL_N+N\). We want to show: if \((a_N)_N\) does \emph{not} converge to \(\infty\), so that \(C\notin M_{n\times r}(A^1_{[0,m]})\) by \Cref{caseILemmaConvergenceCondition} (2) (iii), then neither does \(a_N-a_{n'_N}\), considering only indices with \(a_N<+\infty\). This will imply that \(D\notin\GL_r(A^2_{[0,m]})\) by \Cref{caseILemmaConvergenceCondition} (3) (ii). This is an easy and fun calculus I exercise which the reader is invited to solve for themselves. We proceed in two steps: First, we may assume that the finite values of \(a_N\) are bounded below -- for if they weren't, there would be a subsequence \(a_{s(N)}\) of successively lower minima such that \(a_{s(N)}\leq a_{N'}\) for each \(N'\leq s(N)\), so \(a_{s(N)}-a_{n'_{s(N)}}\leq 0\), so \(a_N-a_{n'_N}\) has a bounded-above subsequence and can't converge to \(\infty\). So we can assume \(a_N\) is bounded below by some \(M\), and thus \(a_N-a_{n'_N}\leq a_N-M\), which does not converge to \(\infty\) if \(a_N\) does not. In both cases, along the chosen bounded-above subsequence of \(a_N\), we have \(L_N=(a_N-N)/m\to-\infty\). The selected terms of \(D\) therefore have infinitely many distinct \(S\)-exponents but additive weights bounded above, contradicting (3) (ii) for the finitely many entries of \(D\).
\end{proof}

\begin{mythm}\label{mainThmDiscrete}
The functor in \Cref{equivalenceConst} is an equivalence for \(R^+=A^+\) a discrete ring that is a product of perfect characteristic \(p\) valuation rings.
\end{mythm}

\begin{proof}
Consider the composition of functors
\begin{center}
\begin{tikzcd}
\text{finite locally free \(W(A^+)\)-modules} \arrow[d] \\
\text{\(\MCY\)-bundles on \(\Spd(A^+)\)} \arrow[d]      \\
\text{\(W(\MCO^\flat)\)-bundles on \(\Spd(A^+)\)}               
\end{tikzcd}
\end{center}
The first functor is the one defined in \Cref{equivalenceConst}, which we want to show is an equivalence, while the second one is defined as follows: For any map \((A^+,A^+)\to (B,B^+)\) to an affinoid perfectoid, the map of adic spaces
\[\MCY_{[0,1]}(B,B^+)\hookrightarrow \Yco(B,B^+),\]
produces by pullback a vector bundle on the left from every vector bundle on the right. But the left side is an affinoid analytic adic space, and so we get a finite projective module over the ring \(C\) of global sections on \(\MCY_{[0,1]}(B,B^+)\). Finally, there is a map of rings
\[C\longrightarrow W(B)\]
producing a finite projective module on \(W(B)\). The descent data can be handled similarly.\footnote{The reason why we have to argue through \(\MCY_{[0,1]} (A^+((t)))\) is that \(\Yco(A^+((t)))\) is not affinoid, while \(W(A^+((t)))\) does not seem to be a Huber ring, for similar reasons as in the discussion below Proposition 4.2.6 in \cite{SW20}. This problem is responsible for all the bookkeeping in this proof.} \par
Now to show that both arrows are equivalences it suffices to show that the composition is an equivalence and the second arrow is fully faithful. The former is \Cref{KimExtTheorem} (2), the latter will take the rest of this proof.\par
By the usual tensor product trick, it suffices to show that for any integral \(\MCY\)-bundle \(\MCE_{\MCY}\) on \(\Spd(A^+)\) which maps under the second functor to some \(W(\MCO^\flat)\)-module \(\MCE_{W(\MCO^\flat)}\) the induced map of sections
\[F\colon H^0(\Spd(A^+),\MCE_{\MCY})\longrightarrow H^0(\Spd(A^+), \MCE_{W(\MCO^\flat)})\]
is bijective. \par
 Consider the diagram of adic spaces:
    \begin{center}
\begin{tikzcd}
Y^2 \arrow[d, shift right] \arrow[d, shift left] & \cdots \arrow[l, hook] & {Y^2_{[0,2]}} \arrow[d, shift right] \arrow[d, shift left] \arrow[l, hook] & {Y^2_{[0,1]}} \arrow[d, shift right] \arrow[d, shift left] \arrow[l, hook] \\
Y^1                                              & \cdots \arrow[l, hook] & {Y^1_{[0,2]}} \arrow[l, hook]                                              & {Y^1_{[0,1]}} \arrow[l, hook]
\end{tikzcd}
\end{center}
and the diagram of rings:
\begin{center}
\begin{tikzcd}
\cdots \arrow[r, hook] & {A^2_{[0,2]}} \arrow[r, hook]                                              & {A^2_{[0,1]}} \arrow[r, hook]                                              & B^2 \\
\cdots \arrow[r, hook] & {A^1_{[0,2]}} \arrow[u, shift right] \arrow[u, shift left] \arrow[r, hook] & {A^1_{[0,1]}} \arrow[u, shift right] \arrow[u, shift left] \arrow[r, hook] & B^1 \arrow[u, shift right] \arrow[u, shift left]
\end{tikzcd}
\end{center}
The fact that the \(B^i\) are not Huber rings and the \(Y^2_{[0,m]}\) are not affinoid prevent us from combining the two diagrams into one.\par
Let's take a closer look at the first diagram. The integral \(\MCY\)-bundle \(\MCE_{\MCY}\) on \(A^+\) is given by a vector bundle \(\MCE_{Y^1}\) on \(Y^{1}\) together with an isomorphism \(\alpha\) between the two pullbacks of \(\MCE_{Y^1}\) to \(Y^{2}\) satisfying some cocycle condition on \(\MCE_{Y^3}\). Now the following pieces of data are equivalent:
\begin{enumerate}
    \item a section in \(H^0(\Spd(A^+),\MCE_{\MCY})\)
    \item By definition, this is equivalent to a section \(a\in H^0(Y^{1}, \MCE_{Y^1})\) such that the two pullbacks to \(Y^{2}\) are compatible with \(\alpha\).
    \item By analytic gluing, the section \(a\) is equivalent to a family of compatible sections 
    \[a_m\in H^0(Y^{1}_{[0,m]}, \MCE_{Y^1_{[0,m]}}),\]
    of the base changes \(\MCE_{Y^1_{[0,m]}}\) of \(\MCE_{Y^1}\). Since the symmetric domains \(Y^2_{[0,m]}\) form an open cover of \(Y^2\), the sheaf condition shows that the two pullbacks of \(a\) to \(Y^2\) are compatible with \(\alpha\) if and only if for all the \(a_m\) the two pullbacks to \(Y^{2}_{[0,m]}\) are compatible with \(\alpha_m\coloneqq \alpha|_{Y^2_{[0,m]}}\).
    \item Now, as the \(Y^1_{[0,m]}\) are analytic affinoid, the \(\MCE_{Y^1_{[0,m]}}\) are given by finite projective modules \(M_{A^1_{[0,m]}}\). Denote by \(M^{(j)}_{A^2_{[0,m]}}\) the global sections of \(\operatorname{pr}_j^*\MCE_{Y^1_{[0,m]}}\), for \(j=1,2\). Writing \(M_{A^1_{[0,m]}}=e(A^1_{[0,m]})^s\) for an idempotent matrix \(e\), taking pullback and global sections yields
    \[M^{(j)}_{A^2_{[0,m]}}=i_j^A(e)(A^2_{[0,m]})^s\cong M_{A^1_{[0,m]}}\otimes_{A^1_{[0,m]},i_j^A}A^2_{[0,m]}.\]
    Thus these modules are finite projective. Now our compatible sections from (3) are equivalent to compatible elements \(a_m\in M_{A^1_{[0,m]}}\), such that for each \(a_m\) the images under the two maps
    \[M_{A^1_{[0,m]}}\longrightarrow M^{(j)}_{A^2_{[0,m]}},\qquad j=1,2,\]
    are compatible with \(\alpha_m\colon M^{(1)}_{A^2_{[0,m]}}\isom M^{(2)}_{A^2_{[0,m]}}\), the map on global sections induced by the restricted descent isomorphism.
\end{enumerate}
On the other hand, elements in \(H^0(\Spd(A^+),\MCE_{W(\MCO^\flat)})\) are by definition given by sections of some finite projective \(B^1=W(A^+((T^{\pRoots})))\)-module \(M_{B^1}\) whose two pullbacks to the finite projective \(B^2\)-module \(M_{B^2}\) are compatible with \(\alpha_{B^2}\). Now to show that our map \(F\colon H^0(\Spd(A^+),\MCE_{\MCY})\to H^0(\Spd(A^+), \MCE_{W(\MCO^\flat)})\) is a bijection, it is enough to show that for each fixed \(m\), the map
\[
\begin{aligned}
G\colon {}&\{\text{sections of \(M_{A^1_{[0,m]}}\) compatible with \(\alpha_m\)}\}\\
&\longrightarrow
\{\text{sections of \(M_{B^1}\) compatible with \(\alpha_{B^2}\)}\}
\end{aligned}
\]
is a bijection.
Indeed, this means we can then successively lift compatible-with-\(\alpha\) sections of \(M_{B^1}\) to a compatible system of compatible-with-\(\alpha\) sections of \(M_{A^1_{[0,m]}}\), which by characterization (4) above glue to an element of \(H^0(\Spd(A^+),\MCE_{\MCY})\).\par
We may assume that the descended \(W(A^+)\)-module has constant rank. The bijectivity of \(G\) is essentially proved in \Cref{caseILemmaGeneralStatement}: assumption (1) is clear in this case, assumption (2) follows from \Cref{KimExtTheorem} (to descend our \(W(\MCO^\flat)\)-bundle to a finite projective \(W(A^+)\)-module) and the fact that \(A^+\) is a product of valuation rings and \((W(A^+), (p))\) is a henselian pair (to show that it is free), and assumption (3) is shown in \Cref{caseILemmaAssumption3}. This proves the bijectivity of the map between global sections and thus the theorem.
\end{proof}

\subsection{Semilinear actions and Galois descent}
We now turn towards the second principal case of our main theorem, where the pseudouniformizer of the integral perfectoid ring is a non-zero divisor. In this subsection, we introduce some basic language for phrasing descent data as continuous semilinear actions (\Cref{proposalGeometricSemilinear}),  and how to read effectivity of the descent datum from such actions.

Let \((R_\infty,R_\infty^+)\) be a complete sheafy Tate Huber pair and \(X_\infty=\Spa(R_\infty,R_\infty^+)\) the associated affinoid analytic adic space. For a profinite set \(S\), let \(\underline S\) denote the associated v-sheaf (as in \cite[\S9.3]{SW20}) and \(C(S,R_\infty)\) the ring of continuous functions with the compact-open topology. For every finite projective \(R_\infty\)-module \(N\), equipped with its natural topology (defined as e.g. the quotient topology from \(R_\infty^n\twoheadrightarrow N\)), one has
\begin{equation*}
C(S,R_\infty)\otimes_{R_\infty} N\simeq C(S,N)
\end{equation*}
as can be seen from the split surjection \(C(S,R_\infty^n)\twoheadrightarrow C(S,N)\). If \(X_\infty\) happens to be affinoid-perfectoid, the profinite-product construction of \cite[\S9.3]{SW20} identifies
\[
X_\infty^\diamond\times\underline S
\simeq\Spa\bigl(C(S,R_\infty),C(S,R_\infty^+)\bigr),
\]
and the space on the right is affinoid perfectoid.

\begin{mydef}\label{proposalProetaleTorsor}(\cite[Definition 10.12 and Lemma 10.13]{etCohDiamonds})
Let \(\MCX\) be a small v-sheaf over \(\Spd(\ZZ_p)\). A v-cover \(X_\infty^\diamond\to\MCX\) by an affinoid analytic adic space \(X_\infty=\Spa(R_\infty,R_\infty^+)\), equipped with a continuous right action of a profinite group \(\Gamma\) over \(\MCX\), is called a \emph{pro-\'{e}tale \(\Gamma\)-torsor} if the canonical map
\begin{equation*}\tag{*}\label{proposalCechProduct}
X_\infty^\diamond\times\underline\Gamma
\xrightarrow{\ \sim\ }
X_\infty^\diamond\times_{\MCX}X_\infty^\diamond,
\qquad (x,\sigma)\longmapsto(x,x\sigma)
\end{equation*}
is an isomorphism. We write the corresponding continuous action on \((R_\infty,R_\infty^+)\) on the left.
\end{mydef}

\begin{myrem}\label{proposalTensorWarning}
If \(\MCX=X^\diamond\) is represented by an affinoid analytic adic space \(X=\Spa(R,R^+)\) as well and the torsor \(X_\infty^\diamond\to \MCX\) is induced by a \(\Gamma\)-equivariant morphism \(X_\infty\to X\) with trivial action on \(X\) \footnote{automatic if \((\cdot)^{\diamond}\) respects morphisms going out of \(X_\infty\), e.g. if it is perfectoid}, there is a natural comparison map
\[
R_\infty\widehat\otimes_R R_\infty\longrightarrow C(\Gamma,R_\infty),
\qquad b_1\otimes b_2\longmapsto
\bigl(\sigma\mapsto b_1\sigma(b_2)\bigr),
\]
where the left side is the ordinary Banach completed tensor product used in the fibre product of adic spaces. If both \(X\) and \(X_\infty\) are affinoid perfectoid, this map is an isomorphism: the adic fiber product \(X_\infty\times_X X_\infty\) is affinoid perfectoid with ring of functions \(R_\infty\widehat\otimes_R R_\infty\), and its associated diamond is the fiber product above.
\end{myrem}

\begin{mydef}\label{proposalSemilinearAction}
A \emph{continuous semilinear action} on a finite projective \(R_\infty\)-module \(N\) is an action 
\begin{align*}
\Gamma\times N&\longrightarrow N\\
(\sigma, a)&\mapsto T_\sigma(a)
\end{align*} 
such that:
\begin{enumerate}
\item Each slice \(T_\sigma\colon N\to N\) is a \(\sigma\)-linear bijection, i.e., an additive bijection satisfying
\[
T_\sigma(bv)=\sigma(b)T_\sigma(v),
\]
\item the action map \(\Gamma\times N\to N\) is continuous.
\end{enumerate}
\end{mydef}

\begin{myprop}\label{proposalGeometricSemilinear}
Assume the situation of \Cref{proposalProetaleTorsor}, and furthermore that \(X_\infty=\Spa(R_\infty,R_\infty^+)\) is affinoid perfectoid. Then the functor
\[
\{\text{v-bundles on }\MCX\}\longrightarrow
\left\{\begin{gathered}
\text{finite projective }R_\infty\text{-modules with}\\
\text{continuous semilinear }\Gamma\text{-action}
\end{gathered}\right\},
\]
given by evaluation on \(X_\infty\) followed by taking global sections is an equivalence of categories.
\end{myprop}
\begin{proof}
Vector bundles on perfectoid spaces satisfy v-descent by \cite[Lemma 17.1.8]{SW20}, so v-bundles on \(\MCX\) are equivalent to vector bundles on \(X_\infty\) with descent datum. Thus, we have to identify descent data and continuous semilinear actions for each finite projective \(R_\infty\)-module \(N\). Every \(R_\infty\)-linear map between finite projective modules is continuous, and by considering internal Homs, it suffices to show that on a fixed \(R_\infty\)-module \(N\), the following data are equivalent:
\begin{enumerate}
	\item A descent datum on the associated vector bundle \(\MCE_\infty\) on \(X_\infty\)
	\item A continuous semilinear action of \(N\) under \(\Gamma\)
\end{enumerate}

We first embed both data types in a common ambient set without continuity condition. Let \(\mathscr A(N)\) be the set of all semilinear \(\Gamma\)-actions on \(N\), and let
\[
\mathscr C(N)\subseteq\mathscr A(N)
\]
be the subset of continuous semilinear actions in the sense of \Cref{proposalSemilinearAction}. For the descent side, \eqref{proposalCechProduct} identifies the overlap \(X_\infty^\diamond\times_{\MCX}X_\infty^\diamond\) with the v-sheaf associated to the affinoid perfectoid space \(\Spa(C(\Gamma,R_\infty),C(\Gamma,R_\infty^+))\). Under this translation, the two projections are given by
\[
p_1^*(b)(\sigma)=b,\qquad p_2^*(b)(\sigma)=\sigma(b).
\]
Thus a descent datum on \(\MCE_\infty\) is an isomorphism
\[
\beta\colon N\otimes_{R_\infty,p_2^*}C(\Gamma,R_\infty)\xrightarrow{\sim}
N\otimes_{R_\infty,p_1^*}C(\Gamma,R_\infty)=C(\Gamma,N),
\]
satisfying the cocycle identity on the global sections of the triple fiber product. Evaluation in \(\sigma\) produces a family of maps
\[\beta_\sigma\colon \sigma^* N\longrightarrow N \]
and interpreting this as a \(\sigma\)-linear map \(N\to N\) for each \(\sigma\) produces an element of \(\mathscr A(N)\). This is the easy direction of \cite[0CDQ]{Stacks}, where this construction is shown to be a bijection for finite \(\Gamma\) (in the setting of rings). Let \(\mathscr D(N)\subset \mathscr A(N)\) be the subset of elements arising through this construction.\par
To see that both subsets are equal, 

we check that continuity of \(T\) and \(\beta\) is already guaranteed by the seemingly weaker condition that the orbit map
\[
\Gamma\longrightarrow N,\qquad \sigma\longmapsto\beta_\sigma(v\otimes1)=T_\sigma(v)
\]
is continuous for each fixed \(v\in N\).

For the action side, let \(n_1,\ldots,n_d\) be generators of \(N\). For some \(v= \sum_{k=1}^d v_kn_k\), the function
\[
T_\sigma(v)=\sum_{k=1}^d\sigma(v_k)T_\sigma(n_k)
\]
is then clearly continuous in both \(\sigma\) and \(v\).\par

For the descent datum, the same consideration shows that for some \(h\in C(\Gamma, R_\infty)\), the term \(h(\sigma)T_\sigma(v)\) varies continuously in \(\sigma\), so we get a well defined \(C(\Gamma,R_\infty)\)-linear map
\begin{align*}
\beta\colon N\otimes_{R_\infty,p_2^*}C(\Gamma,R_\infty)&\longrightarrow
N\otimes_{R_\infty,p_1^*}C(\Gamma,R_\infty)=C(\Gamma,N),\\
v\otimes h&\longmapsto
\bigl(\sigma\mapsto h(\sigma)T_\sigma(v)\bigr).
\end{align*}

Let 
\[\iota\colon X_\infty^\diamond\times_\MCX X_\infty^\diamond\longrightarrow X_\infty^\diamond\times_\MCX X_\infty^\diamond\] be the involution interchanging the two factors. Under \eqref{proposalCechProduct}, it is
\[
\iota(x,\sigma)=(x\sigma,\sigma^{-1}).
\]
Since \(\iota\) exchanges \(p_1\) and \(p_2\), the pullback \(\iota^*\beta\) is a map in the reverse direction of \(\beta\), and the cocycle condition implies 
\[
\beta\circ\iota^*\beta=\id,\qquad
\iota^*\beta\circ\beta=\id.
\]
Thus \(\beta\) is an isomorphism, and so both it and its inverse are continuous, since linear maps between finite projective \(C(\Gamma,R_\infty)\)-modules are continuous in their direct-summand topologies. This gives the required descent datum.
\end{proof}

This finishes the translation from geometric descent data (if \(X_\infty\) is affinoid perfectoid) to continuous semilinear \(\Gamma\)-actions. We now want to go one step further and express semilinear actions as continuous cocycles. This requires us to "stabilize" the action by comparing it against a "model" of the underlying module \(N\):

\begin{myconst}(reference semilinear action)\label{proposalReferenceAction}
Let \(R\to R_\infty\) be a continuous homomorphism of complete Tate rings and let \(\Gamma\times R_\infty\to R_\infty\) be a continuous action by a profinite group such that the image of \(R\) is contained in \(R_\infty^\Gamma\). Then, for any finite projective \(R\)-module \(P\), there is a canonical continuous semilinear action on \(P\otimes_R R_\infty\), given by
\[
s_\sigma(v\otimes b)=v\otimes\sigma(b).
\]
Now let \(N\) be any finite projective \(R_\infty\)-module with continuous semilinear action \(T\) and assume there exists an identification
\[
N\simeq P\otimes_R R_\infty.
\]
Then the canonical action on \(P\otimes_R R_\infty\) becomes a \emph{reference action} for \(T\). Define a map of sets
\begin{align*}
\Gamma &\longrightarrow \Aut_{R_\infty}(N),\\
\sigma&\mapsto a(\sigma)=T_\sigma s_\sigma^{-1}
\end{align*}
Direct calculation gives the \emph{cocycle identity}
\[a(\sigma\tau)=a(\sigma)\sigma(a(\tau)),\]
with the convention of \cite[\S2.2]{Sen}. Conversely, any \(a\) satisfying this identity makes \(T_\sigma\coloneqq a(\sigma)s_\sigma\) an action. Conjugating \(T\) by \(c\in\Aut_{R_\infty}(N)\), that is, replacing \(T_\sigma\) by \(c^{-1}T_\sigma c\), changes the cocycle to the cohomologous cocycle
\[
a'(\sigma)=c^{-1}a(\sigma)\sigma(c).
\]
\end{myconst}
Doing the converse calculation, we arrive at:

\begin{mylemma}
The above constructions define a bijection
\[\{\text{\emph{continuous semilinear \(\Gamma\)-actions on \(N\)}}\} \longrightarrow \{\text{\emph{continuous cocycles \(\Gamma\to \Aut_{R_\infty}(N)\)}}\}.\]
With the reference action fixed, two actions \(T,S\) are sent to cohomologous cocycles if and only if there exists \(c\in\Aut_{R_\infty}(N)\) such that \(S_\sigma=c^{-1}T_\sigma c\) for every \(\sigma\in\Gamma\).
\end{mylemma}

Note that the trivial cocycle corresponds to the reference semilinear action, not to the action fixing every element of \(N\).

For the last result of the subsection, we will assume the following geometric setup:
\begin{itemize}
    \item \(K\) is a complete discretely valued field of mixed characteristic \((0,p)\).
    \item \(L/K\) is an infinite Galois extension, \(\widehat L\) its completion, and \(\Gamma\coloneqq\Gal(L/K)\).
    \item \((B,B^+)\) is a complete Huber pair over \((K,\MCO_K)\)
	\item \((\hat B_\infty,\hat B_\infty^+)\coloneqq(B,B^+)\widehat\otimes_{(K,\MCO_K)}(\widehat L,\MCO_{\widehat L})\). This gives \(\hat B_\infty\) an induced action by \(\Gamma\).
	\item For \(K_i/K\) a finite Galois subextension of \(L/K\), put \(B_i=B\otimes_K K_i\) and \(\Gamma_i=\Gal(L/K_i).\)
\end{itemize}
Neither \(B\) nor \(\hat B_\infty\) is assumed perfectoid.\par
By Ax--Sen--Tate, 
\[\hat B_\infty^{\Gamma_i}=B_i,\qquad \hat B_\infty^\Gamma=B.\]
A standard coefficient-approximation-argument also yields 
\[
(Q\otimes_{B_i}\hat B_\infty)^{\Gamma_i}=Q
\]
for any finite projective \(B_i\)-module \(Q\), with the reference action on the left.\par
This implies the extension of scalars functor
\[
\{\text{\emph{finite projective \(B\)-modules}}\} \longrightarrow
\left\{\begin{gathered}
\text{\emph{finite projective \(\hat B_\infty\)-modules with}}\\
\text{\emph{continuous semilinear \(\Gamma\)-action}}
\end{gathered}\right\}
\]
is fully faithful. We will call objects in the essential image \emph{effective}. We can detect effectivity at finite levels:

\begin{myprop}\label{proposalFiniteEffectivity}
For a  finite projective \(\hat B_\infty\)-module \(N\) with continuous semilinear \(\Gamma\)-action \(T\), the following are equivalent:
\begin{enumerate}
\item \((N,T)\) is effective over \(B\), i.e., there is a finite projective \(B\)-module \(M\) and a \(\Gamma\)-equivariant isomorphism \(M\otimes_B \hat B_\infty\simeq N\), with the reference semilinear action on the left.
\item For some \(i\), there are a finite projective \(B_i\)-module \(Q_i\) and a \(\Gamma_i\)-equivariant isomorphism
\[
Q_i\otimes_{B_i}\hat B_\infty\simeq N,
\]
where the action on the left is the reference action.
\end{enumerate}
When these hold, the descended module is canonically \(M=N^\Gamma\). 
\end{myprop}

\begin{proof}
Condition (1) implies (2) by taking \(Q_i=M\otimes_B B_i\). Conversely, the invariant calculation above identifies \(N_i\coloneqq N^{\Gamma_i}\) with \(Q_i\). Since \(\Gamma_i\) is normal, \(N_i\) inherits a semilinear action of \(\Gamma/\Gamma_i\). Finite Galois descent along \(B_i/B\), as in \cite[0CDQ and 058S]{Stacks} (this needs the invariant calculations above) gives a finite projective module
\[
M=N_i^{\Gamma/\Gamma_i}=N^\Gamma,
\qquad M\otimes_B B_i\simeq N_i.
\]
Composing with the equivariant multiplication isomorphism \(N_i\otimes_{B_i}\hat B_\infty\simeq N\) proves (1).
\end{proof}

Finally, we translate this effectivity criterion into the cocycle formulation.
\begin{mycor}\label{proposalLocalTriviality}
Let \((N,T)\) be as above. Choose a finite projective \(B_i\)-module \(P_i\) and an identification \(N\simeq P_i\otimes_{B_i}\hat B_\infty\), not assumed equivariant. Let
\[
a\colon\Gamma_i\longrightarrow\Aut_{\hat B_\infty}(N)
\]
be the continuous cocycle corresponding to \(T|_{\Gamma_i}\) relative to this model. Then the following are equivalent:
\begin{enumerate}
\item \((N,T)\) is effective over \(B\) 
\item \(a|_{\Gamma_j}\) is cohomologous to the trivial cocycle for some \(j\geq i\)
\end{enumerate}
Explicitly, the second condition amounts to the existence of \(c\in\Aut_{\hat B_\infty}(N)\) such that
\[
c^{-1}a(\sigma)\sigma(c)=1
\quad\text{for every }\sigma\in\Gamma_j.
\]
\end{mycor}

\begin{proof}
\((2)\implies (1)\) is \Cref{proposalFiniteEffectivity}. For \((1)\implies (2)\), let \(M\) be the descended module and put \(M_i=M\otimes_B B_i\). The modules \(M_i\) and \(P_i\) then become isomorphic over \(\hat B_\infty\). To apply the other direction of \Cref{proposalFiniteEffectivity}, we have to show that they are already isomorphic at some finite level \(j\) with \(K_i\subseteq K_j\). Choose mutually inverse isomorphisms \(u,v\) after base change to \(\hat B_\infty\) and, using density of \(\bigcup_{j\geq i}B_j\subset \hat B_\infty\), approximate both on some finite level \(j\) by
\[
u_j\colon P_i\otimes_{B_i}B_j\longrightarrow M_i\otimes_{B_i}B_j,
\qquad
v_j\colon M_i\otimes_{B_i}B_j\longrightarrow P_i\otimes_{B_i}B_j
\]
so that \(v_ju_j\) and \(u_jv_j\) differ from their respective identities by operators of norm less than one. The standard Neumann series argument thus shows that the two compositions are isomorphisms, thus the same must be true for \(u_j\) (and \(v_j\), of course).
\end{proof}
\subsection{A descent result via Sen operators}
Heavily inspired by the proof of \cite[Theorem 2.3.8]{PR}, we will use the results of \cite{Sen} and the translations proved in the last subsection to prove \Cref{SenDescent}. Roughly speaking, the theorem tells us when effectivity of a descent datum of vector bundles along certain pro-étale covers of adic spaces can be tested both generically and pointwise. Consider the following field theoretic setup:

\begin{itemize}
    \item \(K\) a complete discretely valued field of mixed characteristic \((0,p)\) with perfect residue field, algebraic closure \(\overline K\), and Galois group \(G_K\coloneqq\Gal(\overline K/K)\),
    \item \(\chi_K\colon G_K\to\ZZ_p^\times\) the \(p\)-adic cyclotomic character, characterized by \(\sigma(\zeta)=\zeta^{\chi_K(\sigma)}\) for every \(\sigma\in G_K\) and \(\zeta\in\mu_{p^\infty}\),
    \item \(K_\infty\coloneqq\overline K^{\ker\chi_K}=K(\mu_{p^\infty})\), the cyclotomic extension cut out by \(\chi_K\), obtained by adjoining all \(p\)-power roots of unity to \(K\), with completion \(\hat K_\infty\),
    \item \(\Gamma\coloneqq \Gal(K_{\infty}/K)\), equipped with the induced injective character \(\chi\colon\Gamma\hookrightarrow\ZZ_p^\times\).
\end{itemize}
We fix furthermore 
\begin{itemize}
\item \((B,B^+)\) a complete Huber pair over \((K,\MCO_K)\) with base-changed versions
\[B_\infty\coloneqq B\otimes_K K_\infty,\quad B^+_\infty\coloneqq B^+\otimes_{\MCO_K}\MCO_{K_\infty}, \quad (\hat B_\infty,\hat B^+_\infty)\coloneqq (B,B^+)\hat\otimes_{(K,\MCO_K)} (\hat K_\infty,\hat \MCO_{K_\infty}).\]
\item A finite projective \(\hat B_\infty\)-module \(\hat M_\infty\) with a continuous semilinear \(\Gamma\)-action
\[
\Gamma \times \hat M_{\infty} \longrightarrow \hat M_{\infty},
\qquad (\sigma,v)\mapsto T_\sigma(v)
\]
as in \Cref{proposalSemilinearAction} (\(\Gamma\) acts on \(\hat B_\infty=B\widehat\otimes_K\hat K_\infty\) through the second factor).
\item \(X=\Spa(B,B^+)\), and for \(x\in X\) the non-completed and completed residue fields
\[\kappa(x)=\Frac(B/\supp(x)), \qquad k_x\coloneqq\widehat{\kappa(x)}.\]
\item The base changes
\[
\hat k_{\infty, x}\coloneqq k_x\widehat\otimes_K\hat K_\infty
\simeq k_x\widehat\otimes_B\hat B_\infty,
\qquad
\hat M_{\infty,x}\coloneqq\hat M_\infty\otimes_{\hat B_\infty}\hat k_{\infty, x},
\]
with an induced semilinear action
\[
\Gamma \times \hat M_{\infty,x} \longrightarrow \hat M_{\infty,x},
\qquad T_{\sigma,x}(v\otimes h)=T_\sigma(v)\otimes\sigma(h).
\]
Note that \(\hat k_{\infty, x}\) need not be a field. 
\end{itemize}

We will prove the following theorem:

\begin{mythm}\label{SenDescent}
    In the above setup, assume there exists an open subset \(U\subseteq \Spa(B,B^+)\) such that:
    \begin{enumerate}
        \item The evaluation map
        \[
        B\longrightarrow\prod_{x\in U}k_x,
        \qquad b\longmapsto\bigl(b(x)\bigr)_{x\in U},
        \]
        is injective. This holds in particular if \(B\) is uniform and either:
        \begin{enumerate}
            \item \(U\subseteq \Spa(B,B^+)\) is dense or
            \item \(U=\Spa(B,B^+)-V(f)\) for \(f\in B\) a non-zero divisor.\footnote{In contrast to the situation for schemes, it is not clear that b) implies a).}
        \end{enumerate}
        \item For every \(x\in U\), the action \((\hat M_{\infty,x},T_x)\) is effective over \(k_x\) in the sense of the last section, i.e., for every \(x\in U\) there exist a finite-dimensional \(k_x\)-vector space \(M_x\) and a \(\Gamma\)-equivariant isomorphism
        \[
        M_x\otimes_{k_x}\hat k_{\infty, x}\xrightarrow{\ \sim\ }\hat M_{\infty,x},
        \]
        where the reference action on the source is \(\sigma(m\otimes h)=m\otimes\sigma(h)\).
    \end{enumerate}
    Then there exist a finite projective \(B\)-module \(M\) and a \(\Gamma\)-equivariant isomorphism \(M\otimes_B\hat B_\infty\simeq\hat M_\infty\), with the reference action on the source.
\end{mythm}

\begin{myrem}(Role of pointwise effectivity)
\begin{enumerate}
    \item The action \(T\) is given globally; only its effectivity is tested at the points of \(U\).
    \item The extension-of-scalars functor from finite projective \(B\)-modules to continuous semilinear \(\hat B_\infty\)-modules is fully faithful by \Cref{proposalFiniteEffectivity}. Applying loc. cit. to \(k_x\) in place of \(B\), we get similar full faithfulness for the functor from finite-dimensional \(k_x\)-vector spaces to continuous semilinear \(\hat k_{\infty, x}\)-modules. Thus, an effective global descent \(M\) specializes canonically to the pointwise descent data. In this sense, pointwise effectivity is really a property instead of a datum here.
\end{enumerate}
\end{myrem}

To make use of the formulation of descent data using cocycles in \Cref{proposalReferenceAction}, we have to find a ``model'' of \(\hat M_\infty\) over some finite subextension, by which we mean the following:

\begin{mylemma}\label{SenFindingP}
    There is a finite subextension \(K\subset K'\subset K_\infty\) and a finite projective \(B'=B\otimes_K K'\)-module \(P\) such that
    \[\hat M_\infty\cong P\otimes_{B'} \hat B_{\infty}.\]
\end{mylemma}

Note that the isomorphism \(\hat M_\infty\cong P\otimes_{B'} \hat B_{\infty}\) need not respect the action \(T\), so it does not yet solve the descent problem.

\begin{proof}
As a first step, we want to find a finite projective \(B_\infty\)-module \(P_\infty\) whose completion (in other words, its base change to \(\hat B_\infty\)) is isomorphic to \(\hat M_\infty\). \cite[Corollary 2.1.22 (c)]{BC22} gives us such a \(P_\infty\), under the condition that \(B_\infty\) contains a non-unital open subring \(B_\infty^I\) such that: 
\begin{enumerate}
\item \(B_\infty^I\) is \emph{henselian} as a non-unital ring. This property is discussed in \cite[§2.1.1]{BC22}, but a sufficient condition for a non-unital ring to be henselian is if it can be realized as the ideal of a henselian pair.
\item The subspace topology on \(B_\infty^I\) is linear, i.e., it admits a neighbourhood basis of \(0\) consisting of ideals in \(B_\infty^I\).
\end{enumerate}

We will show that these properties are satisfied for \(B^I_\infty\coloneqq \varpi B_{\infty,0}\), where \(\varpi\in B\) is a pseudouniformizer (it will then automatically be a pseudouniformizer for \(B_\infty\)) and \(B_{\infty,0}\subset B_\infty\) is a ring of definition, which we assume to be obtained by base change from a ring of definition \(B_0\) of \(B\). Now \(\varpi B_{\infty,0}\) clearly satisfies the second property above. For the henselian property, note that \((B_0,\varpi B_0)\) is a henselian pair since \(B_0\) is \(\varpi\)-adically complete. But \(B_0\to B_{\infty,0}\) is the base change of the integral morphism \(\MCO_K\to \MCO_{K_\infty}\) and hence an integral morphism itself, and henselian pairs are preserved by integral morphisms (\cite[09XK]{Stacks}), so \((B_{\infty,0},\varpi B_{\infty,0})\) is a henselian pair and thus \(\varpi B_{\infty,0}\) is henselian in the above sense.\par
This proves the existence of \(P_\infty\). Going further, we have an isomorphism
\[B_\infty=\varinjlim_{K\subseteq K'\subset K_\infty} B\otimes_K K',\]
with the colimit running over all finite subextensions. This translates into a filtered colimit of the (equivalence classes of) finite projective modules, with the transition maps given by base change: As each finite projective \(B_\infty\)-module can be described by finitely many generators and relations, there must be some \(B'\coloneqq B\otimes_K K'\) where all of these relations can be defined. Thus we find a \(B'\) and a finite projective \(B'\)-module \(P\) which after base change to \(B_\infty\) becomes isomorphic to \(P_\infty\), and after completing becomes isomorphic to \(\hat M_\infty\), as desired.
\end{proof}

In \cite[\S2.3]{Sen}, Sen developed a theory to test triviality of a cocycle on open normal
subgroups of a cyclotomic Galois group by the vanishing of a certain endomorphism. Let’s briefly recall the construction.\par
Let \(E\) be a not necessarily commutative associative unital \(B\)-algebra, with \(B\) acting centrally, which is finite projective\footnote{This finiteness and projectivity assumption is convenient here but unnecessary in Sen's setup, which allows general Banach \(K\)-algebras.} as a \(B\)-module. We equip \(E\) with its natural Banach topology. The case to keep in mind here is the endomorphism algebra of a finite projective \(B\)-module. The base change
\[
E_{\hat B_\infty}\coloneqq E\otimes_B\hat B_\infty
\]
has a natural \(\Gamma\)-action through the second factor. Let \(\alpha\colon\Gamma\to E_{\hat B_\infty}^\times\) be a continuous cocycle in the convention of \Cref{proposalReferenceAction}. For \(K''/K\) a sufficiently large finite subextension of \(K_\infty/K\), define
\[
\Gamma''\coloneqq\Gal(K_\infty/K''), \qquad B''=B\otimes_K K'',\qquad E''=E\otimes_B B''=E\otimes_K K''.
\]
Sen's construction then gives \(m\in E_{\hat B_\infty}^\times\) and a continuous homomorphism
\[
\rho\colon\Gamma''\longrightarrow(E'')^\times,
\qquad \rho(\sigma)=m^{-1}\alpha(\sigma)\sigma(m).
\]
Since \(\Gamma''\) acts trivially on \(E''\), this is just an ordinary group homomorphism.

After possibly enlarging \(K''\) and correspondingly restricting \(\rho\), we may assume that \(\chi(\Gamma'')\subseteq1+p^2\ZZ_p\) and that \(\rho(\Gamma'')\) lies in a neighbourhood of \(1\) where the logarithm and exponential are inverse. Choose a topological generator \(\gamma\) of \(\Gamma''\simeq\ZZ_p\) and define the \emph{finite-level operator}
\begin{equation*}
\Theta_{\alpha,m,K''}\coloneqq
\frac{\log\rho(\gamma)}{\log\chi(\gamma)}\in E''.
\end{equation*}
For \(z\in\ZZ_p\), continuity and the homomorphism property give
\[
\log\rho(\gamma^z)=z\log\rho(\gamma),\qquad
\log\chi(\gamma^z)=z\log\chi(\gamma).
\]
Indeed, the identities hold for integral powers and extend by continuity. Thus the quotient \(\log\rho(\sigma)/\log\chi(\sigma)\) is constant for \(1\ne\sigma\in\Gamma''\), and
\[
\Theta_{\alpha,m,K''}=
\lim_{\substack{\sigma\to1\\\sigma\ne1}}
\frac{\log\rho(\sigma)}{\log\chi(\sigma)}\in E'',\qquad
\rho(\sigma)=\exp\bigl(\log\chi(\sigma)\Theta_{\alpha,m,K''}\bigr).
\]
This is the existence argument in \cite[\S2.3, p.~27]{Sen}.

The \emph{Sen operator} now is the conjugate
\begin{equation*}
\phi_\alpha\coloneqq m\Theta_{\alpha,m,K''}m^{-1}
\in E_{\hat B_\infty}.
\end{equation*}

{\emergencystretch=1em Since \(E''\hookrightarrow E_{\hat B_\infty}\) is injective, this definition implies that \(\phi_\alpha=0\) if and only if \(\Theta_{\alpha,m,K''}=0\).\par}

\begin{mythm}\emph{(\cite[\S2.4, 2.5]{Sen})}\label{SenFundamentalThm}
The construction above enjoys the following properties:
\begin{enumerate}
    \item For fixed \(\chi\), \(\phi_\alpha\) exists and is independent of the auxiliary choices \(K''\), \(m\), and \(\rho\).
    \item It is compatible with continuous base change \(B\to D\) between commutative Banach \(K\)-algebras, including the finite-level operators for the compatible choices.
    \item \(\phi_\alpha=0\) if and only if \(\Theta_{\alpha,m,K''}=0\) if and only if \(\alpha\) is cohomologous to the trivial cocycle on some open normal subgroup of \(\Gamma\).
\end{enumerate}
\end{mythm}

\begin{proof}[Proof of \Cref{SenDescent}]
First, let's confirm the ``in particular'' statements. If \(B\) is uniform, evaluation at all points of \(X\) is injective by \cite[Theorem 5.2.1]{SW20}. Now if \(U\) is dense, a function vanishing at its points vanishes on \(X\), since its zero locus is closed. If \(U=X\setminus V(f)\) and \(g\) vanishes at all points of \(U\), then \(gf\) vanishes at all points of \(X\), so \(gf=0\), which implies \(g=0\) if \(f\) is a non-zero divisor.\par
For the rest of the proof, the idea is to use property (3) from \Cref{SenFundamentalThm} and \Cref{proposalLocalTriviality} to equate effectivity of descent with vanishing of \(\phi_\alpha\), and injectivity of the evaluation map in \Cref{SenDescent} (1) to imply global vanishing of \(\phi_\alpha\) from pointwise vanishing via the finite-level operator \(\Theta\). For the details, we will first have to fix yet some more notation, before proceeding in several steps:

\begin{itemize}
\item \emph{Finite-step model.} Using \Cref{SenFindingP}, choose a finite subextension \(K'/K\), a finite projective \(B'=B\otimes_K K'\)-module \(P\), and an identification
\[
\hat M_\infty\simeq P\otimes_{B'}\hat B_\infty.
\]
Put \(\Gamma'=\Gal(K_\infty/K')\) and
\[
E=\End_{B'}(P),\qquad
E_{\hat B_\infty}=E\otimes_{B'}\hat B_\infty
\simeq\End_{\hat B_\infty}(\hat M_\infty).
\]
\item \emph{Cocycle.}
Let \(s\) be the reference action associated with \(P\) and let
\begin{align*}
\alpha\colon\Gamma'&\longrightarrow E_{\hat B_\infty}^\times,\\
\sigma&\mapsto T_\sigma s_\sigma^{-1}
\end{align*}
be the associated cocycle as in \Cref{proposalReferenceAction}.
\item \emph{\'{E}tale cover.} For the chosen \(K'\), let
\[
\pi'\colon X'=\Spa(B',B'^+)\longrightarrow X,
\]
where \(B'^+\) is the integral closure of \(B^+\) in \(B'\), and define \(U'=\pi'^{-1}(U)\).
The map \(B\to B'\) is finite \'{e}tale.
\item \emph{Notation for points.} Throughout this proof, \(x\) denotes a point of the original \(U\), while \(y\) denotes a point of \(U'\) above \(x\). Write \(k_y=\widehat{\kappa(y)}\) for the completed residue field and set
\begin{align*}
\hat k_{\infty, y}&=k_y\widehat\otimes_{K'}\hat K_\infty,\qquad
\hat M_{\infty, y}&=\hat M_\infty\otimes_{\hat B_\infty}\hat k_{\infty, y}.
\end{align*}
Put \(E_y=\End_{k_y}(P\otimes_{B'}k_y)\simeq E\otimes_{B'}k_y\). The natural map \(\hat k_{\infty,x}\to\hat k_{\infty,y}\) is \(\Gamma'\)-equivariant, so condition (2) of \Cref{SenDescent} implies effectivity of \((\hat M_{\infty,y},T_y)\) over \(k_y\) by base change.
\item \emph{Sen operator.}
Note that \(K'\) is still discretely valued and \(K_\infty/K'\) is cut out by the cyclotomic character of \(G_{K'}\), so we will be able to apply Sen's construction to this field extension.
Choose a finite subextension \(K''/K'\) and \(m\in(E\otimes_{B'}\hat B_\infty)^\times\) as in Sen's construction, so that \(\rho(\sigma)=m^{-1}\alpha(\sigma)\sigma(m)\) takes values in \((E\otimes_{K'}K'')^\times\) on \(\Gal(K_\infty/K'')\). Put \(B''=B'\otimes_{K'}K''=B\otimes_K K''\), \(E''=E\otimes_{K'}K''\), and \(E''_y=E_y\otimes_{K'}K''\). Abbreviate
\[
\Theta=\Theta_{\alpha,m,K''}\in E\otimes_{K'}K'',\qquad
\phi=\phi_\alpha=m\Theta m^{-1}
\]
for the finite-level and full Sen operators.
\end{itemize}

 With all this in place, the actual proof is quite straightforward:

\begin{enumerate}
\item (\(\Theta_y=0\) for every \(y\in U'\)).
Fix \(y\in U'\) above \(x\in U\), and let \(\Theta_y\) be the image of \(\Theta\) in \(E''_y\); by \Cref{SenFundamentalThm} (2), this is the finite-level operator for the specialized cocycle with the induced choices. Then by applying \Cref{proposalLocalTriviality} (with base field \(K'\), base algebra \(k_y\) and model \(P\otimes_{B'}k_y\)) to the effective descent datum \((\hat M_{\infty,y},T_y)\), the cocycle \(\alpha_y\) is cohomologous to the trivial cocycle on an open subgroup of \(\Gamma'\), and by \Cref{SenFundamentalThm} (3), \(\Theta_y=0\).
\item (\(\Theta=0\)).
As \(K\to K'\) is finite, \(- \otimes_K K'\) is exact and commutes with products, so condition (1) of \Cref{SenDescent} lifts to \(K'\), giving an injection
\[
B'\lhook\joinrel\longrightarrow
\prod_{x\in U}(k_x\otimes_K K').
\]
Since \(B'/B\) is finite \'{e}tale, each \(k_x\otimes_K K'\) is a product of finite separable field extensions detected by the completed residue fields \(k_y\) at points above \(x\), so we furthermore get injectivity of
\[
B'\longrightarrow\prod_{y\in U'}k_y
\]
and, by finite projectivity of \(B''\) as a \(B'\)-module and of \(E''\) as a \(B''\)-module, also of
\[E''\longrightarrow\prod_{y\in U'}E''_y,\]
showing that \(\Theta=0\), since its image is \((\Theta_y)_y=0\).

\item (\((\hat M_\infty,T)\) is effective over \(B\)).
By \Cref{SenFundamentalThm} (3), \(\Theta=0\) implies that \(\phi=0\), which implies that \(\alpha\) is cohomologous to the trivial cocycle on an open normal subgroup of \(\Gamma'\). After restriction, it is therefore cohomologous to the trivial cocycle on \(\Gamma_j=\Gal(K_\infty/K_j)\) for some finite subextension \(K_j/K\) containing \(K'\). Applying \Cref{proposalLocalTriviality} to the original tower \(K_\infty/K\) and model \(P\) proves effectivity of the original descent datum over \(B\) and thus the theorem.
\end{enumerate}
\end{proof}
\subsection{\texorpdfstring{The equivalence away from \(p=0\)}
  {The equivalence away from p=0}}
The goal of this subsection is to use Sen's formalism established above to prove the following proposition, which can be regarded as a rational version of our main theorem under extra assumptions on \(A^+\).\footnote{This extra assumption guarantees that \(p\) gets inverted in a transcendental way in \(\Spa(A^+)\dot\times\Spa(\QQ_p)\). Dropping this assumption requires some subtle reflections on meromorphicity conditions as in \cite{GI23}; in \cite{Ans22}, this is carried out for a single point \(\Spec(k)\).}

\begin{myprop}\label{analyticLocuss}
    Let \(A^+\) be an integral perfectoid topological ring with pseudouniformizer \(\varpi\) such that \((A,A^+)\coloneqq(A^+[1/\varpi],A^+)\) is a perfectoid Huber pair. Then the functor 
    \begin{center}
    \begin{tikzcd}
    \emph{vector bundles on } \Spa(A^+)\dot\times \Spa(\QQ_p) \arrow[r] & \emph{rational }\MCY\emph{-bundles on } \Spd(A^+),
    \end{tikzcd}
    \end{center}
    is an equivalence of categories.
\end{myprop}

\begin{proof}
Full faithfulness follows from \Cref{proposalGeometricSemilinear,proposalFiniteEffectivity} and \Cref{baseChangeFunctorsProperties} (2), using the perfectoid covers below. For essential surjectivity, we use the formalism laid out in the previous subsection.
\(\Yoc(A,A^+)\) has an analytic cover by \((U^m)_{m\in \ZZ_{>0}}\) where
\[U^{m}\coloneqq \Spa(B^m,B^{m+})\coloneqq\Yoc(A,A^+)_{[\varpi^\flat]^m\leq p\neq 0}\longrightarrow \Yoc(A,A^+)\]
with
\[B^m\coloneqq \Ainf(A^+)\langle\frac{[\varpi^\flat]^m}{p}\rangle[1/p].\]
In the notation of \cite{SW20}, \(U^{m}=\MCY_{[1/m,\infty]}\). To apply \Cref{SenDescent}, take
\[K=\QQ_p, \qquad K_\infty=\QQ_p(\mu_{p^\infty}).\]
As in loc. cit., write \(\hat K_\infty\) for its completion and \(\hat\MCO_\infty=\ZZ_p[\mu_{p^\infty}]^{\wedge p}\) for the ring of integers of \(\hat K_\infty\) and form the completed tensor product of Huber pairs
\[(\hat B^m_\infty,\hat B^{m+}_\infty)\coloneqq
(B^m,B^{m+})\hat\otimes_{(\QQ_p,\ZZ_p)}(\hat K_\infty,\hat\MCO_\infty).\]
Then
\[U_\infty^m\coloneqq\Spa(\hat B^m_\infty,\hat B^{m+}_\infty)\longrightarrow U^m\]
is a pro-étale \(\Gamma\)-torsor, where \(\Gamma=\Gal(K_\infty/\QQ_p)\cong\ZZ_p^\times\).\footnote{Compare the similar perfectoid covers in the proof of \Cref{baseChangeFunctorsProperties} on full faithfulness of the functor from \(\MCY\)-bundles to v-bundles: there we adjoined all \(p\)-power roots of \(p\), while here we use the cyclotomic extension. While the argument there would have worked equally well with both perfectoid covers, here we really need the cyclotomic version to place ourselves in the context of the Sen operators from the last section.}
Since
\[C_\infty\coloneqq\bigl(W(A^{+\flat})\hat\otimes_{\ZZ_p}\hat\MCO_\infty\bigr)[1/p]\]
is perfectoid by \cite[proof of Theorem 2.3.8]{PR}, its rational localization
\[\hat B^m_\infty\cong C_\infty\langle [\varpi^\flat]^m/p\rangle\]
is perfectoid by \cite[Theorem 6.1.10]{SW20}, so \(U_\infty^m\) is affinoid perfectoid.\par
Let \(\MCE\) be a v-bundle on \(\Yoc(A,A^+)^\diamond\) arising from a rational \(\MCY\)-bundle. We want to apply \Cref{SenDescent} with \(K=\QQ_p\), \(K_\infty=\QQ_p(\mu_{p^\infty})\), \((B,B^+)=(B^m,B^{m+})\), \(\hat M_\infty=\MCE|_{U_\infty^m}\), carrying its continuous semilinear \(\Gamma\)-action by \Cref{proposalGeometricSemilinear}, and \(U=U^m\setminus V([\varpi^\flat])\).\par
Condition (1)(b) holds: the sousperfectoid algebra \(B^m\) is uniform by \cite[Proposition 6.3.4]{SW20}, and \([\varpi^\flat]\) is a non-zero divisor in \(B^m\). For pointwise effectivity, evaluating the rational \(\MCY\)-bundle on \((A^\flat,A^{+\flat})\) gives a vector bundle on \(\Yoo(A,A^+)\), with
\[\Yoo(A,A^+)\cap U^m=U.\]
Over \(U\), the natural action of \(\Gamma\) on its pullback agrees with the one on \(\MCE(U^m_\infty)\), so its fibres provide the effective descents required by condition (2). Thus \Cref{SenDescent} gives a finite projective \(B^m\)-module \(M^m\) whose associated v-bundle is \(\MCE|_{(U^m)^\diamond}\). By full faithfulness, these bundles glue for varying \(m\) to a vector bundle on \(\Yoc(A,A^+)\) with associated v-bundle \(\MCE\).
\end{proof}

\begin{myrem}
Note that we have only really used that we can test effectivity of the descent \emph{generically}; on the other hand, in \cite[Theorem 2.3.8]{PR}, it was only used that we can test effectivity \emph{pointwise}.
\end{myrem}
\subsection{\texorpdfstring{Case II: products of points with \(\varpi\) a non-zero divisor}
  {Case II: products of points with ϖ a non-zero divisor}}
In this section, we will show the following special case of \Cref{equivalenceThm}:

\begin{mythm}\label{mainThmAnalytic}
The functor 
\[\text{\(\OpP\)-bundles on } A^+
\longrightarrow\text{integral \(\MCY\)-bundles on }\Spd(A^+)\]
described in \Cref{equivalenceConst} is an equivalence of categories if 
\[A^+=\prod_{j\in J} V_j\]
where \(V_j\) is a valuation ring with algebraically closed field of fractions that is complete with respect to some non-zero divisor \(\varpi_j\). Here \(A^+\) is understood to carry the \(\varpi = (\varpi_j)_{j\in J}\)-adic topology.
\end{mythm}

We will call such an \(A^+\) a \emph{product of points}\footnote{The exact meaning of this term varies in the literature}. It is automatically integral perfectoid, and \((A^+[1/\varpi],A^+)\) is a perfectoid Huber pair.\par

For ease of notation, we will introduce the following temporary definition. It will only be used in this subsection:

\begin{mydef}
The category of ``\(\MCY_{[0,\infty]}\)-bundles on \(\Spd(A,A^+)\)" has objects consisting of triples \((\MCE_{(0,\infty]},\MCE_{[0,\infty)},\alpha)\) where 
\begin{enumerate}
    \item \(\MCE_{(0,\infty]}\) is a rational \(\MCY\)-bundle on \(\Spd(A^+)\),
    \item \(\MCE_{[0,\infty)}\) is an integral \(\MCY\)-bundle on \(\Spd(A,A^+)\),
    \item \(\alpha\) is an isomorphism of rational \(\MCY\)-bundles on \(\Spd(A,A^+)\) between the relevant restrictions of \(\MCE_{(0,\infty]}\) and \(\MCE_{[0,\infty)}\).
\end{enumerate}
The notion of morphism is the obvious one.
\end{mydef}

\begin{proof}(of \Cref{mainThmAnalytic}, assuming \Cref{lemmaAnchütz} and \Cref{lemmaComplicated})
There is a commutative diagram of categories

\begin{center}
\begin{tikzcd}
\text{vector bundles on }\Spec(\Ainf(A^+)) \arrow[d]  \arrow[r]  & \text{integral \(\MCY\)-bundles on }\Spd(A^+) \arrow[d] \\
{\text{vector bundles on }\Ycc(A,A^+)} \arrow[r] & \text{``\(\MCY_{[0,\infty]}\)-bundles on \(\Spd(A,A^+)\)"}.
\end{tikzcd}
\end{center}

We need to show that the top arrow is an equivalence.\par
The bottom arrow is defined by understanding a vector bundle on \(\Ycc(A,A^+)\) as a triple \((\MCE_{(0,\infty]},\MCE_{[0,\infty)},\alpha)\) where 
\begin{enumerate}
    \item \(\MCE_{(0,\infty]}\) is a vector bundle on \(\Yoc(A,A^+)\).
    \item \(\MCE_{[0,\infty)}\) is a vector bundle on \(\Yco(A,A^+)\).
    \item \(\alpha\) is an isomorphism of vector bundles on \(\Yoo(A,A^+)\) between the relevant restrictions of \(\MCE_{(0,\infty]}\) and \(\MCE_{[0,\infty)}\).
\end{enumerate}

Passing each of these three through the functors in \Cref{baseChangFunctors}, we get an object in ``\(\MCY_{[0,\infty]}\)-bundles on \(\Spd(A,A^+)\)", and \Cref{analyticLocuss} and \Cref{baseChangeFunctorsProperties} (3) show that this is an equivalence of categories.\par
The left vertical arrow is an equivalence of categories by \Cref{lemmaAnchütz}.\par
Finally, the right vertical arrow is defined by producing from an integral \(\MCY\)-bundle \(\MCE\) on \(\Spd(A^+)\) a triple \((\MCE_{(0,\infty]},\MCE_{[0,\infty)},\alpha)\) by, respectively, \Cref{baseChangFunctors} (4), pullback of \(\MCY\) bundles along maps of v-sheaves and compatibility of those two constructions, mentioned at the end of \Cref{baseChangFunctors}. We show that it is fully faithful in \Cref{lemmaComplicated}. It follows from the other two equivalences that it is essentially surjective and thus an equivalence.
\end{proof}

While \Cref{lemmaAnchütz} has been proved by Gleason, building on a result by Kedlaya, the proof of \Cref{lemmaComplicated} will take up the remainder of this subsection.

\begin{mythm}\label{lemmaAnchütz}\emph{(\cite[Proposition 2.1.17]{GleaPhD})}
    Let \(A^+\) be a product of points as described above. Then every vector bundle of constant rank on \(\Ycc(A,A^+)\) is free. In particular, as
    \[\Gamma(\Ycc(A,A^+),\MCO_{\Ycc(A,A^+)})=\Ainf(A^+),\]
    the base change functor from vector bundles of constant rank on \(\Spec(\Ainf(A^+))\) to vector bundles of constant rank on \(\Ycc(A,A^+)\) is a (non-exact) equivalence of categories.
\end{mythm}

\begin{mylemma}\label{sufficesForStructureSheaf}
    Let \(Z=\Spa(B,B^+)\) be an affinoid analytic adic space with an open subset \(\tilde Z\) such that
\[\Gamma(Z,\MCO_Z)\longrightarrow \Gamma(\tilde Z, \MCO_Z)\]
is a bijection (resp. injection). Then
\[\Gamma(Z,\MCE)\longrightarrow \Gamma(\tilde Z, \MCE)\]
is a bijection (resp. injection) for every vector bundle \(\MCE\) on \(Z\).
\end{mylemma}
\begin{proof}
    This is a standard argument in algebraic geometry. Since \(Z\) is affinoid and analytic,
    \[M\mapsto \left(U\mapsto (H^0(U,\MCO_Z)\otimes_{B} M)\right)\]
    defines an equivalence between finite projective \(B\)-modules and vector bundles on \(Z\). Cover \(\tilde Z\) by open affinoids \(U_i\). Then one computes:
\begin{align*}H^0(\tilde Z,\MCE)
&=\ker\left(\prod_i H^0(U_i,\MCE)\to \prod_{i,j} H^0(U_i\cap U_j,\MCE)\right)\\
&=\ker\left(\prod_i H^0(U_i,\MCO_Z)\otimes_{B} M\to \prod_{i,j} H^0(U_i\cap U_j,\MCO_Z)\otimes_{B} M\right)\\
&=\ker\left(\prod_i H^0(U_i,\MCO_Z)\to \prod_{i,j} H^0(U_i\cap U_j,\MCO_Z)\right)\otimes_{B} M\\
&=H^0(\tilde Z,\MCO_Z)\otimes_B M,
\end{align*}
where the third equality uses finite projectiveness of \(M\) to commute the tensor product with (possibly infinite) products and kernels. This immediately shows the statement for bijectivity; the statement for injectivity follows because tensoring with \(M\) is exact.
\end{proof}

\begin{myprop}\label{lemmaComplicated}
Let \(A^+\) be a product of points as in \Cref{mainThmAnalytic}. Then the functor
\[\emph{integral \(\MCY\)-bundles on }\Spd(A^+)\longrightarrow 
\emph{``\(\MCY_{[0,\infty]}\)-bundles on \(\Spd(A,A^+)\)"}
\]

is fully faithful.
\end{myprop}

\begin{proof} First note that, by passing to internal homs, it suffices to show equality of global sections, i.e., for \(\MCE\) an integral \(\MCY\)-bundle on \(\Spd(A^+)\), that the map
\[H^0_{\MCY_{\ZZ_p}}(\Spd(A^+),\MCE)\longrightarrow
H^0_{\MCY_{\QQ_p}}(\Spa(A^+),\MCE)
    \times_{H^0_{\MCY_{\QQ_p}}(\Spa(A,A^+),\MCE)}
    H^0_{\MCY_{\ZZ_p}}(\Spa(A,A^+),\MCE)\]

is an isomorphism. Here (and just for this proof), \(H^0_{\MCY_\MCS}(\cdot, \MCE), \MCS\in\{\ZZ_p,\QQ_p\}\) denote global sections of integral/rational \(\MCY\)-bundles in the obvious way.\par
Consider a two-truncated hypercover 
\[\bigsqcup_j \Spd(B_{j},B_{j}^+)\rightrightarrows \Spd(B,B^+)\longrightarrow \Spd(A^+)\] 
by perfectoid spaces.\footnote{\(\Spd(A^+)\) is quasi-compact, but not quasi-separated, so we cannot expect to choose \(\sqcup_j \Spd(B_{j},B_{j}^+)\) as a single affinoid perfectoid.}

This gives rise to a diagram of v-stacks

\begin{center}
    \begin{tikzcd}
    {\MCY_1}\coloneqq{\MCX_1} \times_{\MCX} \MCY \arrow[d, shift right] \arrow[d, shift left] \arrow[r] 
    & \MCX_1\coloneqq (\sqcup_j\Spd(B_{j},B_{j}^+)) \times \Spd(\ZZ_p)\arrow[d, shift left] \arrow[d, shift right] 
    \\  {\MCY_0}\coloneqq\MCX_0 \times_{\MCX} \MCY \arrow[d] \arrow[r]
    & { \MCX_0 \coloneqq \Spd(B,B^+) \times \Spd(\ZZ_p)} \arrow[d]         \\\MCY\coloneqq{\Ycc(A,A^+)^\diamond} \arrow[r]     
    & \MCX\coloneqq \Spd(A^+)\times\Spd(\ZZ_p).                                 
    \end{tikzcd}
\end{center}

Write 
\begin{align*}
    X_0&=\Yco(B,B^+)\\
    X_1&=\bigsqcup_j \Yco(B_j,B^+_j)\\
    Y_i& = X_i-V(p,\varpi).
\end{align*}
Then \(\MCX_i=X_i^\diamond\), \(\MCY_i=Y_i^{\diamond}\).\footnote{Note for the definition of \(\Yco(B,B^+)\) that \(\varpi\) is a pseudouniformizer of \(A^+\), not of \((B,B^+)\).}\par
Suppose we could find a choice of \((B,B^+),(B_j,B_j^+)_{j\in J}\) such that:
\begin{enumerate}
        \item The map \[H^0(X_0,\MCE)\longrightarrow H^0(Y_0, \MCE)\] is an isomorphism.
        \item The map \[H^0(X_1,\MCE)\longrightarrow H^0(Y_1, \MCE)\] is injective.
\end{enumerate}

Then bijectivity of the map of global sections from the beginning of the proof would follow. Indeed, there are bijections:

\begin{enumerate}[left=1.1cm]
    \item[] \(H^0_{\MCY_{\ZZ_p}}(\Spd(A^+),\MCE)\)
    \item[\(\cong\)] \(\left\{a\in H^0(X_0,\MCE) \text{ such that the two pullbacks to \(X_1\) agree} \right\}\)
    \item[\(\cong\)] \( \left\{a\in H^0(Y_0,\MCE) \text{ such that the two pullbacks to \(Y_1\) agree} \right\}\)
    \item[\(\cong\)] the set of all \((a,b)\in H^0(Y_0-V(p),\MCE) \times H^0(Y_0-V(\varpi),\MCE)\) such that:
    \begin{itemize}
        \item The two pullbacks of \((a,b)\) to \(H^0(Y_1-V(p),\MCE)\times H^0(Y_1-V(\varpi),\MCE)\) agree
        \item \(a|_{Y_0-V(p\varpi)}=b|_{Y_0-V(p\varpi)}\) 
    \end{itemize}
    \item[\(\cong\)] \(H^0_{\MCY_{\QQ_p}}(\Spa(A^+),\MCE)
    \times_{H^0_{\MCY_{\QQ_p}}(\Spa(A,A^+),\MCE)}
    H^0_{\MCY_{\ZZ_p}}(\Spa(A,A^+),\MCE)\),
\end{enumerate}
where the second isomorphism would follow from the properties (1) and (2) above. Here (and in the rest of the proof) we slightly abuse notation by also writing \(\MCE\) for the vector bundle on the analytic adic space \(X_0=\Yco(B,B^+)\) obtained by evaluating the integral \(\MCY\)-bundle at \(\Spd(B,B^+)\).\par
What is left to show now is that there are indeed \((B,B^+),(B_j,B_j^+)_{j\in J}\) such that the \(X_i,Y_i\) satisfy conditions (1) and (2) above. \par
Condition (2) actually holds for any choice of \((B_j,B_j^+)_{j\in J}\); it does not even have to be a cover. It is enough to prove injectivity of the map of sections of \(\MCE\) associated to
\[Y_0^j\coloneqq(\Spa(B_j,B_j^+)\dot\times \Spa(\ZZ_p))-V(p,\varpi)\longrightarrow X_0^j\coloneqq\Spa(B_j,B_j^+)\dot\times \Spa(\ZZ_p) \]
for every \(j\) individually. Choose a pseudouniformizer \(\varpi_{B_j}^{}\in B_j^+\) defining a continuous map \(\kappa\colon X_0^j=\Yco(B_j,B^+_j)\to [0,\infty)\) as in \cite[Figure 5 in section 12.1]{SW20} and cover \(X_0^j=\Yco(B_j,B^+_j)\) by rational open subsets
\[\MCY_{[0,1]}=U_{p\leq \varpi^{}_{B_j}\neq 0}, \qquad\MCY_{[1,n]}= U_n\coloneqq U_{\varpi^{}_{B_j}\leq p\neq 0, \;p^n\leq\varpi^{}_{B_j} \neq 0},\quad n\in \ZZ_{>0},\]
where we write \(\MCY_I\) as a shorthand for \(\MCY_I(B_j,B_j^+)\). Now \(Y_0^j\subset X_0^j\) contains \(X_0^j-V(p)=\Yoo\), so for each \(\MCY_{[1,n]}\), the inclusion \(Y_0^j\cap \MCY_{[1,n]} \hookrightarrow \MCY_{[1,n]}\) is an isomorphism and it suffices to show injectivity on sections of \(\MCE\) along \(\MCY_{[0,1]}\cap Y_0^j\to \MCY_{[0,1]}\). But \(\MCY_{[0,1]}\) is affinoid analytic, given as 
\[\MCY_{[0,1]}=\Spa\left(W(B^+_j)\langle  \frac{p}{\varpi^{}_{B_j}}\rangle[1/\varpi_{B_j}]\right),\]
so we can use \Cref{sufficesForStructureSheaf} to reduce to \(\MCE\) being the structure sheaf. Now consider the composition
\[\MCY_{[1/2,1]}
\hookrightarrow \MCY_{(0,1]} = \MCY_{[0,1]} -V(p)
\hookrightarrow \MCY_{[0,1]}-V(p,\varpi) = \MCY_{[0,1]} \cap Y_0^j
\hookrightarrow \MCY_{[0,1]},\]
giving rise to maps of global sections in the other direction. We want to show injectivity of the map of global sections corresponding to the last inclusion, and so it suffices to show injectivity of the map of global sections corresponding to the whole composition. But this can be easily calculated by writing both sides of
\[W(B_j^+)\langle \frac{p}{\varpi^{}_{B_j}}\rangle[1/\varpi_{B_j}] \hookrightarrow W(B_j^+)\langle \frac{p}{\varpi_{B_j}^{}}, \frac{\varpi_{B_j}^2}{p} \rangle[1/\varpi_{B_j}]\]
as (localizations of) rings of Witt vectors with some convergence condition.\par
Let's now have a look at the choice of \((B_0,B^+_0)\) and condition (1) above. Define 
\[\tilde A^+\coloneqq A^+[[t^{1/p^\infty}]].\]
Then \(\Spd(\tilde A^+)\to \Spd(A^+)\) is a v-cover of v-sheaves. Now \(\Spa(\tilde A^+)\) is not analytic, but \(\Spa(\tilde A^+)-V(\varpi,t)\) is, and it still covers \(\Spd(A^+)\) --- this is exactly as in the proof of \Cref{coverByRepresentable}, just with added p-power roots of \(t\). Now
\[U\coloneqq\Spa(\tilde A^+)-V(\varpi,t)\]
has a rational open cover by the two affinoid perfectoid subsets \(U_{t\leq \varpi\neq 0}\) and \(U_{\varpi\leq t\neq 0}\). Define
\[\Spa(B_0,B_0^+) \coloneqq (U_{t\leq \varpi\neq 0} \sqcup U_{\varpi\leq t\neq 0})\]
so that
\[X_0=\Yco(B_0,B^+_0)=U_{t\leq \varpi\neq 0}\dot\times\Spa(\ZZ_p)\sqcup U_{\varpi\leq t\neq 0}\dot\times \Spa(\ZZ_p).\]
Now as \(Y_0=X_0-V(p,\varpi)\) and \(\varpi\) vanishes nowhere on 
\[U_{t\leq \varpi\neq 0}\dot\times\Spa(\ZZ_p)\subset X_0,\]
this space is already contained in \( Y_0\). Thus, similarly to what we did for \(\MCY_{[1,n]}(B_j,B^+_j)\) above, it is enough to show equality of global sections of \(\MCE\) along
\[(U_{\varpi\leq t\neq 0}\dot\times \Spa(\ZZ_p))-V(p,\varpi)\hookrightarrow U_{\varpi\leq t\neq 0}\dot\times \Spa(\ZZ_p).\]
The right side is not affinoid, but rather it has a cover by rational open affinoid subsets
\[\Spa(C_n,C_n^+)\coloneqq (U_{\varpi\leq t\neq 0}\dot\times \Spa\ZZ_p)_{p^n\leq [t]\neq 0},\]
where \([t]\) denotes the Teichmüller lift, so we need to show that the restriction map
\[\Gamma(\Spa(C_n,C_n^+),\MCE)\longrightarrow \Gamma(\Spa(C_n,C_n^+)-V(p,[\varpi]),\MCE)\]
is a bijection for each \(n\). For readability, let's fix \(n\), write \(q=p^n\) and omit \(n\) from the index. As taking Witt vectors commutes with completions, one computes
\[C^+=\Ainf(\tilde A^+)\langle \frac{[\varpi]}{[t]},\frac{q}{[t]} \rangle,\qquad C=C^+[[t]^{-1}].\]
As \([t]\) is invertible in \(C\), one has \(V(q,[\varpi])=V(q/[t],[\varpi]/[t])\). Setting \(a=q/[t]\) and \(b=[\varpi]/[t]\), this places us in the situation of \Cref{codimTwoSubset}. Indeed, \((C,C^+)\) is a complete sheafy Tate-Huber pair, since it represents the rational affinoid open just described. Both \(q/[t]\) and \([\varpi]/[t]\) are non-zero divisors and as \(([t],q,[\varpi])\) and \(([t],[\varpi],q)\) are clearly regular sequences in \(B^+\coloneqq\Ainf(\tilde A^+)\), \Cref{regularSequenceInLocalization} below implies that \((q/[t], [\varpi]/[t])\) and \(([\varpi]/[t],q/[t])\) are regular sequences in \(C^+/([t])\).
\end{proof}

\begin{mylemma}\label{regularSequenceInLocalization}
    Let \(B^+\) be a topological ring carrying the \((t,a_1,\dots, a_n)\)-adic topology for \((t, a_1, \dots, a_n)\) a regular sequence in \(B^+\). Then \((\frac{a_1}{t},\dots,\frac{a_n}{t})\) is a regular sequence in 
    \[ B^+[ \frac{a_1}{t},\dots,\frac{a_n}{t}]^{\wedge}/(t)\]
\end{mylemma}

\begin{proof}
One has
\begin{align*}
        B^+[ \frac{a_1}{t},\dots,\frac{a_n}{t}]^{\wedge}
    &=(B^+[X_1,\dots,X_n]/I)^{\wedge (t,a_1,\dots,a_n)}\\
    &=(B^+[X_1,\dots,X_n]/I)^{\wedge t}
\end{align*}
for
\[I=(tX_1- a_1,\dots, tX_n-a_n),\]
with the two completions agreeing because in the quotient ring \(t\) divides the \(a_i\). Taking the quotient ring of this by \(t\), one gets
\[(B^+/(t, a_1,\dots a_n))[X_1,\dots,X_n],\]
in which \((X_1,\dots,X_n)=(\frac{a_1}{t},\dots,\frac{a_n}{t})\) is clearly a regular sequence.
\end{proof}

\begin{myprop}\label{codimTwoSubset}
    Let \((C,C^+)\) be a complete sheafy Tate-Huber pair with pseudouniformizer \(t\) and consider an ideal generated by two elements \(a,b\in C^+\) such that there exists a ring of definition \(C_0\) in which \((t,a,b)\) and \((t,b)\) form regular sequences (in the classical sense, e.g., \cite[0AUH]{Stacks}). Then restriction induces an isomorphism
    \[\Gamma(\Spa(C,C^+),\MCE)\xrightarrow{\sim}\Gamma(\Spa(C,C^+)-V(a,b),\MCE)\]
    for any vector bundle \(\MCE\) on \(\Spa(C,C^+)\).
\end{myprop}

\begin{proof}
Using \Cref{sufficesForStructureSheaf}, we can assume that \(\MCE\) is the structure sheaf. Define \(Z\coloneqq\Spa(C,C^+)\) and \(\tilde Z\coloneqq Z-V(a,b)\) and cover it by the rational subsets\footnote{Note that the condition \(a\leq b\neq 0\) alone need not define a rational subset, as \((a,b)C_0\subset C_0\) need not be open; see \cite[Definition 3.1.1]{SW20}.}
\[U_n\coloneqq Z_{t^n\leq a\neq 0},\qquad V_n\coloneqq Z_{t^n\leq b\neq 0},\quad n\in \ZZ_{>0},\]
which have local sections
\[\Gamma(U_n,\MCO_Z)=C\left\langle\frac{t^n}{a}\right\rangle,\qquad
\Gamma(V_n,\MCO_Z)=C\left\langle\frac{t^n}{b}\right\rangle.
\]
Now by the sheaf condition and the increasing open cover \(\tilde Z=\bigcup_n(U_n\cup V_n)\), it suffices to show that
\[C\left\langle\frac{t^n}{a}\right\rangle
\cap_{C\left\langle\frac{t^n}{a},\frac{t^n}{b}\right\rangle}
C\left\langle\frac{t^n}{b}\right\rangle=C\]
for every \(n\), with the intersection symbol understood as an equalizer along the natural maps and equality interpreted as the natural map from right to left being an isomorphism.
We proceed in several reduction steps:
\begin{enumerate}
    \item It would be enough to show that
    \[C\left\langle\frac{t}{a}\right\rangle
    \cap^{}_{C\left\langle\frac{t}{a},\frac{t}{b}\right\rangle} C\left\langle\frac{t}{b}\right\rangle = C\]
    or equivalently,
    \[\Gamma(U_1\cup V_1,\MCO_Z)=C,\]
    which somewhat surprisingly already holds true.\footnote{The relevant picture here is the following: We want to exclude the subset \(V(a,b)\), which should be thought of as a codimension-two subset, from \(Z\). We approximate \(V(a)\) by ``tubes'' of the form \(Z-U_n\), and similarly for \(V(b)\) and \(Z-V_n\). Now, surprisingly, just \(U_1\cup V_1\) is big enough to produce the same global sections as all of \(Z\), although its complement may be much larger than a codimension-two subset.}
    To see that this would indeed be enough, replace \(t\) by \(t^n\): the topology is unchanged, and \((t^n,a,b)\) and \((t^n,b)\) are still regular sequences by \cite[07DV]{Stacks}. Thus the same argument gives the assertion for every \(n\), without requiring injectivity of the restriction maps.

    \item Furthermore, it suffices to show that
    \[C\left\langle\frac{t}{a}\right\rangle_0
    \cap^{}_{C\left\langle\frac{t}{a},\frac{t}{b}\right\rangle_0} C\left\langle\frac{t}{b}\right\rangle_0=C_0\]
    for suitable rings of definition \((\cdot)_0\). Indeed, suppose we would find some
    \[f\in C\left\langle\frac{t}{a}\right\rangle
    \cap^{}_{C\left\langle\frac{t}{a},\frac{t}{b}\right\rangle}
    C\left\langle\frac{t}{b}\right\rangle
    =C\left\langle\frac{t}{a}\right\rangle_0[t^{-1}]
    \cap^{}_{C\left\langle\frac{t}{a},\frac{t}{b}\right\rangle_0[t^{-1}]} C\left\langle\frac{t}{b}\right\rangle_0[t^{-1}]\]
    not contained in \(C=C_0[t^{-1}]\). Then multiplication with \(t^n\) for \(n\) sufficiently large would produce an element in
    \[C\left\langle\frac{t}{a}\right\rangle_0
    \cap^{}_{C\left\langle\frac{t}{a},\frac{t}{b}\right\rangle_0}
    C\left\langle\frac{t}{b}\right\rangle_0\]
    not contained in \(C_0\). Injectivity likewise follows by exactness of localization at \(t\).

    \item Now let \(C_0\) be a ring of definition such that \((t,a,b)\) and \((t,b)\) are regular sequences. It is open and hence closed in the complete ring \(C\), so
    \[C_0=\varprojlim_n C_0/t^n.\]
    On the other side of the equation,
    \[C\left\langle\frac{t}{a}\right\rangle_0\coloneqq C_0[\frac t a]^{\wedge t}\]
    is a ring of definition, and similarly for the other two rings.
    Again we write each ring as a completion, which we can pull out of the equalizer:
    \[C\left\langle\frac{t}{a}\right\rangle_0
    \cap^{}_{C\left\langle\frac{t}{a},\frac{t}{b}\right\rangle_0}
    C\left\langle\frac{t}{b}\right\rangle_0
    =\varprojlim_n\left((C_0[\frac t a]/t^n)
    \cap_{C_0[\frac t a,\frac t b]/t^n}(C_0[\frac t b]/t^n)\right).\]
    Thus it is enough to show that for every \(n\geq1\)
    \[(C_0[\frac t a]/t^n)\cap_{C_0[\frac t a,\frac t b]/t^n}(C_0[\frac t b]/t^n)=C_0/t^n.\]

\end{enumerate}
We first do the case \(n=1\). We have presentations:
    \begin{align*}
    C_0[\frac t a]&=C_0[X]/(aX-t),\\
    C_0[\frac t b]&=C_0[Y]/(bY-t),\\
    C_0[\frac t a,\frac t b]&=C_0[X,Y]/(aX-t,bY-t).
    \end{align*}
	Indeed, there are obvious surjections from right to left which become isomorphisms after inverting \(t\). Thus, to prove that they are isomorphisms, it suffices to show that the rings on the right are \(t\)-torsion free. We will show this for the third ring, the other two cases are similar and easier. Writing \(R\coloneqq C_0/tC_0\), note first that \((aX,bY)\) is a regular sequence in \(R[X,Y]\) since both \((b)\) and \((a,b)\) are regular sequences in \(R\); this is clear for each step, except maybe regularity of \(bY\) inside
\[R[X,Y]/(aX)\cong R[Y]\oplus \bigoplus_{i>0}(R/(a))[Y]X^i,\]
where it follows by separately examining both summands. Now assume\[tH=0\in C_0[X,Y]/(aX-t,bY-t).\] Thus there are \(L,M\in C_0[X,Y]\) such that
    \[tH=(aX-t)L+(bY-t)M\in C_0[X,Y].\]
Reducing this modulo \(t\) and using that \((aX,bY)\) is a regular sequence in the quotient ring, we get for some \(N\in C_0[X,Y]\)
\[\overline L-bY\overline N = \overline M+aX\overline N=0\in R[X,Y].\]
Write
\[L' = L-bY N,\qquad  M'=M+aX N.\]
These are divisible by \(t\), and a calculation gives
    \[(aX-t)L'+(bY-t)M'=tH+t(bY-aX)N.\]
Now every term above is divisible by \(t\), and cancelling it shows \(H\in (aX-t, bY-t)\), so \(H=0\) in the quotient. This shows \(t\)-torsion freeness and thus the presentations above. Thus
    \begin{align*}
    C_0[\frac t a]/t&=R[X]/(aX),\\
    C_0[\frac t b]/t&=R[Y]/(bY),\\
    C_0[\frac t a,\frac t b]/t&=R[X,Y]/(aX,bY).
    \end{align*}
Comparing coefficients proves the case \(n=1\) of the equalizer from step (3).\par
The statement for general \(n\) follows by induction. A compatible pair modulo \(t^n\) reduces modulo \(t\) to a constant; subtract a lift in \(C_0\) and divide both components by \(t\). The quotients are compatible modulo \(t^{n-1}\) because \(C_0[\frac t a,\frac t b]\) is \(t\)-torsion-free, so induction proves surjectivity. Likewise, a constant mapping to zero is divisible by \(t\); cancellation in \(C_0[\frac t a],\allowbreak C_0[\frac t b]\) and induction give divisibility by \(t^n\), proving injectivity. Taking the limits in (3) and over the increasing cover \((U_n\cup V_n)_n\) finishes the proof.

\end{proof}
\subsection{Finishing the proof}
\label{secFinishing}

We need one more lemma and another proposition to reduce our main theorem to the two principal cases discussed earlier:

\begin{mylemma}\label{IisVarpi}
    Let \(A^+=\prod_i V_i\) be a product of valuation rings and let \(I\subseteq A^+\) be a finitely generated ideal containing \(p\). 
    
    \begin{enumerate}
        \item There is some \(\varpi\in A^+\) such that the \(I\)-adic topology agrees with the \(\varpi\)-adic one. 
        \item If we furthermore assume that the \(V_i\) have algebraically closed field of fractions, one can choose \(\varpi\) so that \(\varpi^p|p\).
    \end{enumerate}
\end{mylemma}

\begin{proof}
    Let \(g=\{a_1,\dots,a_n\}\) be a set of generators of \(I\) and assume without loss of generality \(a_1=p\). Then, by \(I\)-completeness, for any fixed \(V_i\), all \(a_j\in g\) whose image in \(V_i\) is non-zero generate the same topology on \(V_i\). For each \(V_i\), let \(a^{(i)}\in V_i\) be an element in the image of \(g\) with maximal valuation, i.e., so that for each \(j\in \{1,\dots,n\}\) one has \(a^{(i)}|a_j\) in \(V_i\), and define \(\varpi_i\coloneqq a^{(i)}\) (resp. \(\varpi_i\coloneqq(a^{(i)})^{1/p}\) in case (2)). Then the preceding consideration shows that the \(I\)-adic topology on \(\prod_i V_i\) agrees with the \((\varpi_i)_i\)-adic one. In case (2), as we have assumed that \(p\in g\) and that \(a^{(i)}\) has maximal valuation among all elements in \(g\), we also have \(\varpi^p|p\).
\end{proof}

\begin{myprop}\label{lemmaCoveringg}
    Let \(R^+\) be a complete topological ring carrying the \(I\)-adic topology for some finitely generated ideal \(I\) containing \(p\). Then there is a map of topological rings \(R^+\to A^+=A_{\disc}^+\times A_{\ana}^+\) where
\begin{enumerate}
\item \(A_{\disc}^+\) is a perfect ring of characteristic \(p\) that is a product of valuation rings carrying the discrete topology.
\item \(A_{\ana}^+\) is a product of points in the sense of the previous subsection, i.e., a product of valuation rings with algebraically closed fields of fractions, equipped with the \(\varpi\)-adic topology and complete with respect to it for some non-zero divisor \(\varpi\).
\item \(R^+\to A^+\) is an \(I\)-complete arc cover of rings in the sense of \Cref{defArcSite}.
\item \(\Spd(A^+)\to \Spd(R^+)\) is a v-cover of v-sheaves.
\item One has \[\Spd(A^+\otimes^{\wedge I}_{R^+} A^+)=\Spd(A^+)\times_{\Spd(R^+)} \Spd(A^+).\]
\end{enumerate} 
\end{myprop}

In particular, properties (1) and (2) imply that \(A^+\) is integral perfectoid.

\begin{proof}
Consider the set\footnote{The argument for finding the smallest valuation ring in each equivalence class just below shows that each such class has a representative which is \(\kappa\coloneqq \max(\omega, |R^+|)\)-small, so the equivalence classes indeed form a set, not a proper class.} \(J\) of equivalence classes of maps \(R^+\to V\) with \(V\) an \(I\)-complete valuation ring. In each of these equivalence classes we can find a smallest valuation ring by taking any representative \(v\colon R^+\to V\) and considering 
\[\Frac(\im(v))\cap V\subset \Frac(V),\]
with \(\im(v)\) the image of \(v\).\par
For every \(j\in J\) let \(V_j\) be the \(I\)-completion of the absolute integral closure of the associated smallest valuation ring as above. Then 
\[A^+\coloneqq \prod_{j\in J} V_j\]
is a product of valuation rings with algebraically closed field of fractions, and thus, by \Cref{IisVarpi}, there is some \(\varpi=(\varpi_j)_{j\in J}\in A^+\) such that \(\varpi^p\mid p\) and for which the \(I\)-adic topology agrees with the \(\varpi\)-adic one. As \(A^+\) is complete with respect to the former, it is also complete with respect to the latter.\par
Let 
\[J_{\disc}\coloneqq \{j\in J\mid \varpi_j=0\in V_j\}\]
and
\[J_{\ana}\coloneqq J - J_{\disc}\]
and set 
\[A^+_{\disc}\coloneqq \prod_{j\in J_{\disc}} V_j,
\qquad A^+_{\ana} = \prod_{j\in J_{\ana}}V_j.\]
Then the \(\varpi\)-adic (and hence also the \(I\)-adic) topology on \(A_{\disc}^+\) is clearly discrete, and \((\varpi_j)_{j\in J_{\ana}} \in A^+_{\ana}\) is a non-zero divisor. This shows that \(A^+\) satisfies (1) and (2). \par
For (3), note that every map \(R^+\to V\) to an \(I\)-complete valuation ring \(V\) (so in particular to one of rank one) factors through some \(V_i\), at least after exchanging \(V\) for its \(I\)-completed absolute integral closure, which is an extension of valuation rings. Via projection to the \(i\)-th component, we get a factorization over \(R^+\to A^+=\prod_j V_j\).\par
For (4) it suffices to show that for all maps to products of points in the adic spaces sense
\[(R^+,R^+)\longrightarrow\left((\prod_{j'\in J'} C_{j'}^+)[((\varpi_{j'})_{{j'}\in J'})^{-1}],\prod_{{j'}\in J'}C_{j'}^+\right)\eqqcolon (S,S^+)\]
with \(C_{j'}^+[1/\varpi_{j'}]\) algebraically closed there exists a lift as a map of Huber pairs

\begin{center}
\begin{tikzcd}
{( A^+, A^+)} \arrow[rd, dotted]
&    

\\
{(R^+,R^+)} \arrow[r] \arrow[u]                  &
{(S,S^+).}
\end{tikzcd}
\end{center}
By what we have seen in (3), we can lift every map \(R^+\to C^+_{j'}\) to a map of rings \(A^+\to C^+_{j'}\), producing a map of pairs of rings 
\[(A^+,A^+)\longrightarrow (S^+,S^+)\longrightarrow(S,S^+).\]
To see that this map is continuous, note that the topology on \(A^+\) is generated by the image of \(I\) under the map \(R^+\to A^+\). Hence we need to show that some power of the image of \(I\) under the composition \(R^+\to A^+\to S^+\) is contained in the ideal generated by \((\varpi_{j'})_{{j'}\in J'}\in S^+\). But this is exactly the condition for \(R^+\to S^+\) to be continuous.\par
Condition (5) follows from the fact that \((A^+\otimes^{\wedge I}_{R^+} A^+,A^+\otimes^{\wedge I}_{R^+} A^+)\) is a Huber pair, and it is straightforward to see that it satisfies the universal property of the fibre product on the right side. Here it is essential that \(A^+\) carries the \(I\)-adic topology.
\end{proof}

\begin{proof} (of \Cref{equivalenceThm}):
Consider a map \(R^+\to A^+\coloneqq A^+_{\disc} \times A^+_{\ana}\) from the previous proposition. Then \(A^+_{\disc}\) and \(A^+_{\ana}\) are of the form dictated in the subsections covering case I (\Cref{mainThmDiscrete}) and case II (\Cref{mainThmAnalytic}), respectively. As both integral \(\MCY\)-bundles and \(\OpP\)-bundles behave as expected with disjoint unions of schemes/small v-sheaves, we then get:
\begin{itemize}
    \item[] \(\{\OpP\)-bundles on \(A^+\}\)
    \item[\(\cong\)] \(\{\OpP\)-bundles on \(A_{\disc}^+\}\times \{\OpP\)-bundles on \(A_{\ana}^+\}\) 
    \item[\(\cong\)] \(\{\)integral \(\MCY\)-bundles on \(\Spd(A_{\disc}^+)\}\times \{\)integral \(\MCY\)-bundles on \(\Spd(A_{\ana}^+)\}\) 
    \item[\(\cong\)] \(\{\)integral \(\MCY\)-bundles on \(\Spd(A^+)=\Spd(A_{\disc}^+)\sqcup \Spd(A^+_{\ana})\}\),
\end{itemize}
Showing the equivalence for \(A^+\). To descend this result to \(R^+\), apply \Cref{lemmaCoveringg} to \(A^+\otimes^{\wedge I}_{R^+}A^+\) to obtain a cover
\[A^+\otimes^{\wedge I}_{R^+}A^+\longrightarrow A_1^+,\]
and again to obtain a cover
\[A_1^+\otimes^{\wedge I}_{R^+}A_1^+\longrightarrow A_2^+\]
to get a truncated cosimplicial object
\begin{equation*}\tag{*}\label{finishingHypercover}
    A^+ \rightrightarrows A^+_1  \rrrarrows A^+_2. \qquad 
\end{equation*}
Now by property (3) of \Cref{lemmaCoveringg} this is an \(I\)-complete arc hypercover of \(R^+\), and so by \Cref{arcDescentForCrysNew}, \(\OpP\)-bundles on \(R^+\) correspond to descent data on \eqref{finishingHypercover}. On the other hand, by properties (4) and (5) of \Cref{lemmaCoveringg}, \(\Spd(\cdot)\) applied to \eqref{finishingHypercover} is also a hypercover of \(v\)-sheaves and thus by \Cref{vDescentForVBundles} (2) integral \(\MCY\)-bundles on \(R^+\) correspond to descent data of integral \(\MCY\)-bundles on \eqref{finishingHypercover}. Thus, by the equivalence of the two categories for \(A^+\) (and with the same argument also for \(A_1^+\) and \(A_2^+\)), they also agree on \(R^+\).\par
This leaves us with exactness, i.e., that a sequence of \(\OpP\)-bundles on \(A^+\) is exact if and only if the associated sequence of \(\MCY\)-bundles is. By \Cref{lemmaCoveringg} in combination with \Cref{exactnessArcLocallyGeneral} (2) and \Cref{exactnessIsVLocal}, we can assume that our ring still has the same form as above. By \Cref{exactnessArcLocallyGeneral} (1), exactness of \(\OpP\)-bundles is the same as exactness of finite projective \(\Ainf(A^+)\)-modules. Using the product-of-points form of \(A^+\), a short exact sequence of such modules splits, so its associated sequence of integral \(\MCY\)-bundles is short exact. It remains to prove the converse, for which we again use the special shape of \(A^+\) to treat the discrete and analytic factors separately.\par 
\emph{The discrete case \(A^+=A^+_{\disc}\).}\footnote{since we don't know \((W(A^+_{\disc}),W(A^+_{\disc}))\) to be sheafy, see \Cref{discreteWittSheafiness} we have to do this step by hand.} Let 
\[ M_1\longrightarrow M_2\longrightarrow M_3\]
 be finite projective \(\Ainf(A^+)=W(A^+)\)-modules whose associated integral \(\MCY\)-bundles form a short exact sequence. By the full faithfulness proved above, these module maps form a complex. For each characteristic-\(p\) perfectoid Huber pair \((B,B^+)\) with a map
\[\Spd(B,B^+)\longrightarrow \Spd(A^+),\]
restrict the associated bundles to the \(p=0\)-locus along
\[\Spa(B,B^+)\hookrightarrow\Yco(B,B^+),\]
obtaining a short exact sequence of vector bundles. By compatibility of the different functors between bundles, it is the sequence associated to
\[M_i\otimes_{W(A^+)}B=(M_i/pM_i)\otimes_{A^+} B.\]
Letting \((B,B^+)\) vary, the v-bundles associated to \(M_i/pM_i\) form a short exact sequence of v-bundles on \(\Spd(A^+)\). By the exact tensor equivalence of \cite[Corollary 5.3]{Kim}, the sequence
\[0\longrightarrow M_1/pM_1\longrightarrow M_2/pM_2\longrightarrow M_3/pM_3\longrightarrow0\]
of \(A^+\)-modules is also exact. Using that \(p\) lies in the Jacobson ideal of \(W(A^+)\), the result now follows using the Nakayama-splitting trick already seen in the proof of \Cref{exactnessArcLocally}.\par
{\emergencystretch=1em \emph{The analytic case \(A^+=A^+_{\ana}\).} Now put \(A^+=A^+_{\ana}\), and consider finite projective \(\Ainf(A^+)\)-modules whose associated integral \(\MCY\)-bundles form a short exact sequence. As in the discrete case, full faithfulness ensures that the module maps form a complex. The space \(\Spa(\Ainf(A^+))\) exists as an adic space by \Cref{liuAinf}. We can cover \(\Spd(A^+)\) by a representable \(\Spa(B,B^+)\) as in \Cref{coverByRepresentable}.\par}
\[\Spd(B,B^+)\times \Spd(\ZZ_p)\longrightarrow \Spd(A^+)\times \Spd(\ZZ_p)\]
is still a cover, so 
\begin{align*}
    |\Spa(B,B^+)\dot\times \Spa(\ZZ_p)|&\isom|\Spd(B,B^+)\times\Spd(\ZZ_p)|\\
    &\twoheadrightarrow |\Spd(A^+)\times \Spd(\ZZ_p)|\\
    &= |\left(\Spa(A^+)\dot\times \Spa(\ZZ_p)\right)^{\diamond}|\\
    &\twoheadrightarrow |\Spa(A^+)\dot\times \Spa(\ZZ_p)|\\
    &=|\Spa(\Ainf(A^+))|,
\end{align*}
with the second twoheadarrow a surjection by \cite[Proposition 18.2.2]{SW20}\footnote{Cf. the argument in \Cref{YcoSurjectiveness}}. By \Cref{exactnessYBundleEquivalentDefinition}, the evaluated complex of vector bundles on
\[\Spa(B,B^+)\dot\times \Spa(\ZZ_p)=\Yco(B,B^+)\]
is short exact. Hence \Cref{exactnessAdicLocalCorollary} implies short exactness of the associated vector bundles on \(\Spa(\Ainf(A^+))\). By \Cref{adicPointsFromClassicalPoints}, each maximal ideal of \(\Ainf(A^+)\) is the support of an adic point. The module complex is therefore exact over every residue field at a maximal ideal, and the finite-projective criterion used in \Cref{exactnessArcLocally} implies that it is exact.
\end{proof}

\section{Adding a group action and Frobenius structure}
\subsection{Group action}
For the remainder of the paper, \(\MCG\) will denote an affine flat group scheme of finite type over \(\ZZ_p\) (a generalization to, e.g., finite extensions of \(\ZZ_p\) is straightforward) and \(\Rep(\MCG)\) the category of \(\MCG\)-representations on finite free \(\ZZ_p\)-modules.

A \(\MCG\)-bundle will always be defined as some exact tensor functor. This is justified by the following result:

\begin{mythm}\emph{\cite[Example 1.7 (1), Proposition 1.6]{Wed23}}
Let \(A\) be a regular noetherian ring of dimension \(\leq 1\) and let \(\MCG\) be any affine flat group scheme of finite type over \(A\). Then for any \(A\)-scheme \(X\), there is an equivalence of groupoids between \(\MCG\)-torsors on \(X\) and exact tensor functors \(\Rep(\MCG)\to \emph{Vec}(X)\).
\end{mythm}

Let's use this to define \(\MCG\)-structure versions of both our formal and adic objects.

\begin{mydef}\label{defGCrys}
A \emph{\(\MCG\)-\(\OpP\)-bundle} over an \(I\)-complete ring \(R^+\) is given by the following equivalent pieces of data:
\begin{enumerate}
    \item A rule sending each \((R^+\to S^+)\in (R^+)_{\Prism}^{\perf}\) to a \(\MCG\)-bundle on \(\Ainf(S^+)\), compatible with pullbacks.
    \item An exact tensor functor from \(\Rep(\MCG)\) to the category of \(\OpP\)-bundles over \(R^+\), with \((R^+)_{\Prism}^{\perf}\) being given the \(I\)-complete arc topology.
\end{enumerate}
\end{mydef}

The equivalence follows from \Cref{exactnessArcLocallyGeneral}(1). By \Cref{arcFinerThanFlat} one could also replace the \(I\)-complete arc topology by the \((p, \overline I)\)-complete flat topology, see \Cref{defPrismaticSite} for the definition of these topologies.

\begin{mydef}\label{defGYBundle}
An integral (resp. rational) \emph{\(\MCG\)-\(\MCY\)-bundle} on a small v-sheaf \(\MCX\) is one of the following equivalent pieces of data:
	\begin{enumerate}
		\item Associating to any map from a representable \(\Spa(B,B^+)\to \MCX\) a \(\MCG\)-bundle on 
        \[\Yco(B,B^+)=\Spa(B,B^+)\dot\times\Spa(\ZZ_p)\]
        resp.
        \[\Yoo(B,B^+)=\Spa(B,B^+)\dot\times\Spa(\QQ_p),\]
        compatible with pullbacks.
        \item An exact tensor functor from \(\Rep(\MCG)\) to the category of integral/rational \(\MCY\)-bundles on \(\MCX\).
	\end{enumerate}
\end{mydef}

Here, the equivalence follows from \Cref{exactnessVBundleYBundleEquivalence}. Using characterization (2) in both \Cref{defGCrys} and \Cref{defGYBundle}, the exactness of the equivalence in \Cref{equivalenceThm} immediately implies the following:

\begin{mycor}\label{GEquivalenceTheorem}
Let \(R^+\) be as in \Cref{equivalenceThm}. For any affine flat group scheme \(\MCG/\ZZ_p\) of finite type, \Cref{equivalenceConst} defines an equivalence of groupoids
\[\text{\emph{\(\MCG\)-\(\OpP\)-bundles on }}R^+
\longrightarrow\text{\emph{integral \(\MCG\)-\(\MCY\)-bundles on }}\Spd(R^+).\]
\end{mycor}
\subsection{Frobenius structure}
Recall that each prism comes equipped with a Frobenius endomorphism \(\phi\). These endomorphisms assemble to a Frobenius endomorphism of sheaves of rings \(\OpP\to\OpP\), also denoted \(\phi\).

\begin{mydef}
Let \(R^+\) be an \(I\)-complete ring. A \emph{perfect-prismatic \(F\)-crystal} on \(R^+\) is one of the following equivalent pieces of data:
\begin{enumerate}
    \item A rule associating to each \((R^+\to S^+)\in (R^+)_{\Prism}^{\perf}\) a BKF-module on \(S^+\), compatible with pullbacks.
    \item An \(\OpP\)-bundle \(\MCE\) together with a Frobenius-linear isomorphism 
    \[\Phi\colon(\phi^*\MCE)\otimes_{\OpP}\OpP[1/d]\isom \MCE\otimes_{\OpP}\OpP[1/d].\]
\end{enumerate}
\end{mydef}

These definitions agree by quasi-compactness of \(R^+\) in \((R^+)_{\Prism}^{\perf}\).

\begin{myrem}
    If \(R^+\) is an integral perfectoid ring, the category of perfect-prismatic \(F\)-crystals is equivalent to the category of BKF-modules on \(R^+\).
\end{myrem}

\begin{myprop}\label{arcDescentForFCrys}
Let \(R^+\) be an \(I\)-complete ring for some finitely generated ideal \(I\subset R^+\) containing \(p\). Sending an \(I\)-complete \(R^+\)-algebra \(S^+\) to the category of perfect-prismatic \(F\)-crystals on \(S^+\) satisfies \(I\)-complete arc descent.
\end{myprop}

\begin{proof}
By \Cref{arcDescentForCrysNew}, associating to any \(I\)-complete \(R^+\)-algebra \(S^+\) the category of \(\OpP\)-bundles on \(S^+\) satisfies \(I\)-complete arc descent. In particular, for any \(\OpP\)-bundle \(\MCE\) on \(S^+\), the rule
\[(S'^+)\mapsto \MCE(S'^+)=\Hom(\MCO_{\Prism}^{\perf}(S'^+),\MCE(S'^+))\]
is again an \(I\)-complete arc sheaf on \((S^+)^{\perf}_{\Prism}\). Furthermore, \((S^+)_{\Prism}^{\perf}\) trivially satisfies the condition in \cite[0739]{Stacks} (all covers are finite and fibre products exist by \Cref{defArcSite}), so
\[\MCE[1/d]=\colim (\MCE\xrightarrow{\cdot d}\MCE\xrightarrow{\cdot d} \dots)\]
is again a sheaf. Now given a descent datum of perfect-prismatic \(F\)-crystals over some \(I\)-complete arc hypercover 
\[S^+\longrightarrow S_1^+\rightrightarrows S_2^+\rrrarrows S_3^+,\]
we first descend the underlying \(\OpP\)-bundles on \(S^+_i\) to an \(\OpP\)-bundle \(\MCE\) on \(S^+\) via \Cref{arcDescentForCrysNew}. Then the extra Frobenius-linear structure of the perfect-prismatic \(F\)-crystal on \(S^+_i\) can be described as an (invertible) element in 
\[((\phi^*\MCE)^{\vee}\otimes_{\MCO_{\Prism}^{\perf}} \MCE)[1/d](S_i^+).\]
But this forms a sheaf by our previous observation, so we can descend these sections to a section in 
\[((\phi^*\MCE)^{\vee}\otimes_{\MCO_{\Prism}^{\perf}} \MCE)[1/d](S^+)\]
giving the \(\OpP\)-bundle \(\MCE\) the structure of a perfect-prismatic \(F\)-crystal.
\end{proof}

Let's turn our attention to the adic side again. As in \Cref{notationVSheaves}, we will use \(X_\diamond\) to denote a v-sheaf with a structure map to \(\Spd(\ZZ_p)\) (not necessarily coming from an adic space).

\begin{mydef}(Scholze-Weinstein shtukas)\label{DefShtukas}
\begin{enumerate}
    \item By a \emph{shtuka}\footnote{In \cite{SW20}, this would be called a \(\GL_n\)-shtuka, but we say just shtuka both for brevity and to emphasize the difference to the \emph{groupoid} of \(\GL_n\)-shtukas defined below.} over some affinoid perfectoid \(\Spa(R,R^+)\) of characteristic \(p\) we mean the following data:
		\begin{itemize}
    		\item An untilt \((R^\sharp,R^{+\sharp})\) over \(\ZZ_p\),
    		\item A vector bundle \(\MCE\) on \(\Yco(R,R^+)\),
    		\item An isomorphism of vector bundles 
            \[\Phi\colon(\phi^*\MCE)|_{\Yco-\Spa(R^\sharp,R^{+\sharp})}
            \longrightarrow\MCE|_{\Yco-\Spa(R^\sharp,R^{+\sharp})}\]
            that is meromorphic with respect to a distinguished element \(\xi\) cutting out 
            \[\Spa(R^\sharp,R^{+\sharp})\subset \Yco(R,R^+).\]
		\end{itemize} 
    \item A \emph{shtuka with fixed leg} over a perfectoid Huber pair \((B,B^+)\), not necessarily of characteristic \(p\), is a shtuka over \((B^\flat,B^{+\flat})\) where the untilt \((B^{\flat\sharp},B^{+\flat\sharp})\) that comes as part of the datum is given by \((B,B^+)\).
\end{enumerate}
\end{mydef}

\begin{mydef}
    Let \(X_\diamond\) be a v-sheaf together with a map to \(\Spd(\ZZ_p)\). A \emph{family of shtukas with fixed leg} or just a \emph{shtuka} for short, over \(X_\diamond\) is given by a rule associating to any map from a representable \(\Spa(B,B^+)\to X_\diamond\) with associated untilt \(\Spa(B^\sharp,B^{+\sharp})\) a shtuka with fixed leg on \(\Spa(B^\sharp,B^{+\sharp})\), compatible with pullbacks.
\end{mydef}

\begin{myrem}\label{niceShtukaDescription}
If \(X_\diamond=X^\diamond\) for an adic space \(X\), the above definition simplifies to a rule associating to any map \(\Spa(B,B^+)\to X\) from a (not necessarily characteristic \(p\)) affinoid perfectoid a shtuka with fixed leg on \((B,B^+)\).
\end{myrem}

\begin{myprop}\label{vDescentForShtukas}
Associating to a small v-sheaf \(X_\diamond\) over \(\Spd(\ZZ_p)\) the category of shtukas with fixed leg on \(X_\diamond\) satisfies descent along v-covers that are compatible with the structure maps to \(\Spd(\ZZ_p)\).
\end{myprop}

\begin{proof}
For shtukas without fixed leg and affinoid perfectoids, this is the well-known statement shown in \cite{SW20} that shtukas descend along v-covers of affinoid perfectoids. If we assume the cover to be compatible with the maps to \(\Spd(\ZZ_p)\), the descended shtuka automatically respects the fixed leg structure (as the map to \(\Spd(\ZZ_p)\) just describes the position of the leg). Finally, the extension to general small v-sheaves is formal and works in the same way as, e.g., \Cref{arcDescentForCrysNew} or \Cref{vDescentForVBundles}.\par
Alternatively, one can show the statement directly using \Cref{vDescentForVBundles}(2) and the argument from \Cref{FEquivalenceThm}.
\end{proof}


\begin{myconst}\label{FEquivalenceConst}
The functor in \Cref{equivalenceConst} extends to a functor
\[\text{perfect-prismatic \(F\)-crystals on }R^+
\longrightarrow\text{shtukas with fixed leg on }\Spd(R^+).\]
Indeed, a map from a representable
\[\Spa(B,B^+)\longrightarrow\Spd(R^+)\]
supplies an untilt \((B^\sharp,B^{+\sharp})\) and a map \(R^+\to B^{+\sharp}\). The map
\[\Yco(B,B^+)\longrightarrow \Spa(\Ainf(B^{+\sharp}))\]
then commutes with Frobenius morphisms and sends the distinguished element\[d\in\Ainf(B^{+\sharp})=W(B^+)\]to a distinguished element \(\xi\) in the global sections of \(\Yco(B,B^+)\), so that we can pull back our Frobenius-linear isomorphism.
\end{myconst}

\begin{mythm}\label{FEquivalenceThm}
Let \(R^+\) be a complete topological ring carrying the \(I\)-adic topology for some finitely generated ideal \(I\subset R^+\) containing \(p\). Then the functor just defined is an exact tensor equivalence.
\end{mythm}

\begin{proof}
Given \Cref{equivalenceThm}, we need to show that for an \(\OpP\)-bundle \(\MCE\) on \(R^+\), every compatible collection of isomorphisms 
\[\Phi_{(B,B^+)}\colon(\phi^*\MCE)|_{\Yco(B,B^+)}[1/\xi]\longrightarrow \MCE|_{\Yco(B,B^+)}[1/\xi]\]
running over maps from representables \(\Spa(B,B^+)\to \Spd(R^+)\), with associated untilts \((B^\sharp,B^{+\sharp})\), descends to a compatible family of isomorphisms
\[(\phi^*\MCE)|_{\Ainf(B^{+\sharp})}[1/d_{B^{+\sharp}}]\longrightarrow \MCE|_{\Ainf(B^{+\sharp})}[1/d_{B^{+\sharp}}].\]
We can assume by \Cref{lemmaCoveringg} and the descent results \Cref{arcDescentForFCrys} and \Cref{vDescentForShtukas} that \(R^+\) is integral perfectoid. Fix a distinguished element \(d\in\Ainf(R^+)\), and on each test object take \(\xi\) to be its image. Cover \(\Spd(R^+)\) by a single affinoid perfectoid \(\Spa(B_0,B_0^+)\), possible by \Cref{coverByRepresentable}. Then by the meromorphicity condition in \Cref{DefShtukas}, there is some \(n\) such that multiplying the map
\[\Phi_{ (B_0, B_0^+)}\colon\phi^*\MCE|_{\Yco (B_0, B_0^+)-V(\xi)}\longrightarrow \MCE|_{\Yco (B_0, B_0^+)-V(\xi)}\]
by \(\xi^n\) yields a map
\[\xi^n\Phi_{ (B_0, B_0^+)}\colon\phi^*\MCE|_{\Yco (B_0, B_0^+)}\longrightarrow \MCE|_{\Yco (B_0, B_0^+)}.\]
Using compatibility with base change, one sees that in fact for the same \(n\) and each map \(\Spa(B,B^+)\to \Spd(R^+)\) we get
\[\xi^n\Phi_{ (B, B^+)}\colon\phi^*\MCE|_{\Yco (B, B^+)}\longrightarrow \MCE|_{\Yco (B, B^+)}.\]
But the collection of these maps for varying \((B,B^+)\) is just a map of integral \(\MCY\)-bundles \(\phi^*\MCE\to \MCE\), which, by \Cref{equivalenceThm}, descends to a map
\[d^n\Phi\colon \phi^*\MCE|_{\Ainf(R^+)}\longrightarrow \MCE|_{\Ainf(R^+)}\]
of \(\OpP\)-bundles. Applying the same procedure to \(\Phi_{(B,B^+)}^{-1}\) we can see that 
\[\Phi\colon \phi^*\MCE|_{\Ainf(R^+)}[1/d]\longrightarrow \MCE|_{\Ainf(R^+)}[1/d]\]
is an isomorphism, so it indeed defines the structure of a perfect-prismatic \(F\)-crystal.
\end{proof}

\begin{mycor}\label{BKFEquivalence}
    Let \(R^+\) be an integral perfectoid ring. Then there is an exact tensor equivalence:
    \[\emph{BKF-modules on }R^+\longrightarrow \emph{shtukas with fixed leg on }\Spd(R^+)\]
\end{mycor}
\subsection{Putting both together}
Let \(\MCG\) be an affine flat group scheme of finite type over \(\ZZ_p\), as fixed at the beginning of the section.

\begin{mydef}
    A \emph{perfect-prismatic \(\MCG\)-crystal} on an \(I\)-complete ring \(R^+\) is one of the following equivalent pieces of data:
    \begin{enumerate}
        \item {\emergencystretch=1em An exact tensor functor from \(\Rep(\MCG)\) to the category of perfect-prismatic \(F\)-crystals on \(R^+\).\par}
        \item For any map \((R^+\to S^+)\in (R^+)_{\Prism}^{\perf}\), an exact tensor functor from \(\Rep(\MCG)\) to the category of BKF-modules on \(S^+\), compatible with pullbacks.
    \end{enumerate}
\end{mydef}

These are easily seen to be equivalent as the exact structure on perfect-prismatic \(F\)-crystals is the same as that of the underlying \(\OpP\)-bundle. The same argument shows that the following two definitions are equivalent:

\begin{mydef}
    A \emph{(family of) \(\MCG\)-shtuka(s) with fixed leg} on \(X_\diamond\) is one of the following equivalent pieces of data:
    \begin{enumerate}
        \item An exact tensor functor from \(\Rep(\MCG)\) to the category of shtukas with fixed leg on \(X_\diamond\).
        \item For any map from a representable \(\Spa(B,B^+)\to X_\diamond\) with untilt \(\Spa(B^\sharp,B^{+\sharp})\), an exact tensor functor from \(\Rep(\MCG)\) to the category of shtukas with fixed leg on \(\Spa(B^\sharp,B^{+\sharp})\), compatible with pullbacks.
    \end{enumerate}
\end{mydef}

The following now follows immediately from the first characterization in both definitions above.

\begin{mycor}\label{newFGEquivalenceThm}
       Let \(R^+\) be as in \Cref{FEquivalenceThm}. The functor in that theorem induces an equivalence of groupoids
            \[\emph{perfect-prismatic \(\MCG\)-crystals on }R^+
            \longrightarrow\emph{\(\MCG\)-shtukas with fixed leg on }\Spd(R^+).\]
\end{mycor}
\newpage

\par\bigskip\noindent\textbf{Acknowledgements.}\enspace
First and foremost, I want to thank my Ph.D. advisor Torsten Wedhorn, who has taught and mentored me since my first contact with algebraic geometry seven years ago and who encouraged me to study prismatic \(F\)-crystals even before there was an established notion of such things in the literature. \par
I also want to sincerely thank the anonymous referee, who found essential gaps in the first version of the paper and whose thorough reading and useful comments contributed substantially to the final form of the article.\par
Furthermore, I want to thank Christopher Lang and Ian Gleason for helpful conversations and useful corrections.\par

\par\bigskip\noindent\textbf{Funding.}\enspace
This project was funded by the Deutsche Forschungsgemeinschaft (DFG, German Research Foundation) TRR 326 Geometry and Arithmetic of Uniformized Structures, project number 444845124.

\bibliographystyle{plainurl}
\begingroup
\small
\bibliography{bibliography}
\endgroup

\end{document}